\documentclass[12pt]{article}
\usepackage{amsfonts}

\usepackage{amsmath,amssymb,amsthm,latexsym,graphicx,citesort,color}

\newtheorem{theorem}{Theorem}
\newtheorem{conjecture}[theorem]{Conjecture}

\newtheorem{proposition}[theorem]{Proposition}
\newtheorem{definition}[theorem]{Definition}
\newtheorem{corollary}[theorem]{Corollary}
\newtheorem{remark}[theorem]{Remark}
\newcommand{\RR}{\mathbb{R}}

\newcommand{\R}{\mathbb{R}}

\newcommand{\W}{\mathcal{W}}
\newcommand{\M}{\mathcal{M}}

\newcommand{\mS}{\mathcal{S}}

\begin{document}

\author{P. Kuchment and L. Kunyansky}
\title{Mathematics of Photoacoustic and Thermoacoustic Tomography}
\date{}
\maketitle

\begin{center}{\large Contact information:\\
Peter Kuchment\\
Mathematics Department\\
Texas A \& M University\\
College Station, TX USA 77843-3368\\
phone: (979) 862-3257, FAX: (979) 862-4190\\
e-mail: kuchment@math.tamu.edu\\
Web: http://www.math.tamu.edu/\~{}kuchment
\\

Leonid Kunyansky\\
Department of Mathematics\\
University of Arizona, AZ USA  85721\\
phone: (520)621-4509, FAX:  (520)621-8322\\
e-mail: leonk@math.arizona.edu}
\end{center}

\newpage

\begin{abstract}
This is the manuscript of the chapter for a planned Handbook of Mathematical Methods in Imaging that surveys the mathematical models, problems, and algorithms of the
Thermoacoustic (TAT) and Photoacoustic (PAT) Tomography.
TAT and PAT represent probably the most developed of the several novel ``hybrid'' methods of medical imaging.
These new modalities combine different physical types of
waves (electromagnetic and acoustic in case of TAT and PAT) in such a way that the resolution and contrast of the resulting method are much higher than those achievable using only acoustic or electromagnetic measurements.
\end{abstract}
\tableofcontents

\section{Introduction}
\begin{figure}[ht!]
\begin{center}%
\includegraphics[width=2.5in,height=2.1in]{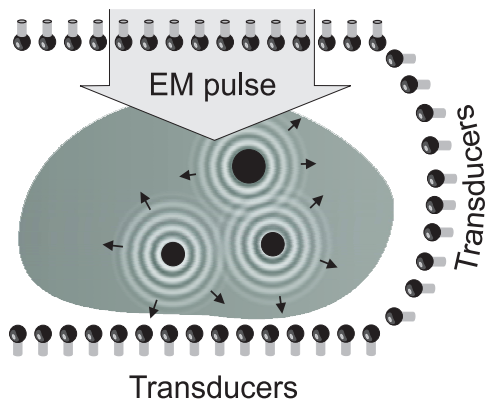}
\caption{TAT/PAT procedure with a partially surrounding acquisition surface.}
\end{center}\label{F:tat}
\end{figure}

We provide here just a very brief description of the TAT/PAT procedure, since
the relevant physics and biology details can be found in another chapter
\cite{Anast_hand} in this volume, as well as in the surveys and books
\cite{MXW_review,Wang_book,CRC}. In TAT (PAT), a short pulse of
radio-frequency EM wave (correspondingly, laser beam) irradiates a biological
object (e.g., in the most common application, human breast), thus causing
small levels of heating. The resulting thermoelastic expansion generates
a pressure wave that starts propagating through the object.
The absorbed EM energy and the initial pressure it creates are much higher in
the cancerous cells than in healthy tissues (see the
discussion of this effect in \cite{Anast_hand,MXW_review,Wang_book,CRC}).
Thus, if one could  reconstruct the initial pressure $f(x)$, the resulting
TAT tomogram would contain highly useful diagnostic information.
    The data for such a reconstruction are obtained by measuring
time-dependent pressure $p(x,t)$ using acoustic transducers located on a
surface $S$ (we will call it the {  observation} or
{  acquisition surface}) completely or partially surrounding the body (see
Fig. \ref{F:tat}).
Thus, although the initial irradiation is electro-magnetic, the actual
reconstruction is based on acoustic measurements. As a result, the high
contrast is produced due to a much higher absorption of EM energy by cancerous
cells (ultrasound alone would not produce good contrast in this case),
while the good (sub-millimeter) resolution is achieved by using ultrasound
measurements (the radio frequency EM waves are too long for high-resolution
imaging). Thus, TAT, by using two types of waves, combines their advantages,
while eliminating their individual deficiencies.

The physical principle upon which TAT/PAT is based was discovered by
Alexander Graham Bell in 1880 \cite{Bell} and its application
for imaging of biological tissues was suggested a century later \cite{Bowen}.
It began to be developed
as a viable medical imaging technique in the middle of 1990s
\cite{Oraev94,Kruger}.

Some of the mathematical foundations of this imaging modality were originally
developed starting in the
1990s for the purposes of the approximation theory \cite{LP1,LP2} (see
\cite{AQ,KuKuTAT} for extensive reviews of the resulting developments),
integral geometry (\cite[Chapter 5]{GGG}, \cite{Gi}),
and sonar and radar \cite{Cheney,LQ,NC}.

One can find recent reviews of the physics, biology, and mathematics issues
of TAT/PAT in
\cite{KuKuTAT,AKK,FR2,FR3,IP,Oraev,Oraev2,Pal_book,PatchSch,BiomDiagn,BiomPhot,CRC,MXW_review,Wang_book}.

TAT/PAT is just one, probably the most advanced at the moment, example of the
several recently introduced hybrid imaging methods, which combine different types of radiation to yield high quality
of imaging unobtainable by single-radiation modalities
(e.g., see \cite{Ammari_book,AmmariEIT,ScherzEIT,Wang_book,KuKuSynth} for other
examples).

\section{Mathematical models of TAT}\label{S:models}

In this section, we describe the commonly accepted mathematical model of the
TAT procedure and the main mathematical problems that need to be addressed.
Since for all our purposes PAT results in the same mathematical model
(although the biological features that TAT and PAT detect are different; see
details in the chapter  \cite{Anastasio}), we will refer to TAT only.

\subsection{Point detectors and the wave equation model}\label{S:wave}

We will mainly assume that point-like omni-directional ultrasound transducers,
located throughout an observation (acquisition) surface $S$, are used to
detect the values of the pressure $p(y,t)$, where $y\in S$ is a detector
location and $t\geq 0$ is the time of the observation. We also denote by
$c(x)$ the speed of sound at a location $x$. Then, it has been argued, that
the following model describes correctly the propagating pressure wave $p(x,t)$
generated during the TAT procedure (e.g.,
\cite{Diebold,Tam,MXW1,Anastasio,Anast_hand}):
\begin{equation}\label{E:wave}
\begin{cases}
    p_{tt}=c^2(x)\Delta_x p, \quad t\geq 0, x\in\RR^3\\
    p(x,0)=f(x), p_t(x,0)=0.
    \end{cases}
\end{equation}
Here $f(x)$ is the initial value of the acoustic pressure, which one needs
to find in order to create the TAT image.
In the case of a closed acquisition surface $S$, we will denote by $\Omega$
the interior domain it bounds. Notice that in TAT the function $f(x)$ is
naturally supported inside $\Omega$. We will see that this assumption about
the support of $f$ sometimes becomes crucial for the feasibility of
reconstruction, although some issues can be resolved even if $f$ has non-zero
parts outside the acquisition surface.

\begin{figure}[ht!]
\begin{center}%
\includegraphics[width=2.4in,height=1.9in]{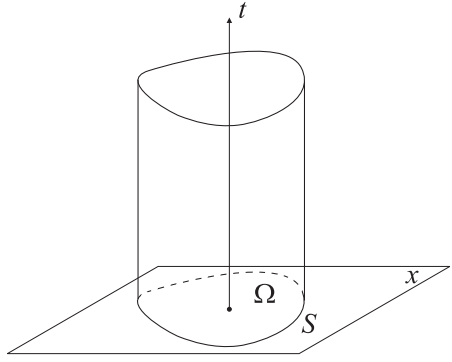}
\end{center}
\caption{The observation surface $S$ and the domain $\Omega$ containing the
object to be imaged.}
\label{F:cylinder}
\end{figure}

The data obtained by the point detectors located on a surface $S$ are
represented by the function
\begin{equation}\label{E:data}
    g(y,t):=p(y,t) \mbox{ for }y\in S, t\geq 0.
\end{equation}
Fig. \ref{F:cylinder} illustrates the space-time geometry of (\ref{E:wave}).

We will incorporate the measured data $g$ into the system (\ref{E:wave}), rewriting it as
follows:
\begin{equation}\label{E:wave_data}
\begin{cases}
    p_{tt}=c^2(x)\Delta_x p, \quad t\geq 0, x\in\RR^3\\
    p(x,0)=f(x), p_t(x,0)=0\\
    p|_{S}=g(y,t), \quad (y,t)\in S\times\RR^+.
    \end{cases}
\end{equation}

Thus, the goal in TAT/PAT is to
{  find, using the data $g(y,t)$ measured by transducers, the initial value
$f(x)$ at $t=0$ of the solution $p(x,t)$ of (\ref{E:wave_data})}.

We will use the following notation:
\begin{definition}\label{D:W}
We will denote by $\W$ the {  forward operator}
\begin{equation}\label{E:forward}
\W: f(x)\mapsto g(y,t),
\end{equation}
where $f$ and $g$ are described in (\ref{E:wave_data}).
\end{definition}

\begin{remark}\indent
\begin{itemize}

  \item The reader should notice that if a different type of detectors is
  used, the system (\ref{E:wave}) stays intact, while the measured data will
  be represented differently from (\ref{E:data}) (see Section
  \ref{S:integrating}). This will correspondingly influence the reconstruction
  procedures.

  \item We can consider the same problem in the space $\R^n$ of any dimension,
  not just in $3D$. This is not merely a mathematical abstraction.  Indeed, in
  the case of the so called integrating line detectors (Section
  \ref{S:integrating}), one deals with the $2D$ situation.
\end{itemize}

\end{remark}

\subsection{Acoustically homogeneous media and spherical means}\label{S:spherical}

If the medium being imaged is acoustically homogeneous (i.e., $c(x)$ equals to
a constant, which we will assume to be equal to $1$ in appropriate units), as
it is approximately the case in breast imaging, one deals with the constant
coefficient wave equation problem
\begin{equation}\label{E:wave_const}
\begin{cases}
    p_{tt}=\Delta_x p, \qquad t\geq 0, x\in\RR^3\\
    p(x,0)=f(x), p_t(x,0)=0\\
    p|_{S}=g(y,t), \quad (y,t)\in S\times\RR^+.
    \end{cases}
\end{equation}
In this case, the well known Poisson-Kirchhoff formulas
\cite[Ch. VI, Section 13.2, Formula (15)]{CH} for the solution of the wave
equation gives in $3D$:
\begin{equation}\label{E:KP}
p(x,t)=a\frac{\partial}{\partial t}\left(t(Rf)(x,t)\right),
\end{equation}
where
\begin{equation}\label{E:mean}
(Rf)(x,r):= \frac{1}{4 \pi}
\int\limits_{|y|=1} f(x+ry)dA(y)
\end{equation}
is the spherical mean operator applied to the function $f(x)$, $dA$ is
the standard area element on the unit sphere in $\RR^3$, and $a$ is a
constant. (Versions in all dimensions are known, see
(\ref{E:KirchPois_even}) and (\ref{E:KirchPois_odd}).) One can derive from
here that knowledge of the function $g(x,t)$ for
$x\in S$ and all $t\geq 0$ is equivalent to knowing the spherical mean $Rf
(x,t)$ of the function $f$ for any points $x\in S$ and any $t\geq 0$. One thus
needs to study the spherical mean operator $R:f\to Rf$, or, more precisely, its
restriction to the points $x\in S$ only, which we will denote by $\M$:
\begin{equation}\label{E:Radon_S}
\M f(x,t):=\frac{1}{4 \pi}
                \int\limits_{|y|=1}f(x+ty)dA(y), \quad x\in S, t\geq 0.
\end{equation}
Due to the connection between the spherical
mean operator and the wave equation, one can choose to work with the former,
and in fact many works on TAT do so. The spherical mean operator $\M$ resembles
the classical Radon transform, the common tool of computed tomography
\cite{Natt_old,Natt_new,Kak}, which integrates functions over planes rather
than spheres. This analogy with Radon transform, although often purely
ideological, rather than technical, provides important intuition and
frequently points in reasonable directions of study. However, when the medium
cannot be assumed to be acoustically homogeneous, and thus $c(x)$ is not
constant, the relation between TAT and integral geometric transforms,
such as Radon transform or spherical mean, to a large extent breaks down,
and thus one has to work with the wave equation directly.

In what follows, we will address both models of TAT (the PDE model and the
integral geometry model) and thus will deal with both forward operators $\W$
and $\M$.

\subsection{Main mathematical problems arising in TAT}\label{S:problems}

We now formulate a list of problems related to TAT which will be
addressed in detail in the rest of the article. (This list is more or
less standard for a tomographic imaging method.)

\begin{description}
\item{\textbf{Sufficiency of the data.}}
The first natural question to ask is: Is the data collected on the observation
surface $S$ sufficient for the unique reconstruction of the initial pressure
$f(x)$ (\ref{E:wave_data})? In other words, is the kernel of the forward
operator $\W$ zero? Or, to put it differently, for which sets $S\in \RR^3$ the
data collected by transducers placed along $S$ determines $f$ uniquely? Yet
another interpretation of this question is through observability of solutions
of the wave equation on the set $S$: does observation on $S$ of a solution of
the problem (\ref{E:wave}) determine the solution uniquely?

    When the speed of sound is constant, and thus the spherical mean model
    applies, the equivalent question is whether the operator $\M$ has zero kernel on an
    appropriate class of functions (say, continuous functions with compact
    support)

    As it is explained in \cite{AQ}, the choice of precise conditions on the
    local function class, such as continuity, is of no importance for the
    answer to the uniqueness question, while  behavior at infinity (e.g.,
    compactness of support) is. So, without loss of generality, when
    discussing uniqueness, one can assume $f(x)$ in (\ref{E:wave_data}) to be
    infinitely differentiable.

\item{\textbf{Inversion formulas and algorithms.}} Since a practitioner needs
to see the actual tomogram, rather than just know its existence, the next
natural question arises: If uniqueness the data collected on $S$ is
established, what are the actual inversion formulas or algorithms? Here again
one can work with smooth functions, in the end extending the formulas by
continuity to a wider class.

\item{\textbf{Stability of reconstruction.}} If we can invert the transform
and reconstruct $f$ from the data $g$, how stable is the inversion? The
measured data are unavoidably corrupted by errors, and stability means that
small errors in the data lead to only small errors in the reconstructed
tomogram.

\item{\textbf{Incomplete data problems.}} What happens if the data is
``incomplete,'' for instance if one can only partially surround the object by
transducers? Does this lead to any specific deterioration in the tomogram, and
if yes, to what kind of deterioration?

\item{\textbf{Range descriptions.}}

The next question is known to be important in analysis of tomographic problems:
What is the range of
the forward operator $\W:f\mapsto g$ that maps the unknown function $f$ to the measured
data $g$? In other words, what is the space of all possible ``ideal'' data
$g(t,y)$ collected on the surface $S$? In the constant speed of sound case,
this is equivalent to the question of describing the range of the spherical
mean operator $\M$ in appropriate function spaces. Such ranges often have
infinite co-dimensions, and the importance of knowing the range of Radon type
transforms for analyzing problems of tomography is well known. For instance,
such information is used to improve inversion algorithms, complete incomplete
data, discover and compensate for certain data errors, etc. (e.g.,
\cite{Leon_Radon,GGG1,GGG,GelfVil,Helg_Radon,Helg_groups,He1,Natt_old,Natt_new,Pal_book,Kuch_AMS05} and references therein).
In TAT, range descriptions are also closely connected with the speed of sound
determination problem listed next (see Section \ref{S:speed} for a discussion
of this connection).

\item{\textbf{Speed of sound reconstruction.}} As the reader can expect, reconstruction procedures require the knowledge of the speed of sound $c(x)$. Thus, the problem arises of the recovery of $c(x)$ either from an additional scan, or (preferably) from
 the same TAT data.

\end{description}

\subsection{Variations on the theme: planar, linear, and circular integrating detectors}\label{S:integrating}

In the described above most basic and well-studied version of TAT, one utilizes point-like broadband transducers to measure the acoustic wave on a surface surrounding the object of interest. The corresponding
mathematical model is described by the system
(\ref{E:wave_data}). In practice, the transducers
cannot be made small enough, since smaller detectors yield weaker signals resulting
in low signal-to-noise ratios. Smaller transducers are also more
difficult to manufacture.

Since finite size of the transducers limits the resolution of the
reconstructed images, researchers have been trying to design alternative
acquisition schemes using receivers that are very thin but long or wide.
Such are $2D$ planar detectors \cite{haltmaier_large,haltmaier_int_and line}
and $1D$ linear and circular \cite{haltmaier_interf,haltmaier_fabri,Grun,Zangerl}
detectors.

We will assume throughout this section that the speed of sound
$c(x)$ is constant and equal to 1.

Planar detectors are made from a thin piezoelectric polymer
film glued onto a flat substrate (see, for example
\cite{Paltauf_review}). Let us assume that the object is
contained within the sphere of radius $R$. If the diameter of
the planar detector is sufficiently large (see
\cite{Paltauf_review} for details), it can be assumed to be
infinite. The mathematical model of such an acquisition
technique is no longer described by (\ref{E:wave_data}). Let us
define the detector plane $\Pi(s,\omega)$ by equation
$x\cdot\omega=s$, where $\omega$ is the unit normal to the
plane and $s$ is the (signed) distance from the origin to the
plane. Then, while the propagation of acoustic waves is still
modeled by (\ref{E:wave}), the measured data
$g_{planar}(s,t,\omega)$ (up to a constant factor which
we will, for simplicity, assume to be equal to 1) can be
represented by the following integral:
\[
g_{planar}(s,\omega,t)=\int\limits_{\Pi(s,\omega)}p(x,t)dA(x)
\]
where $dA(x)$ is the surface measure on the plane. Obviously,
\[
g_{planar}(s,\omega,0)=\int\limits_{\Pi(s,\omega)}p(x,0)dA(x)=\int
\limits_{\Pi(s,\omega)}f(x)dA(x)\equiv F(s,\omega),
\]
i.e. the value of $g$ at $t=0$ coincides with the integral
$F(s,\omega)$ of the initial pressure $f(x)$ over the plane $\Pi(s,\omega)$ orthogonal to $\omega$.

One can show \cite{haltmaier_large,haltmaier_int_and line} that
for a fixed $\omega$, function $g_{planar}(s,\omega,t)$ is the
solution to $1D$ wave equation
\[
\frac{\partial^{2}g}{\partial
s^{2}}=\frac{\partial^{2}g}{\partial t^{2}},
\]
and thus
\begin{align*}
g_{planar}(s,\omega,t)  & =\frac{1}{2}\left[  g_{planar}(s,\omega
,s-t)+g_{planar}(s,\omega,s+t)\right]  \\
& =\frac{1}{2}\left[  F(s+t,\omega)+F(s-t,\omega)\right]  .
\end{align*}
Since the detector can only be placed outside the object, i.e.
$s\geq R$, the term $F(s+t,\omega)$ vanishes, and one obtains
\[
g_{planar}(s,\omega,t)=F(s-t,\omega).
\]
In other words, by measuring $g_{planar}(s,\omega,t)$, one can
obtain values of the planar integrals of $f(x)$. If, as
proposed in \cite{haltmaier_large,haltmaier_int_and line}, one
conducts measurements for all planes tangent to the upper
half-sphere of radius $R$ (i.e. $s=R,\omega\in S_{+}^{2})$,
then the resulting data yield all values of the standard Radon
transform of $f(x)$. Now the reconstruction can be carried out
using one of the many known inversion algorithms for the latter
transform (see \cite{Natt_old,Natt_new,Kak}).

Linear detectors are based on optical detection of acoustic
signal. Some of the proposed optical detection schemes utilize
as the sensitive element a thin straight optical fiber in
combination with Fabry-Perot interferometer
\cite{haltmaier_fabri,Grun}. Changes of acoustic pressure on
the fiber change (proportionally) its length; this elongation,
in turn, is detected by interferometer. A similar idea is used
in \cite{haltmaier_interf}; in this work the role of a
sensitive element is played by a laser beam passing through the
water in which the object of interest is submerged, and thus
the measurement does not perturb the acoustic wave. In both
cases, the length of the sensitive element exceeds the size of
the object, while the diameter of the fiber (or of the laser
beam) can be made extremely small (see \cite{Paltauf_review}
for a detailed discussion), which removes restrictions on
resolution one can achieve in the images.

Let us assume that the fiber (or laser beam) is aligned along
the line
$l(s_{1},s_{2},\omega_{1},\omega_{2})=\{x|x=s_{1}\omega_{1}+s_{2}\omega
_{2}+s\omega\}$, where vectors $\omega_{1}, \omega_{2}$, and
$\omega$ form an ortho-normal basis in $\mathbb{R}^{3}$. Then
the measured quantities $g_{linear}
(s_{1},s_{2},\omega_{1},\omega_{2},t)$ are equal (up to a
constant factor which, we will assume, equals to 1) to the
following line integral:
\[
g_{linear}(s_{1},s_{2},\omega_{1},\omega_{2},t)=\int\limits_{\mathbb{R}^{1}
}p(s_{1}\omega_{1}+s_{2}\omega_{2}+s\omega,t)ds.
\]
Similarly to the case of planar detection, one can show
\cite{haltmaier_fabri,Grun,haltmaier_interf}, that for fixed
vectors $\omega_{1},\omega_{2}$ the measurements
$g_{linear}(s_{1},s_{2},\omega _{1},\omega_{2},t)$ satisfy the
$2D$ wave equation
\[
\frac{\partial^{2}g}{\partial s_{1}^{2}}+\frac{\partial^{2}g}{\partial
s_{2}^{2}}=\frac{\partial^{2}g}{\partial t^{2}}.
\]
The initial values $g_{linear}(s_{1},s_{2},\omega_{1},\omega
_{2},0)$ coincide with the line integrals of $f(x)$ along lines
$l(s_{1},s_{2},\omega_{1},\omega_{2})$. Suppose one makes
measurements for all values of $s_{1}(\tau),$ $s_{2}(\tau)$
corresponding to a curve
$\gamma=\{x|x=s_{1}(\tau)\omega_{1}+s_{2}(\tau)\omega_{2},\tau_{0}\leq\tau
\leq\tau_{1}\}$ lying in the plane spanned by
$\omega_{1},\omega_{2}$. Then one can try to reconstruct the
initial value of $g$ from the values of $g$ on $\gamma$. This
problem is a $2D$ version of (\ref{E:wave_data}) and thus the
known algorithms (see Section \ref{S:reconstr}) are applicable.

In order to complete the reconstruction from data obtained
using line detectors, the measurements should be repeated with
different directions of $\omega$. For each value of $\omega$
the $2D$ problem is solved; the solutions of these problems
yield values of line integrals of $f(x)$. If this is done for
all values of $\omega$ lying on a half circle, the set of the
recovered line integrals of $f(x)$ is sufficient for
reconstructing this function. Such a reconstruction represents
the inversion of the well known in tomography X-ray transform.
The corresponding theory and algorithms can be found, for
instance, in \cite{Natt_old,Natt_new,Kak}.

Finally, the use of circular integrating detectors was
considered in \cite{Zangerl}. Such a detector can be made out of optical fiber combined with an interferometer. In
\cite{Zangerl}, a closed form solution of the corresponding
inverse problem is found. However, this approach is very new
and neither numerical examples, nor reconstructions from real
data have been obtained yet.

\section{Mathematical analysis of the problem}\label{S:analysis}
In this section, we will address most of the issues described in Section
\ref{S:problems}, except the reconstruction algorithms, which will be
discussed in Section \ref{S:reconstr}.

\subsection{Uniqueness of reconstruction}\label{S:unique}
The problem discussed here is the most basic one for tomography: given an
acquisition surface $S$ along which we distribute detectors, is the data
$g(y,t)$ for $y\in S, t\geq 0$ (see (\ref{E:wave_data})) sufficient for a
unique reconstruction of the tomogram $f$? A simple counting of variables
shows that $S$ should be a hyper-surface in the ambient space (i.e., a surface
in $\R^3$ or a curve in $\R^2$). As we will see below, although there are some
simple counter-examples and remaining open problems, for all practical
purposes, the uniqueness problem is positively resolved, and most surfaces $S$
do provide uniqueness. We address this issue for acoustically homogeneous
media first and then switch to the variable speed case.

Before doing so, however, we would like to dispel a  concern that arises when
one looks at the problem of recovering $f$ from $g$ in (\ref{E:wave_data}).
Namely, an impression might emerge that we consider an initial-boundary value
(IBV) problem for the wave equation in the cylinder $\Omega\times\R^+$, and
the goal is to recover the initial data $f$ from the known boundary data $g$.
This is clearly impossible, since according to standard PDE theorems (e.g.,
\cite{Evans,CH}), one can solve this IBV problem for \textbf{arbitrary} choice
of the initial data $f$ and boundary data $g$ (as long as they satisfy simple
compatibility conditions, which are fulfilled for instance if $f$ vanishes
near $S$ and $g$ vanishes for small $t$, which is the case in TAT). This means
that apparently $g$ contains essentially no information about $f$ at all. This
argument, however, is flawed, since the wave equation in (\ref{E:wave_data})
holds in the whole space, not just in $\Omega$. In other words, $S$ is not a
boundary, but rather an observation surface. In particular, considering the
wave equation in the exterior of $S$, one can derive that if $f$ is supported
inside $\Omega$, the boundary values $g$ of the solution $p$ of
(\ref{E:wave_data}) also determine the normal derivative of $p$ at $S$ for all
positive times. Thus, we in fact have (at least theoretically) the full Cauchy
data of the solution $p$ on $S$, which should be sufficient for
reconstruction. Another way of addressing this issue is to notice that if the
speed of sound is constant, or at least non-trapping (see the definition below
in Section \ref{S:unique_var}), the energy of the solution in any bounded
domain (in particular, in $\Omega$) must decay in time. The decay when
$t\to\infty$ together with the boundary data $g$ guarantee the uniqueness of
solution, and thus uniqueness of recovery $f$.

These arguments, as the reader will see, play a role in understanding
reconstruction procedures.

\subsubsection{Acoustically homogeneous media}\label{S:unique_const}

We assume here the sound speed $c(x)$ to be constant (in appropriate units,
one can choose it to be equal to $1$, which we will do to simplify
considerations).

In order to state the first important result on uniqueness, let us recall the system (\ref{E:wave_const}), allowing an arbitrary dimension $n$ of
the space:
\begin{equation}\label{E:wave_const_n}
\begin{cases}
    p_{tt}=\Delta_x p, \qquad t\geq 0, x\in\RR^n\\
    p(x,0)=f(x), p_t(x,0)=0\\
    p|_{S}=g(y,t), \quad (y,t)\in S\times\RR^+.
    \end{cases}
\end{equation}

We introduce the following useful definition:
\begin{definition}\label{D:unique}
A set $S$ is said to be {  uniqueness set}, if when used as the acquisition
surface, it provides sufficient data for unique reconstruction of the
compactly supported tomogram $f$ (i.e., the observed data $g$ in
(\ref{E:wave_const_n}) determines uniquely function $f$).
Otherwise, it is called a {  non-uniqueness set}.
\end{definition}
\noindent In other words, $S$ is a uniqueness set if the forward operator $\W$ (or, equivalently, $\M$)
has zero kernel.

We will start with a very general statement about the acquisition
(observation) sets $S$ that provide insufficient information for unique
reconstruction of $f$  (see \cite{AQ} for the proof and references):
\begin{theorem}\label{T:harm}
If $S$ is a non-uniqueness set, then there exists a non-zero harmonic
polynomial $Q$, which vanishes on $S$.
\end{theorem}
This theorem implies, in particular, that all ``bad'' (non-uniqueness)
observation sets are algebraic, i.e. have a polynomial vanishing on them.
Turning this statement around, we conclude that any set $S$ that is a
uniqueness set for harmonic polynomials, is sufficient for unique TAT
reconstruction (although, as we will see in Section \ref{S:incomplete}, this
does not mean practicality of the reconstruction).

The proof of Theorem \ref{T:harm}, which the reader can find in
\cite{AQ,KuKuTAT}, is not hard and in fact is enlightening, but providing it
would lead us too far from the topic of this survey.

We will consider first the case of closed acquisition surfaces, i.e. the ones
that completely surround the object to be imaged. We will address the
general situation afterwards.

\paragraph{Closed acquisition surfaces $S$}\indent

\begin{theorem}\label{T:uniq_closed}(\cite{AQ})
If the acquisition surface $S$ is the boundary of bounded domain $\Omega$
(i.e., a closed surface), then it is a uniqueness set. Thus, the observed data
$g$ in (\ref{E:wave_const_n}) determines uniquely the sought function
$f\in L^2_{comp}(\R^n)$. (The statement holds, even though $f$ is not required
to be supported inside $S$.)
\end{theorem}

{  Proof:} Indeed, since there are no non-zero harmonic functions vanishing
on a closed surface $S$, Theorem \ref{T:harm} implies Theorem
\ref{T:uniq_closed}.\qed

There is, however, another, more intuitive, explanation of why Theorem
\ref{T:uniq_closed} holds true (although it requires somewhat stronger
assumptions, or a more delicate proof than the one indicated below). Namely,
since the solution $p$ of (\ref{E:wave_const_n}) has compactly supported
initial data, its energy is decaying inside any bounded domain, in particular inside
$\Omega$ (see Section \ref{S:unique_var} and \cite{Egorov,Hristova} and
references therein about local energy decay). On the other hand, if there is
non-uniqueness, there exists a non-zero $f$ such that $g(y,t)=0$ for all $y\in
S$ and $t$. This means that we can add homogeneous Dirichlet boundary
conditions $p\mid_S=0$ to (\ref{E:wave_const_n}). But then the standard PDE
theorems \cite{CH,Evans} imply that the energy stays constant in $\Omega$.
Combination of the two conclusions means that $p$ is zero in $\Omega$ for all
times $t$. It is well known \cite{CH} that such a solution of the wave
equation must be identically zero everywhere, and thus $f=0$.

This energy decay consideration can be extended to some classes of non-compactly supported functions $f$ of the $L^p$ classes, leading to the following result of \cite{ABK}:

\begin{theorem}\label{T:ABK}\cite{ABK}
Let $S$ be the boundary of a bounded domain in $\R^n$ and $f\in L^p(\R^n)$.
Then
\begin{enumerate}
\item If $p\leq\frac{2n}{n-1}$ and the spherical mean of $f$ over almost
every sphere centered on $S$ is equal to zero, then $f=0$.

\item The previous statement fails when $p>\frac{2n}{n-1}$ and $S$ is a
sphere.

\end{enumerate}
In other words, a closed surface $S$ is a uniqueness set for functions $f\in
L^p(\R^n)$ when $p\leq\frac{2n}{n-1}$, and might fail to be such when
$p>\frac{2n}{n-1}$.
\end{theorem}
This result shows that the assumption, if not necessarily of compactness of support of $f$, but at
least of a sufficiently fast decay of $f$ at infinity, is important for the uniqueness to hold.

\paragraph{General acquisition sets $S$}\indent

Theorems \ref{T:harm} and \ref{T:uniq_closed} imply the following useful
statement:
\begin{theorem}\label{T:uniq-simple}
If a set $S$ is not algebraic, or if it contains an open part of a closed
analytic surface $\Gamma$, then it is a uniqueness set.
\end{theorem}
Indeed, the first claim follows immediately from Theorem  \ref{T:harm}. The
second one works out as follows: if an open subset of an analytic surface $\Gamma$ is a non-uniqueness set, then by an analytic continuation
type argument (see \cite{AQ}), one can show that the whole $\Gamma$ is such a set.
However, this is impossible, due to Theorem  \ref{T:uniq_closed}.

There are simple examples of non-uniqueness surfaces. Indeed, if $S$ is a
plane in $3D$ (or a line in $2D$, or a hyperplane in dimension $n$) and $f(x)$
in (\ref{E:wave_data}) is odd with respect to $S$, then clearly the whole
solution of (\ref{E:wave_data}) has the same parity and thus vanishes on $S$
for all times $t$. This means that, if one places transducers on a planar $S$,
they might register zero signals at all times, while the function $f$ to be
reconstructed is not zero. Thus, there is no uniqueness of reconstruction when
$S$ is a plane. On the other hand (see \cite{John,CH}), if $f$ is supported
completely on one side of the plane $S$ (the standard situation in TAT), it is uniquely
recoverable from its spherical means centered on $S$, and thus from the observed data $g$.

\begin{figure}[ht!]
\begin{center}%
\includegraphics[width=2.0in,height=1.7in]{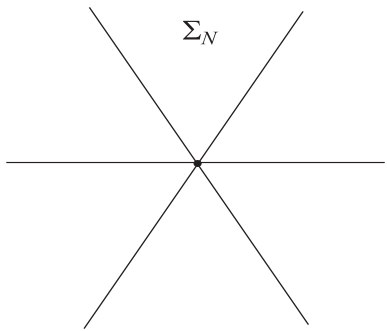}
\end{center}
\caption{Coxeter cross of $N$ lines.}\label{F:coxeter}
\end{figure}

The question arises what are other ``bad'' (non-uniqueness) acquisition surfaces than planes. This issue has been resolved in $2D$ only.
Namely, consider a set of $N$ lines on the plane intersecting at a point and
forming at this point equal angles. We will call such a figure the {
Coxeter cross $\Sigma_N$} (see Fig. \ref{F:coxeter}).
it is easy to construct a compactly supported function that is odd
simultaneously with respect of all lines in $\Sigma_N$. Thus, a Coxeter cross
is also a non-uniqueness set. The following result, conjectured in
\cite{LP1,LP2} and proven in the full generality in \cite{AQ}, shows that, up
to adding finitely many points, this is all that can happen to non-uniqueness
sets:

\begin{theorem}\cite{AQ}\label{T:AQ}
A set $S$ in the plane $\R^2$ is a non-uniqueness set for compactly supported
functions $f$, if and only if it belongs to the union $\Sigma_N\bigcup\Phi$ of
a Coxeter cross $\Sigma_N$ and a finite set of points $\Phi$.
\end{theorem}
Again, compactness of support is crucial for the proof provided in \cite{AQ}.
There are no other proofs known at the moment of this result (see the
corresponding open problem in Section \ref{S:remarks}). In particular, there is
no proven analog of Theorem \ref{T:ABK} for non-closed sets $S$ (unless $S$ is
an open part of a closed analytic surface).

The $n$-dimensional (in particular, $3D$) analog of Theorem \ref{T:AQ} has
been conjectured \cite{AQ}, but never proven, although some partial advances
in this direction have been made in \cite{AmbKuc_inj,FPR}.

\begin{figure}[ht!]
\begin{center}%
\includegraphics[width=1.5in,height=1.9in]{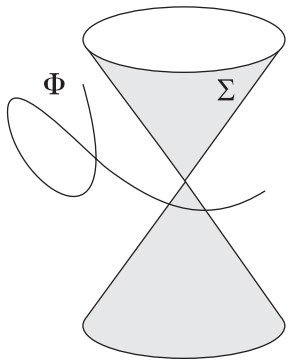}
\end{center}
\caption{The conjectured structure of a most general non-uniqueness set in
$3D$.}\label{F:n-dim}
\end{figure}

\begin{conjecture}\label{C:n-dim}
A set $S$ in $\R^n$ is a non-uniqueness set for compactly supported functions
$f$, if and only if it belongs to the union $\Sigma\bigcup\Phi$, where
$\Sigma$ is the cone of zeros of a homogeneous (with respect to some point in
$\R^n$) harmonic polynomial, and $\Phi$ is an algebraic sub-set of $\R^n$ of
dimension at most $n-2$ (see Fig. \ref{F:n-dim}).
\end{conjecture}

\paragraph{Uniqueness results for a finite observation time}\indent

So far, we have addressed only the question of uniqueness of reconstruction in
the non-practical case of the infinite observation time. There are, however,
results that guarantee uniqueness of reconstruction for a finite time of
observation. The general idea is that it is sufficient to observe for the time
that it takes the geometric rays (see
Section \ref{S:unique_var}) from the interior $\Omega$ of $S$ to reach $S$. In
the case of a constant speed, which we will assume to be equal to $1$, the
rays are straight and are traversed with the unit speed. This means that if
$D$ is the diameter of $\Omega$ (i.e., the maximal distance between two points
in the closure of $\Omega$), then after time $t=D$, all rays coming from
$\Omega$ have left the domain. Thus, one hopes that waiting till time
$t=D$ might be sufficient. In fact, due to the specific initial conditions in
(\ref{E:wave_data}), namely, that the time derivative of the pressure is equal
to zero at the initial moment, each singularity of $f$ emanates two rays, and
at least one of them will reach $S$ in time not exceeding $D/2$. And indeed,
the following result of \cite{FPR} holds:
\begin{theorem}\cite{FPR}\label{T:halftime_const}
If $S$ is smooth and closed surface bounding domain $\Omega$ and $D$ is the
diameter of $\Omega$, then the TAT data on $S$ collected for the time $0\leq
t\leq 0.5D$, uniquely determines $f$.
\end{theorem}
Notice that a shorter collection time does not guarantee uniqueness. Indeed,
if $S$ is a sphere and the observation time is less than $0.5D$, due to the
finite speed of propagation, no information from a neighborhood of the center
can reach $S$ during observation. Thus, values of $f$ in this neighborhood
cannot be reconstructed.

\subsubsection{Acoustically inhomogeneous media}\label{S:unique_var}

We assume that the speed of sound is strictly positive, $c(x)>c>0$, and such
that $c(x)-1$ has compact support, i.e. $c(x)=1$ for large $x$.

\paragraph{Trapping and non-trapping}\indent

We will frequently impose the so-called non-trapping condition on the speed of sound $c(x)$ in $\R^n$. To introduce it, let us consider the Hamiltonian system in $\RR^{2n}_{x,\xi}$ with
the Hamiltonian $H=\frac{c^2(x)}{2}|\xi|^2$:
\begin{equation}\label{E:bichar}
\begin{cases}
x^\prime_t=\frac{\partial H}{\partial \xi}=c^2(x)\xi \\
\xi^\prime_t=-\frac{\partial H}{\partial x}=-\frac 12 \nabla
\left(c^2(x)\right)|\xi|^2 \\
x|_{t=0}=x_0, \quad \xi|_{t=0}=\xi_0.
\end{cases}
\end{equation}
The solutions of this system are called {  bicharacteristics}
and their projections into $\RR^n_x$ are {  rays} (or {  geometric rays}).

\begin{definition}\label{D:nontrap}
We say that the speed of sound $c(x)$ satisfies the {  non-trapping
condition}, if all rays with $\xi_0\neq 0$ tend to infinity when $t \to
\infty$.

The rays that do not tend to infinity, are called {  trapped}.
\end{definition}

A simple example, where quite a few rays are trapped, is the radial parabolic sound speed $c(x)=c|x|^2$.

It is well known (e.g., \cite{Horm}) that singularities of solutions of the
wave equation are carried by geometric rays. In order to make this statement
more precise, we need to recall the notion of a
wave front set $WF(u)$ of a distribution $u(x)$ in $\R^n$. This set carries
detailed information on singularities of $u(x)$.

\begin{definition}\label{D:WF}
Distribution $u(x)$ is said to be \textbf{microlocally smooth near a point
$(x_0,\xi_0)$}, where $x_0,\xi_0\in\R^n$ and $\xi_0\neq 0$, if there is a
smooth ``cut-off'' function $\phi(x)$ such that $\phi(x_0)\neq0$ and that the
Fourier transform $\widehat{\phi u}(\xi)$ of the function $\phi(x)u(x)$ decays
faster than any power $|\xi|^{-N}$ when $|\xi|\to\infty$, in directions that
are close to the direction of $\xi_0$. \footnote{We remind the reader that if
this Fourier transform decays that way in {  all} directions, then $u(x)$ is
smooth (infinitely differentiable) near the point $x_0$.}

The \textbf{wave front set $WF(u)\subset \R^n_x\times(\R^n_\xi\setminus 0)$ of
$u$} consists of all pairs $(x_0,\xi_0)$ such that $u$ is \underline{not}
microlocally smooth near $(x_0,\xi_0)$.
\end{definition}
In other words, if $(x_0,\xi_0)\in WF(u)$, then $u$ is not smooth near $x_0$,
and the direction of $\xi_0$ indicates why it is not: the Fourier transform
does not decay well in this direction. For instance, if $u(x)$ consists of two
smooth pieces joined non-smoothly across a smooth interface $\Sigma$, then
$WF(u)$ can only contain pairs $(x,\xi)$ such that $x\in\Sigma$ and $\xi$ is
normal to $\Sigma$ at $x$.

It is known that the wave front sets of solutions of the wave equation
propagate with time along the bicharacteristics introduced above. This is a
particular instance of a more general fact that applies to general PDEs and
can be found in \cite{Shubin,Horm}. As a result, if after time $T$ all the
rays leave the domain $\Omega$ of interest, the solution becomes smooth
(infinitely differentiable) inside $\Omega$.

One can find simple introduction to the notions of microlocal analysis, such
as the wave front set, for instance in \cite{Str}, and more advanced versions
in \cite{Horm,Shubin}. Applications of microlocal analysis to integral
geometry are discussed in \cite{GreenUhlm,Guill75,Guill85,GS}.


The notion of so called local energy decay, which we survey next,
is important for the understanding of the non-trapping conditions in TAT.

\paragraph{Local energy decay estimates}\indent

Assuming that the initial data $f(x)$ (\ref{E:wave}) is compactly supported
and the speed $c(x)$ is non-trapping, one can provide the so called
\textbf{local energy decay estimates} \cite{Egorov,Vainb,Vainb2}. Namely, in
any bounded domain $\Omega$, the solution $p(x,t)$ of (\ref{E:wave})
satisfies, for a sufficiently large $T_0$ and for any $(k,m)$, the estimate
\begin{equation}\label{E:decay}
    \left|\frac{\partial^{k+|m|}}{\partial^k_t \partial_x^m}\right|\leq
C_{k,m}\nu_k(t)\|f\|_{L^2}, \mbox{ for }x\in\Omega,t>T_0.
\end{equation}
Here $\nu_k(t)=t^{-n+1-k}$ for even $n$ and $\nu_k(t)=e^{-\delta t}$ for odd
$n$ and some $\delta>0$. Any value $T_0$ larger than the diameter of $\Omega$
works in this estimate.

\paragraph{Uniqueness result for non-trapping speeds}\indent

If the speed is non-trapping, the local energy decay allows one to start
solving the problem (\ref{E:wave_data}) from $t=\infty$, imposing zero
conditions at $t=\infty$ and using the measured data $g$ as the boundary
conditions. This leads to recovery of the whole solution, and in particular
its initial value $f(x)$. As the result, one obtains the following simple
uniqueness result of \cite{AK}:
\begin{theorem}\label{T:uniq-variable}\cite{AK}
If the speed $c(x)$ is smooth and non-trapping and the acquisition surface $S$
is closed, then the TAT data $g(y,t)$ determines the tomogram $f(x)$ uniquely.
\end{theorem}
Notice that the statement of the theorem holds even if the support of $f$ is
not completely inside of the acquisition surface $S$.

\paragraph{Uniqueness results for finite observation times}\indent

As in the case of constant coefficients, if the speed of sound is
non-trapping, appropriately long finite observation time suffices for the
uniqueness. Let us denote by $T(\Omega)$ the {\em supremum}
 of the time it takes the ray to reach $S$, over all rays
originating in $\Omega$.
In particular, if $c(x)$ is trapping, $T(\Omega)$ might be infinite.

\begin{theorem}\cite{StefUhlTAT}\label{T:SU}
The data $g$ measured till any time $T$ larger than $T(\Omega)$ is sufficient
for unique recovery of $f$.
\end{theorem}

\subsection{Stability}\label{S:stability}

By stability of reconstruction of the TAT tomogram $f$ from the measured data
$g$ we mean that small variations of $g$ in an appropriate norm lead to small
variations of the reconstructed tomogram $f$, also measured by an appropriate
norm. In other words, small errors in the data lead to small errors in the
reconstruction.

We will try to give the reader a feeling of the general state of affairs with
stability, referring to the literature (e.g.,
\cite{StefUhlTAT,Palam_funk,KuKuTAT,AKQ,HKN}) for further exact details.

We will consider as functional spaces the standard Sobolev spaces $H^s$ of
smoothness $s$. We will also denote, as before, by $\W$ the operator
transforming the unknown $f$ into the data $g$.

Let us recall the notions of \textbf{Lipschitz and H\"{o}lder stability}. An
even weaker \textbf{logarithmic stability} will not be addressed here. The
reader can find discussion of the general stability notions and issues, as
applied to inverse problems, in \cite{Isakov}.

\begin{definition}\label{D:stab}
The operation of reconstructing $f$ from $g$ is said to be \textbf{Lipschitz
stable} between the spaces $H^{s_2}$ and  $H^{s_1}$, if the following estimate holds for some constant $C$:
$$
\|f\|_{H^{s_1}}\leq C \|g\|_{H^{s_2}}.
$$

The reconstruction is said to be \textbf{H\"{o}lder stable} (a weaker
concept), if there are constants $s_1,s_2, s_3, C, \mu>0$, and $\delta >0$ such that
$$
\|f\|_{H^{s_1}}\leq C \|g\|^\mu_{H^{s_2}}
$$
for all $f$ such that $\|f\|_{H^{s_3}}\leq \delta$.
\end{definition}

Stability can be also interpreted in the terms of the singular values $\sigma_j$ of the forward operator $f\mapsto g$ in $L^2$, which have at most power decay when
$j\to\infty$. The faster is the decay, the more unstable the reconstruction
becomes. The problems with singular values decaying faster than any power of
$j$ are considered to be extremely unstable. Even worse are the problems with
exponential decay of singular values (analytic continuation or solving Cauchy
problem for an elliptic operator belong to this class). Again, the book
\cite{Isakov} is a good source for finding detailed discussion of such issues.

Consider as an example inversion of the standard in X-ray CT and MRI Radon transform that integrates a function $f$ over hyper-planes in $\R^n$. It smoothes function by ``adding $(n-1)/2$ derivatives.'' Namely, it maps continuously $H^s$-functions in
$\Omega$ into the Radon projections of class $H^{s+(n-1)/2}$. Moreover, the
reconstruction procedure is Lipshitz stable between these spaces (see
\cite{Natt_old} for detailed discussion).

One should notice that since the forward mapping is smoothing (it ``adds
derivatives'' to a function), the inversion should produce functions that are
less smooth than the data, which is an unstable operation. The rule of thumb
is that the stronger is smoothing, the less stable is inversion (this can be
rigorously recast in the language of the decay of singular values). Thus,
problems that require reconstructing non-smooth functions from infinitely
differentiable (or even worse, analytic) data, are extremely unstable (with
super-algebraic or exponential decay of singular values correspondingly). This
is just a consequence of the standard Sobolev embedding theorems (see, e.g.,
how this applies in TAT case in \cite{nguyen_stab}).

In the case of a constant sound speed and the acquisition surface completely
surrounding the object, as we have mentioned before, the TAT problem can be
recast as inversion of the spherical mean transform $\M$ (see Section
\ref{S:models}). Due to analogy between the spheres centered on $S$ and
hyperplanes, one suspects that  \textbf{inversion of the spherical mean
operator $\M$ is as Lipschitz stable as the inversion of the Radon transform}.
This indeed is the case, \textbf{as long as $f$ is supported inside $S$}, as
can be found in \cite{Palam_funk}. In the cases when closed form inversion
formulas are available (see Section \ref{SS:constantspeed}), this stability
can also be extracted from them. If the support of $f$ does reach outside,
\textbf{reconstruction of the part of $f$ that is outside is unstable} (i.e.,
is not even H\"{o}lder stable, due to the reasons explained in Section
\ref{S:incomplete}).

In the case of \textbf{variable non-trapping speed of sound} $c(x)$, integral
geometry does not apply anymore, and one needs to address the issue using, for
instance, time reversal. In this case, stability follows by solving the wave
equation in reverse time starting from $t=\infty$, as it is done in \cite{AK}.
In fact, \textbf{Lipschitz stability in this case holds for any observation
time exceeding $T(\Omega)$} (see \cite{StefUhlTAT}, where microlocal analysis
is used to prove this result).

The bottom line is that \textbf{ TAT reconstruction is sufficiently stable, as long as the speed of sound is non-trapping}.

However, trapping speed does cause instability \cite{HKN}. Indeed, since some of the rays are trapped inside $\Omega$, the information about some singularities
never reaches $S$ (no matter for how long one collects the data), and thus, as
it is shown in \cite{nguyen_stab}, the reconstruction is not even H\"{o}lder
stable, and the singular values have super-algebraic decay. See also Section
\ref{S:incomplete} below for a related discussion.

\subsection{Incomplete data}\label{S:incomplete}

In the standard X-ray CT, incompleteness of data arises, for instance, if not all projection angles
are accessible, or irradiation of certain regions is avoided, or as in the ROI (region of interest) imaging, only the ROI is irradiated.

It is not that clear what incomplete data means in TAT. Usually one says that one deals with \textbf{incomplete TAT data, if the acquisition surface does
not surround the object of imaging completely.} For instance, in breast
imaging it is common that only a half-sphere arrangement of transducers is
possible. We will see, however, that \textbf{incomplete data effects in TAT can also arise due to trapping, even if the acquisition surface completely surrounds the object}.

The questions addressed here are:
 \begin{enumerate}
   \item Is the collected incomplete data sufficient for \textbf{unique
    reconstruction}?

   \item If yes, does the incompleteness of the data have any effect on
    \textbf{stability and quality of the reconstruction}?
 \end{enumerate}

\subsubsection{Uniqueness of reconstruction}\label{S:limited_unique}

Uniqueness of reconstruction issues can be considered essentially resolved for
incomplete data in TAT, at least in most situations of practical interest. We
will briefly survey here some of the available results. In what follows, the
acquisition surface $S$ is not closed (otherwise the problem is considered to
have complete data).

\paragraph{Uniqueness for acoustically homogeneous media}\indent

In this case, Theorem \ref{T:uniq-simple} contains some useful sufficient
conditions on $S$ that guarantee uniqueness. Microlocal results of
\cite{AQ,LQ,StefUhl}, as well as the PDE approach of \cite{FPR} further
applied in \cite{AmbKuc_inj} provide also some other conditions. We assemble
some of these in the following theorem:
\begin{theorem}\label{T:uniq_incomp}
Let $S$ be a non-closed acquisition surface in TAT. Each of the following
conditions on $S$ is sufficient for the uniqueness of reconstruction of any
compactly supported function $f$ from the TAT data collected on $S$:
\begin{enumerate}
  \item Surface $S$ is not algebraic (i.e., there is no non-zero polynomial
vanishing on $S$).

  \item Surface $S$ is a uniqueness set for harmonic polynomials (i.e., there
is no non-zero harmonic polynomial vanishing on $S$).

  \item Surface $S$ contains an open piece of a closed analytic surface
$\Gamma$.

  \item Surface $S$ contains an open piece of an analytic surface $\Gamma$ separating the space $\R^n$ such that $f$ is supported on one
side of $\Gamma$.

  \item For some point $y\in S$, the function $f$ is supported on one side of the tangent plane $T_y$ to $S$ at $y$.
\end{enumerate}
\end{theorem}

For instance, if the acquisition surface $S$ is just a tiny non-algebraic piece of a surface, data collected on $S$ determines the tomogram $f$ uniquely. However, one realizes that such data is unlikely to be useful for any practical
reconstruction. Here the issue of stability of reconstruction kicks in, as it will be discussed in the stability sub-section further down.

\paragraph{Uniqueness for acoustically inhomogeneous media}\indent

In the case of a variable speed of sound, there still are uniqueness theorems
for partial data \cite{StefUhlTAT,Dustin_uniq}, e.g.

\begin{theorem}\label{T:uniq_incomp_var}\cite{StefUhlTAT}
Let $S$ be an open part of the boundary $\partial\Omega$ of a strictly convex domain $\Omega$ and the smooth speed of sound equals $1$ outside $\Omega$. Then the TAT data
collected on $S$ for a time $T>T(\Omega)$ determines uniquely any function
$f\in H^1_0(\Omega)$, whose support does not reach the boundary.
\end{theorem}
A modification of this result that does not require strict convexity is also
available in \cite{Dustin_uniq}.

While useful uniqueness of reconstruction results exist for incomplete data
problems, all such problems are expected to show instability. This issue is discussed in the sub-sections below. This will also lead to a better understanding of
incomplete data phenomena in TAT.

\subsubsection{``Visible'' (``audible'') singularities}\label{S:visible}

According to the discussion in Section \ref{S:unique_var}, the singularities
(the points of the wave front set $WF(f)$ of the function $f$ in
(\ref{E:wave_data})) are transported with time along the bi-characteristics
(\ref{E:bichar}). Thus, in the $x$-space they are transported along the
geometric rays. These rays may or may not reach the acquisition surface $S$, which triggers the introduction of the following notion:

\begin{definition}
A phase space point $(x_0,\xi_0)$ is said to be \textbf{``visible''} (sometimes the word
\textbf{``audible''} is used instead), if the corresponding ray (see (\ref{E:bichar})) reaches in finite time the observation surface $S$.

A region $U\subset \R^n$ is said to be in the \textbf{visibility zone}, if all
points $(x_0,\xi_0)$ with $x_0\in U$ are visible.
\end{definition}

\begin{figure}[ht!]
\begin{center}
\begin{tabular}{ccc}
\includegraphics[width=1.2in,height=1.2in]{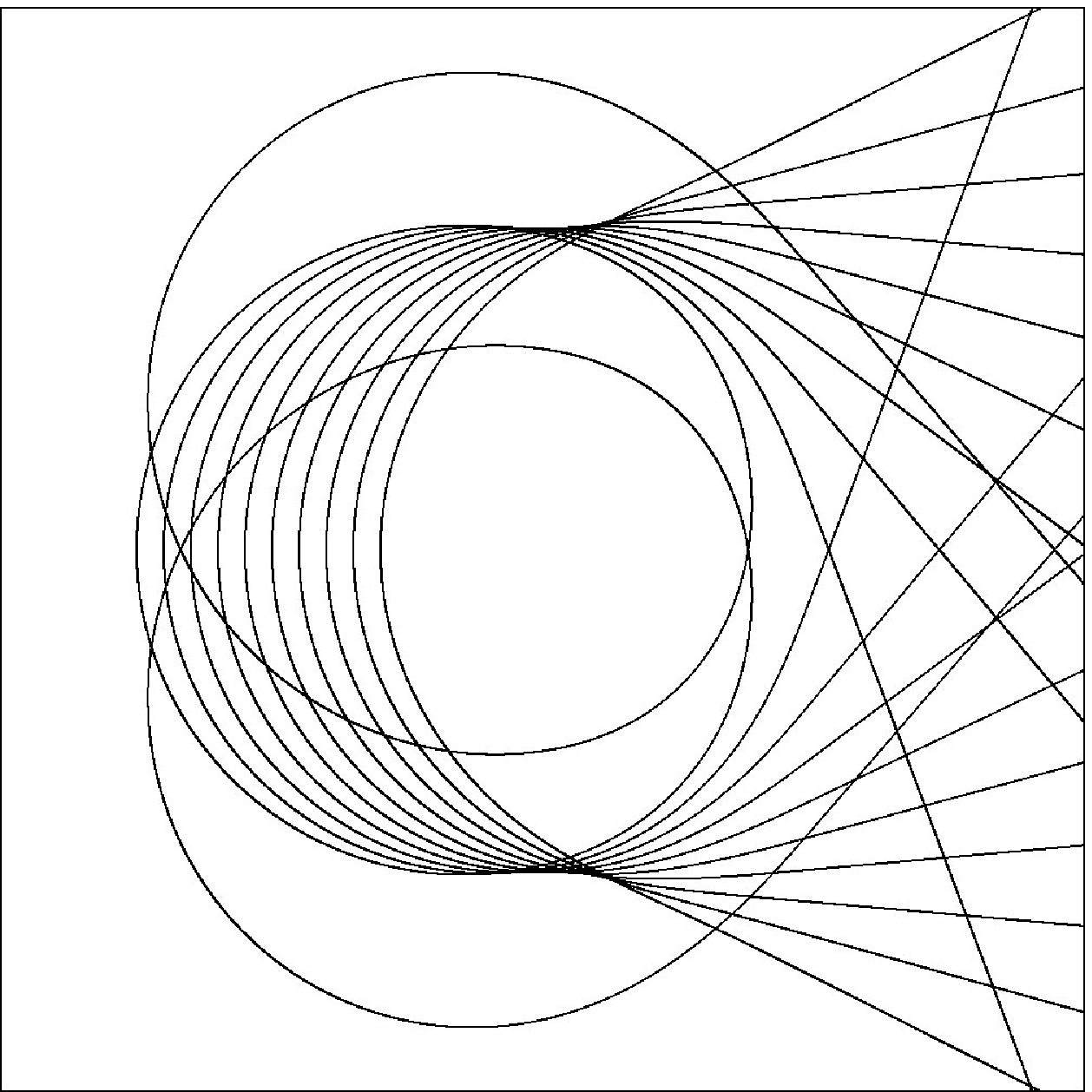} &
\includegraphics[width=1.2in,height=1.2in]{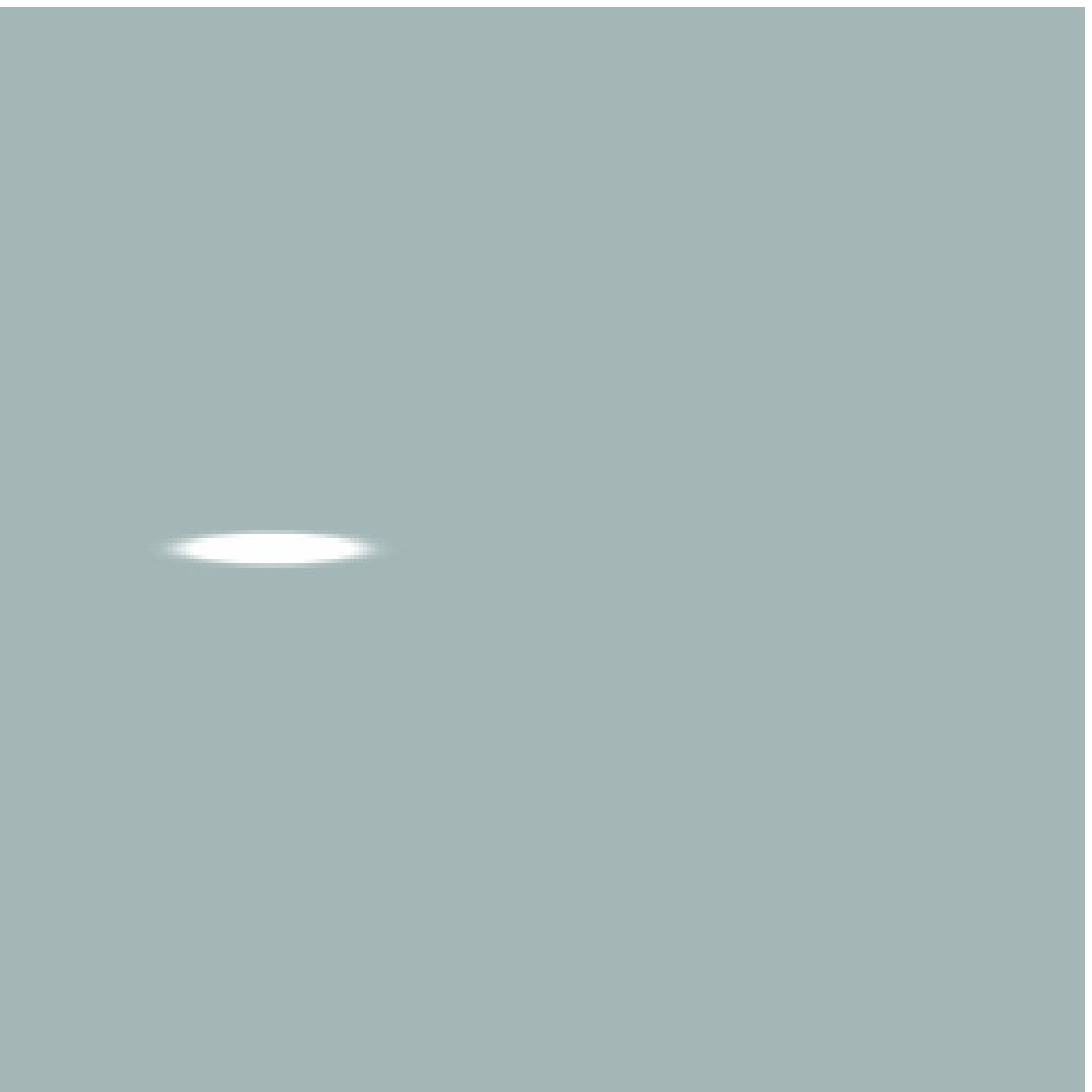} &
\includegraphics[width=1.2in,height=1.2in]{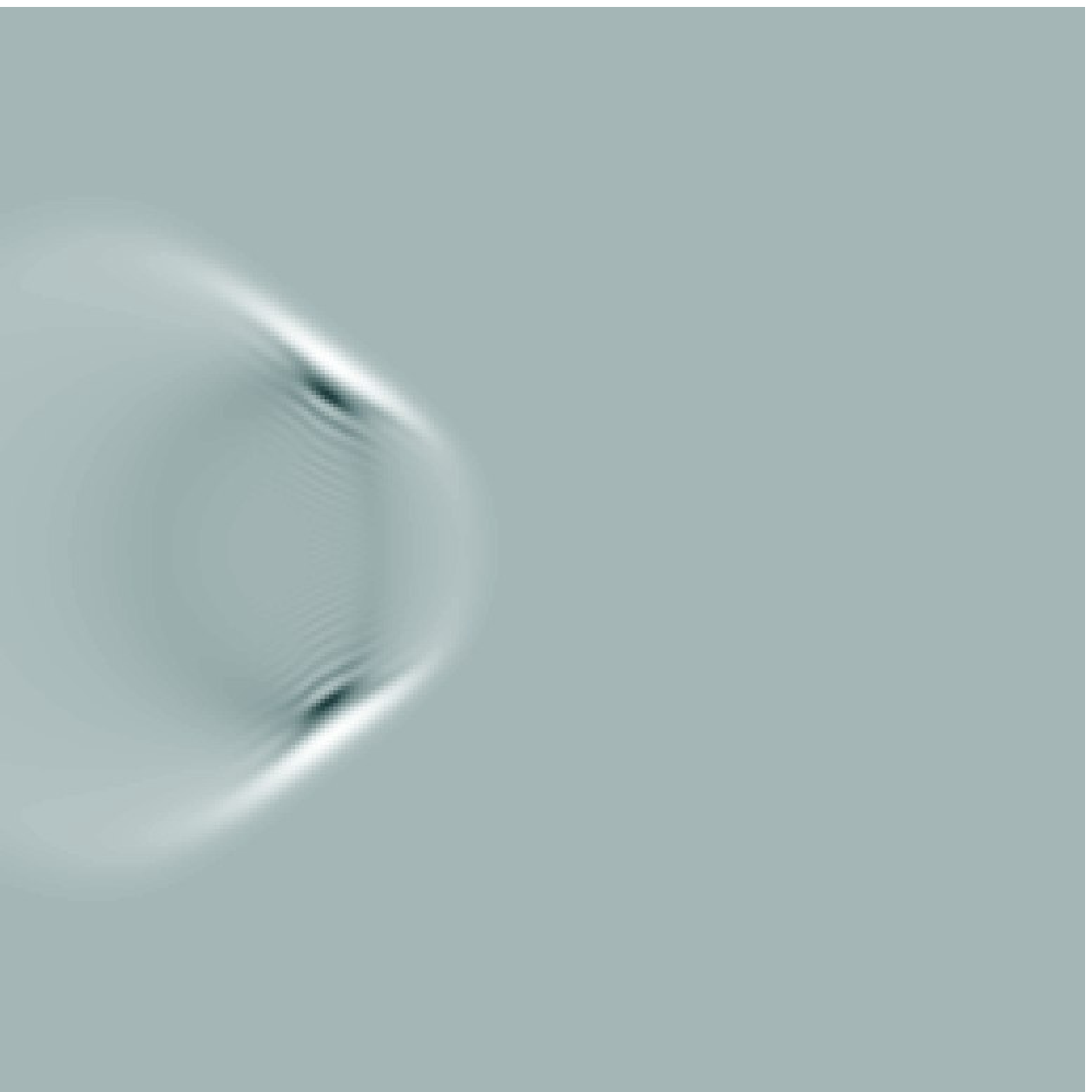}\\
(a) & (b) & (c) \\
\\
\includegraphics[width=1.2in,height=1.2in]{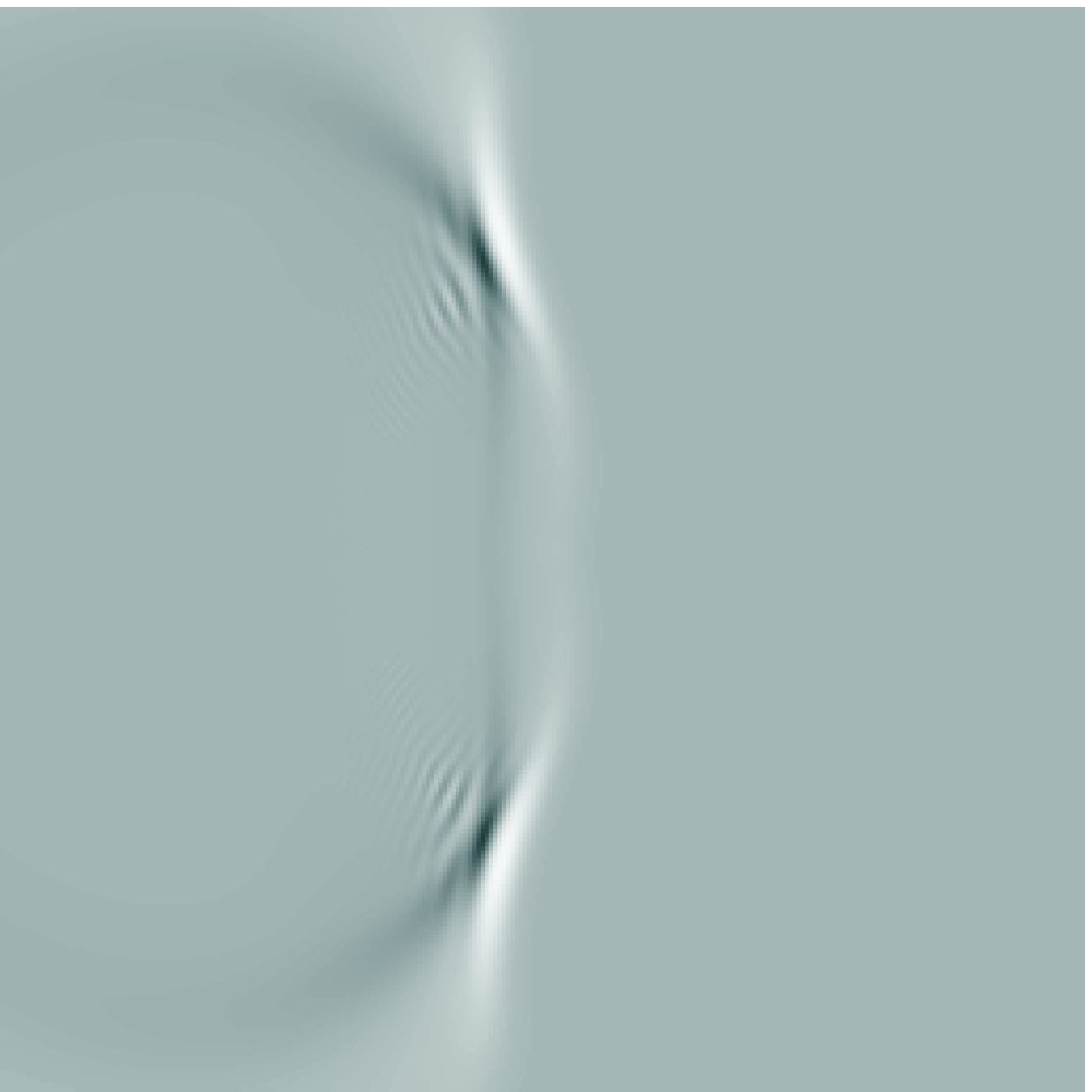} &
\includegraphics[width=1.2in,height=1.2in]{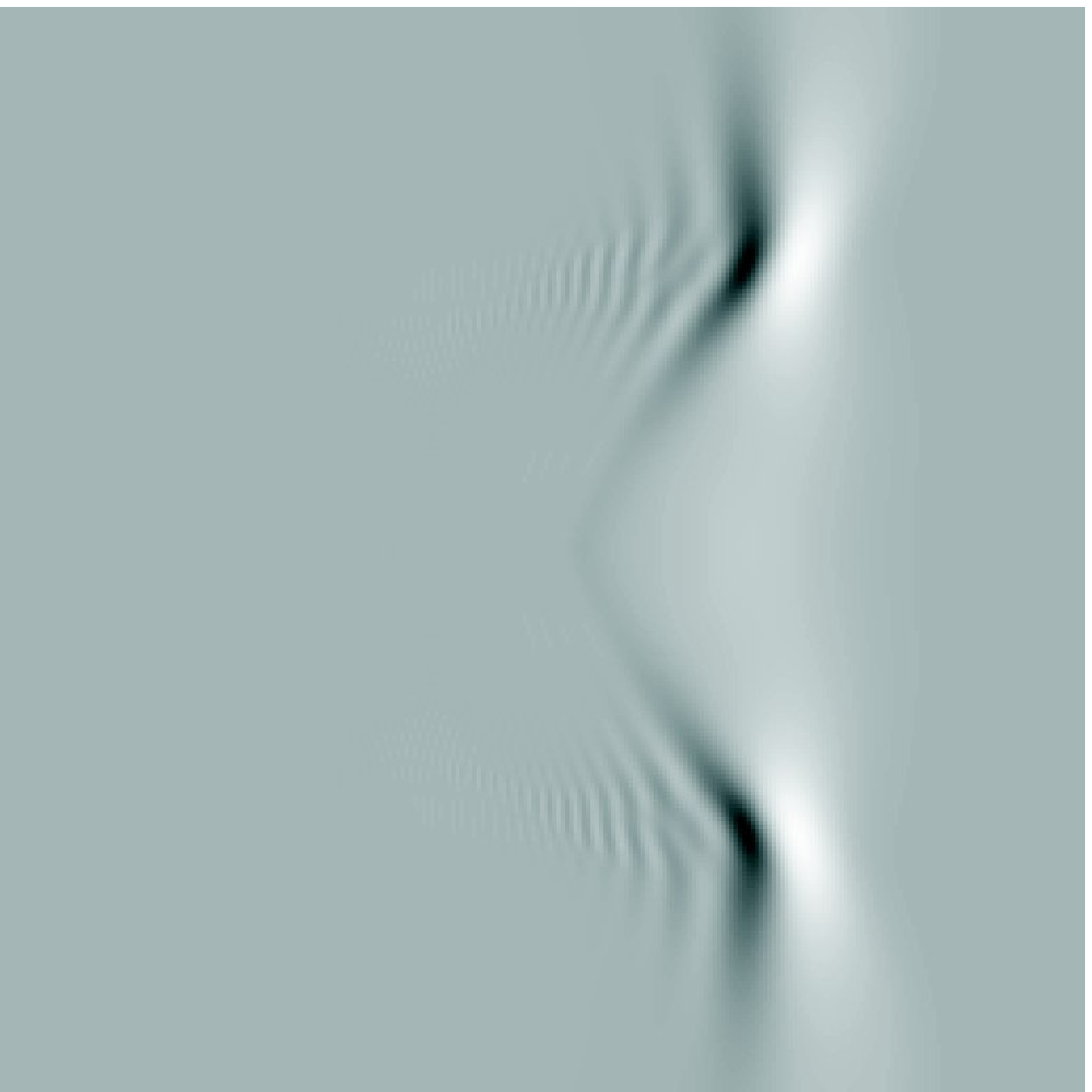} &
\includegraphics[width=1.2in,height=1.2in]{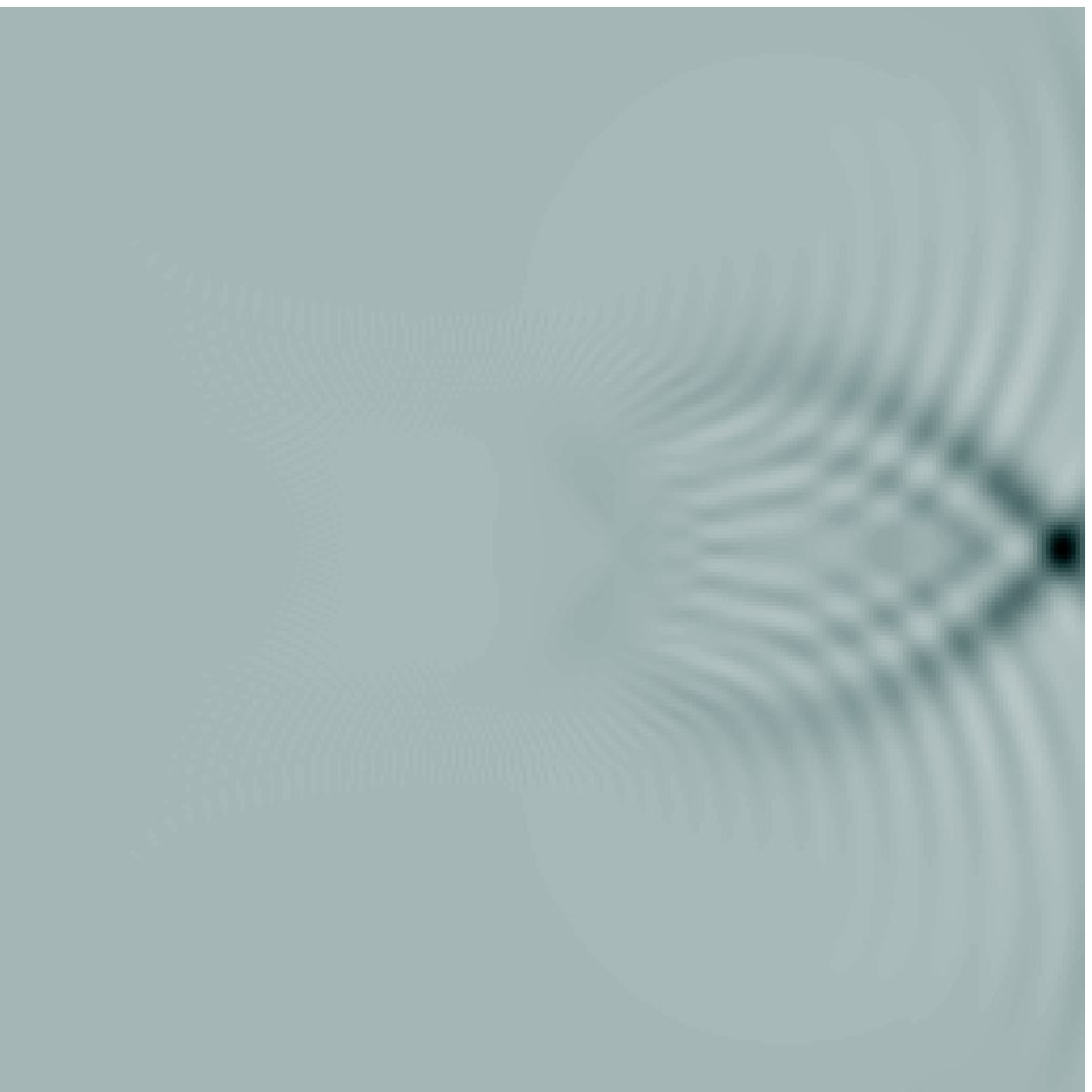} \\
 (d) & (e)& (f)
\end{tabular}
\end{center}
\caption{(a) Some rays starting along the interval $x \in [-0.7,-0.2]$ in the
vertical directions escape on the right; (b) a flat phantom with ``invisible
wavefront''; (c-f)propagation of the flat front: most of the energy of the
signal leaves the square domain through the hole on the right.}\label{F:bent}
\end{figure}

An example of wave propagation through inhomogeneous medium is presented in
Figure \ref{F:bent}. The open observation surface $S$ in this example
consists of the two horizontal and the left vertical sides of the square.
Figure \ref{F:bent}(a) shows some rays that bend, due to
acoustic inhomogeneity, and leave through the opening of the observation surface $S$ (the right side of the square).
Fig. \ref{F:bent} (b) presents a flat phantom, whose
wavefront set creates these escaping rays, and thus is mostly invisible. Then
Fig. \ref{F:bent} (c-f) show the propagation of the corresponding wave front.

\begin{figure}[ht!]
\begin{center}
\begin{tabular}{cc}
\includegraphics[width=1.5in,height=1.5in]{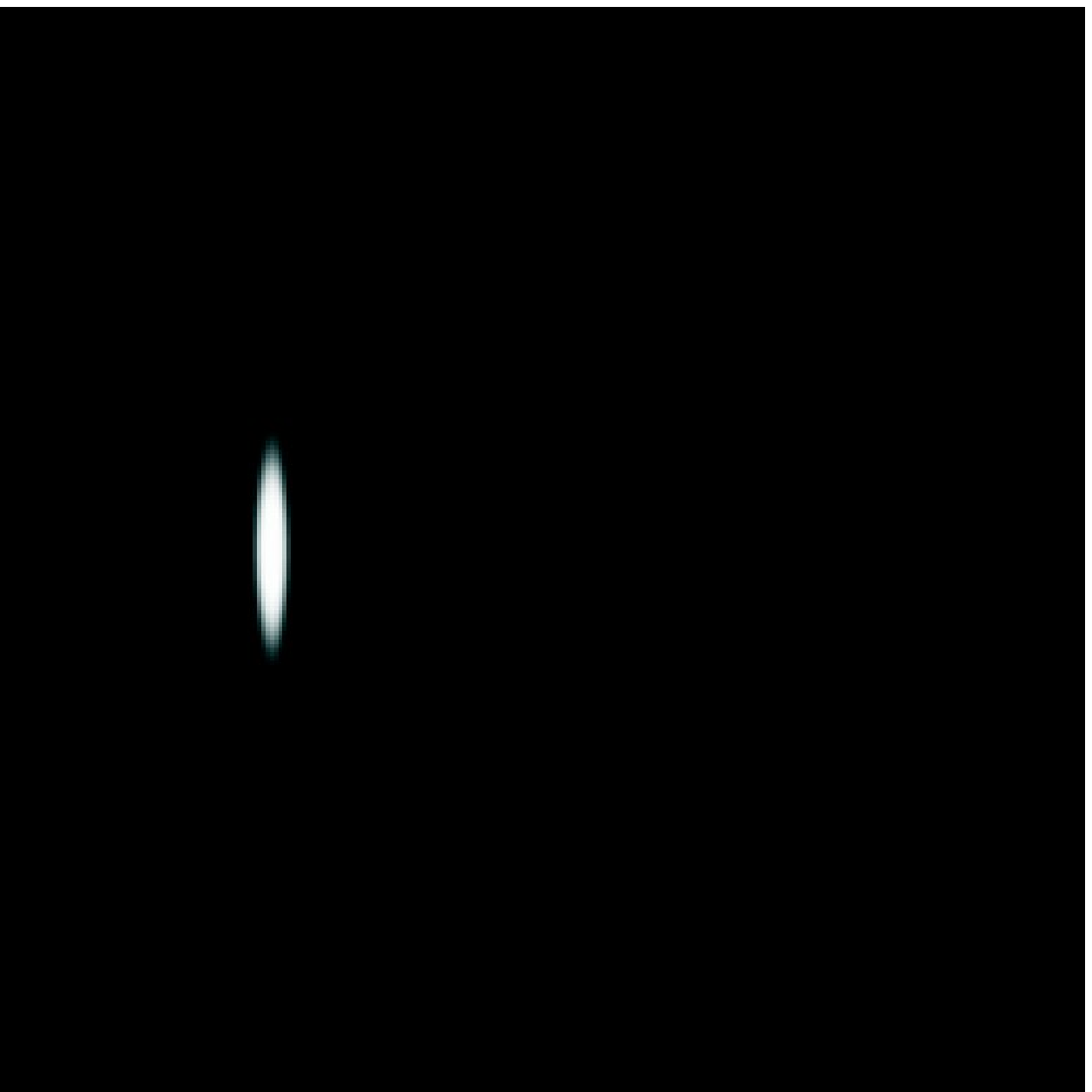} &
\includegraphics[width=1.5in,height=1.5in]{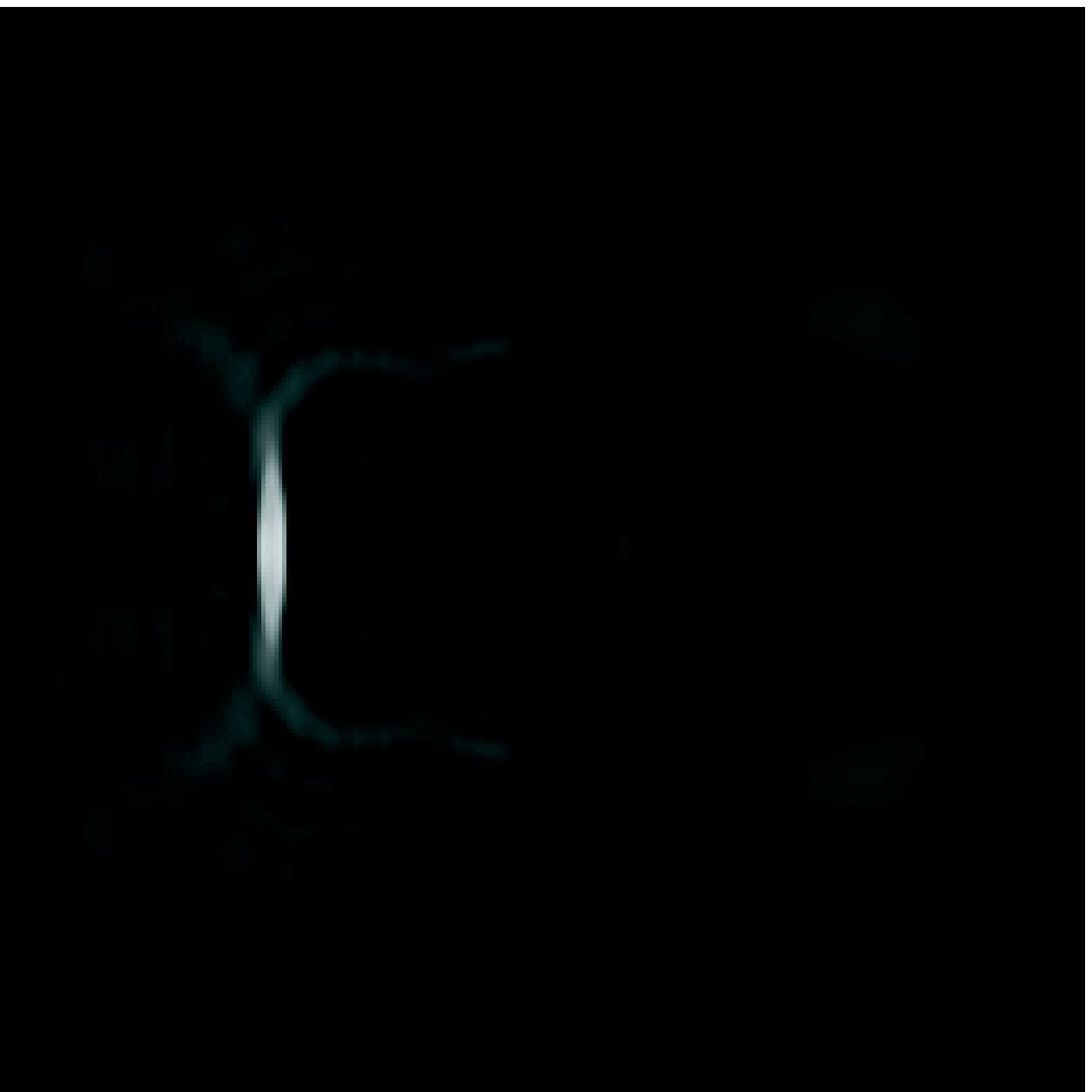} \\
(a) & (b) \\
\\
\includegraphics[width=1.5in,height=1.5in]{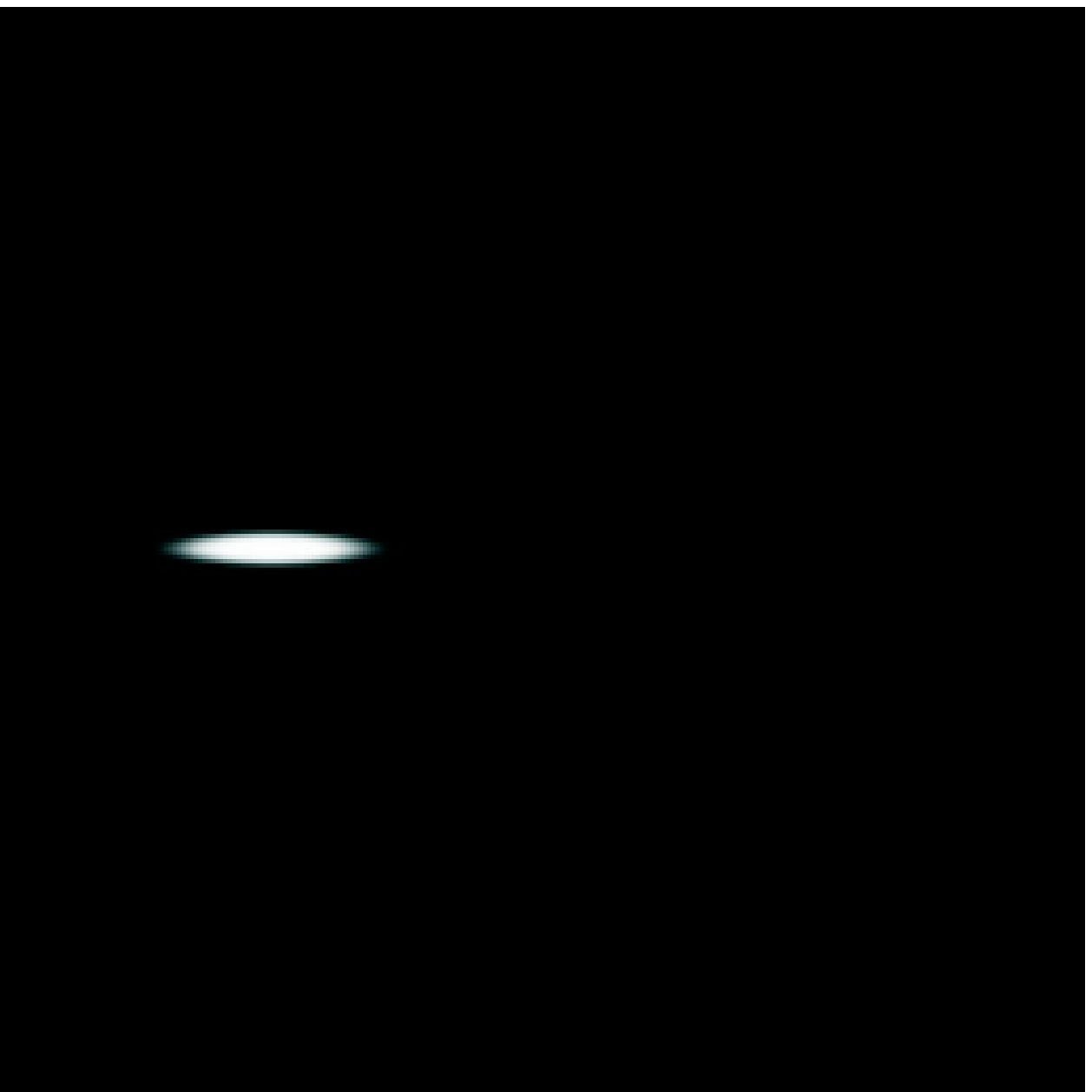} &
\includegraphics[width=1.5in,height=1.5in]{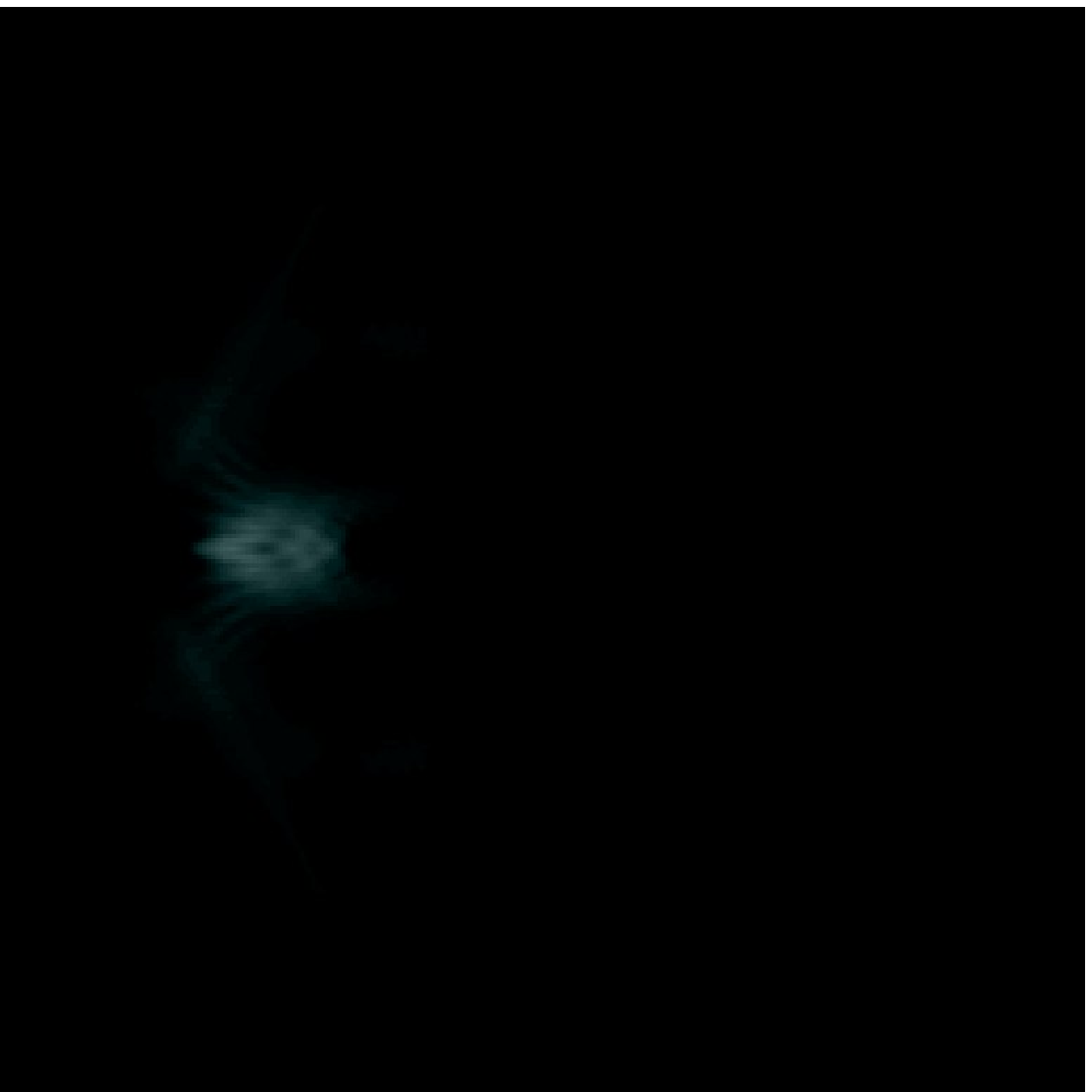} \\
(c) & (d)
\end{tabular}
\end{center}
\caption{Reconstruction with the same speed of sound:
(a-b) phantom with strong vertical fronts and its reconstruction;
(c-d) phantom with strong horizontal fronts and its reconstruction.}
\label{F:escape_rec}
\end{figure}

Since the information about the horizontal boundaries of the phantom escapes,
one does not expect to reconstruct it well. Fig. \ref{F:escape_rec} shows two
phantoms and their reconstructions from the partial data: (a-b) correspond to
the vertical flat phantom, whose only invisible singularities are at its ends.
One sees essentially good reconstruction, with a little bit of blurring at the
endpoints. On the other hand, reconstruction of the horizontal phantom with
almost the whole wave front set invisible, does not work.
The next Fig. \ref{F:square} shows a more complex square phantom, whose
singularities corresponding to the horizontal boundaries are invisible, while
the vertical boundaries are fine. One sees clearly that the invisible parts
have been blurred away.
\begin{figure}[ht!]
\begin{center}
\begin{tabular}{ccc}
\includegraphics[width=1.5in,height=1.5in]{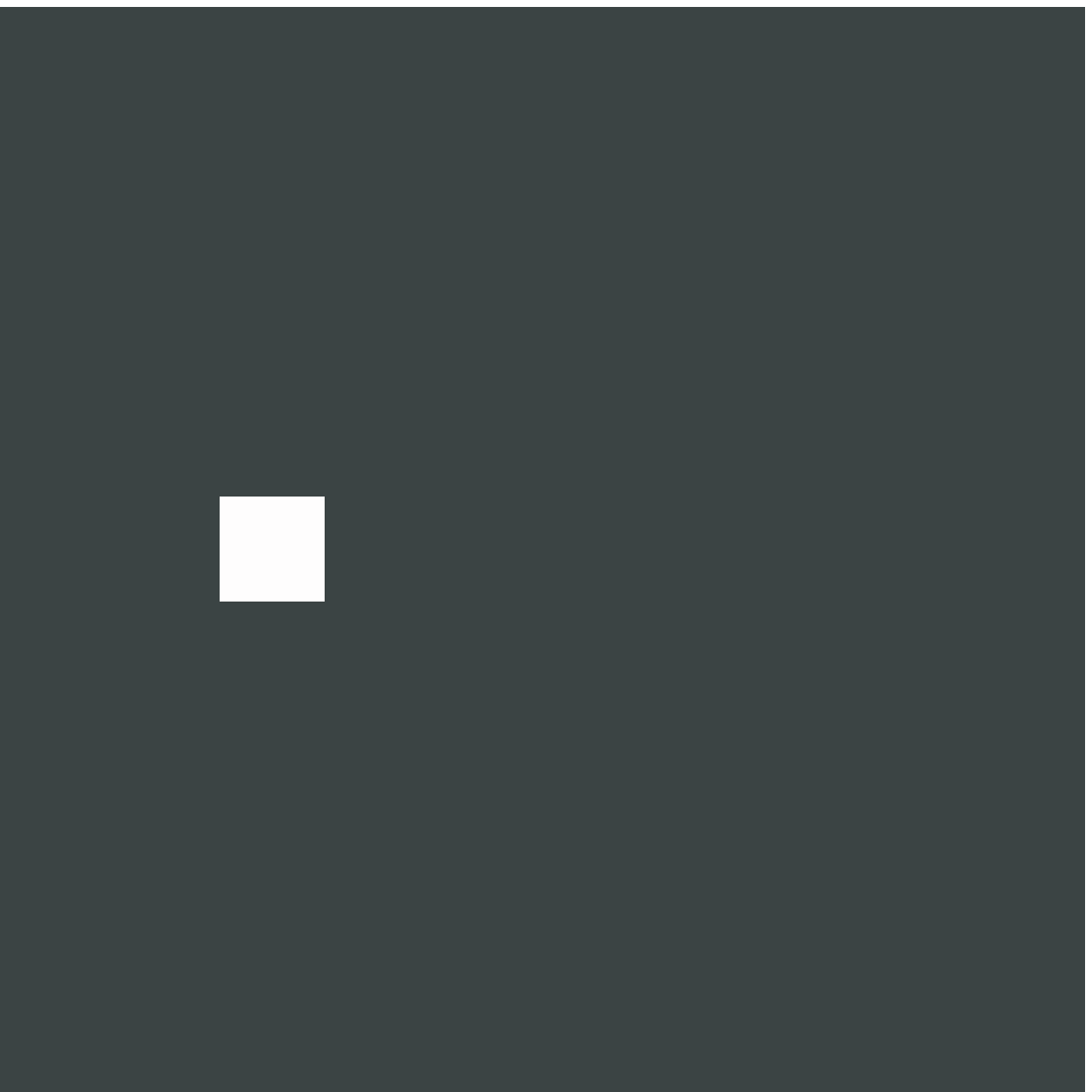} &
\includegraphics[width=1.5in,height=1.5in]{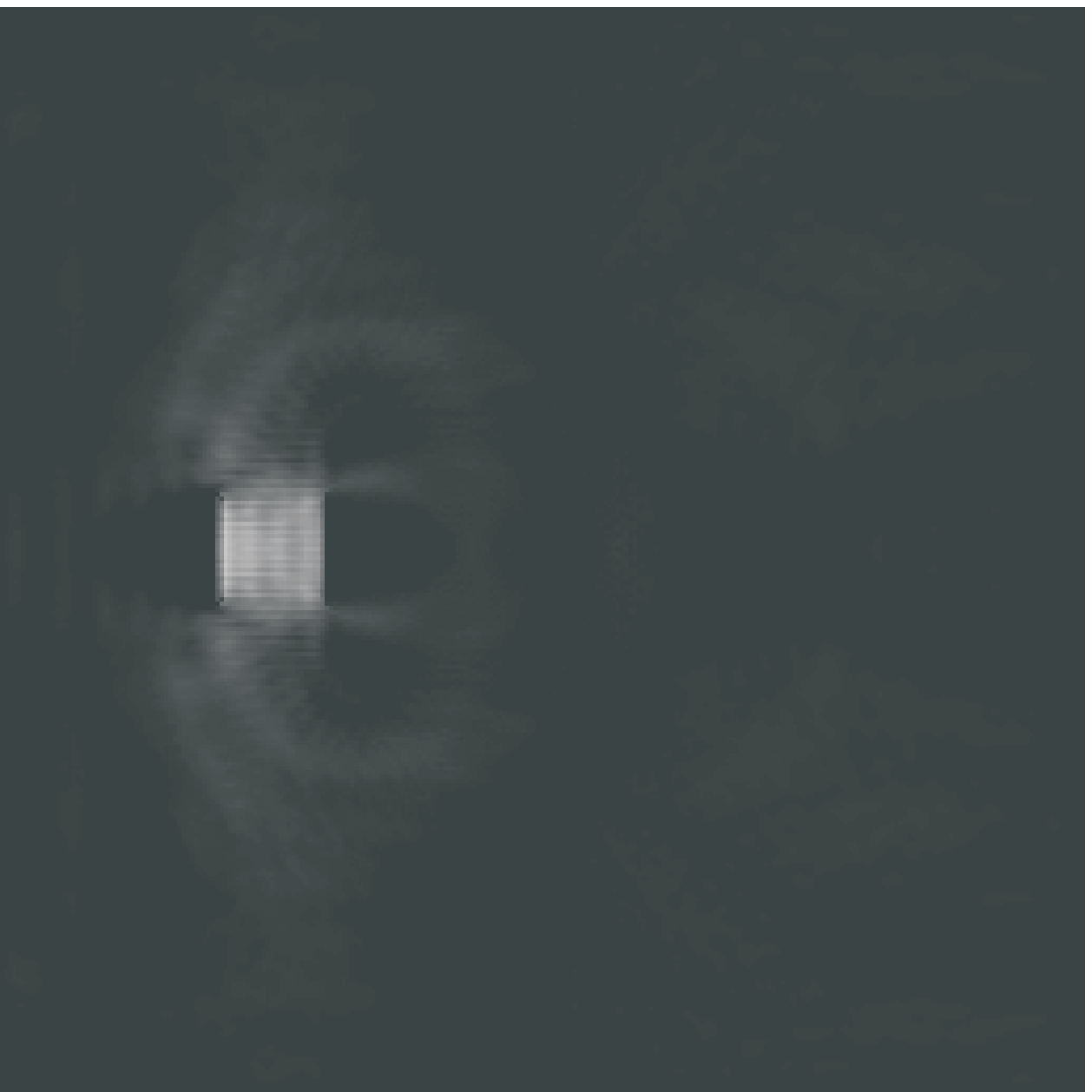} &
\includegraphics[width=1.5in,height=1.5in]{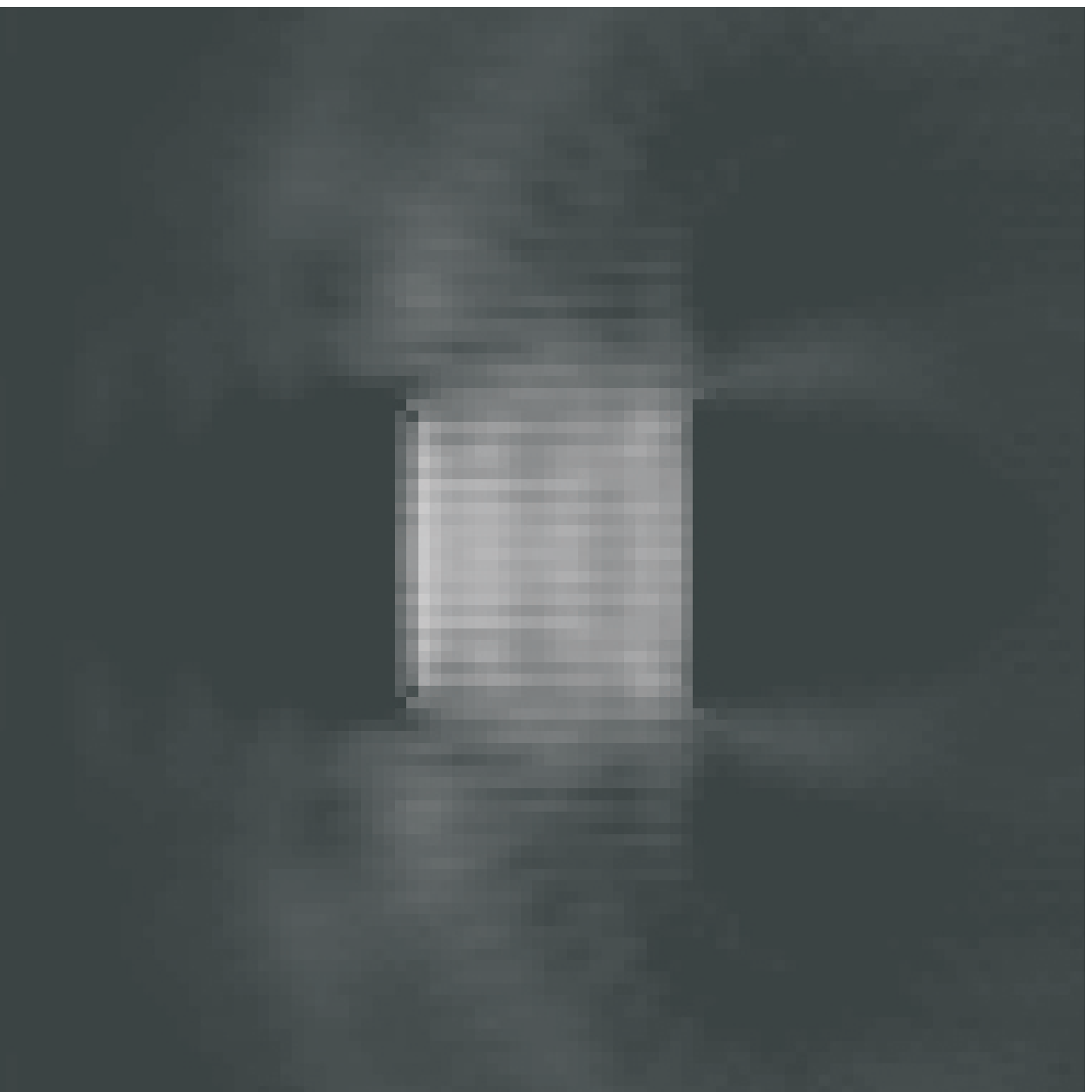} \\
(a) & (b) & (c)
\end{tabular}
\end{center}
\caption{Reconstruction with the same speed of sound:
(a) phantom; (b) its reconstruction; (c) a magnified fragment of
(b).}\label{F:square}
\end{figure}
On the other hand, Fig. \ref{F:incomplete}(a) in Section \ref{S:reconstr} shows
that one can reconstruct an image without blurring and with correct values, if
the image is located in the visibility region. The reconstructed image in this
figure is practically indistinguishable from the phantom shown in Figure \ref{F:2d}(a).

\begin{remark}\label{R:sing_outside}
If $S$ is a closed surface and $x_0$ is a point outside of $S$, there is a vector $\xi_0\neq 0$ such that $(x_0,\xi_0)$ is ``invisible.'' Thus, the visibility zone does not reach outside the closed acquisition surface $S$.
\end{remark}

\subsubsection{Stability of reconstruction for incomplete data problems}\label{SS:stable}\indent

In all examples above, uniqueness of reconstruction held,
but the images were still blurred.
The question arises whether the blurring of ``invisible'' parts is avoidable (after all, the uniqueness theorems seem to claim that ``everything
is visible''). The answer to this is, in particular, the following result of \cite{nguyen_stab}, which is an analog of similar statements in X-ray tomography:
\begin{theorem}\cite{nguyen_stab}\label{T:instab}
If there are invisible points $(x_0,\xi_0)$ in $\Omega\times(\R^n_\xi\setminus
0)$, then inversion of the forward operator $\W$ is not H\"{o}lder stable in
any Sobolev spaces. The singular values $\sigma_j$ of $\W$ in $L^2$ decay
super-algebraically.
\end{theorem}
Thus, having invisible singularities makes the reconstruction severely ill-posed. In particular, according to Remark \ref{R:sing_outside}, this theorem implies the following statement:

\begin{corollary}\label{C:unst_outside}
Reconstruction of the parts of $f(x)$ supported outside the closed observation surface $S$ is unstable.
\end{corollary}

On the other hand,
\begin{theorem}\cite{StefUhlTAT}\label{T:stab}
All visible singularities of $f$ can be reconstructed with Lipschitz stability
(in appropriate spaces).
\end{theorem}
Such a reconstruction of visible singularities can be obtained in many ways,
for instance just by replacing the missing data by zeros (with some smoothing along the junctions with the known data, in order to avoid artifact
singularities). However, there is no hope for stable recovery of the correct values of $f(x)$, if there are invisible singularities.

\subsection{Discussion of the visibility condition}\label{S:vis_discuss}

\paragraph{Visibility for acoustically homogeneous media}\indent

\begin{figure}[th]
\begin{center}
\begin{tabular}
[c]{ccc}
\includegraphics[width=1.5in,height=1.5in]{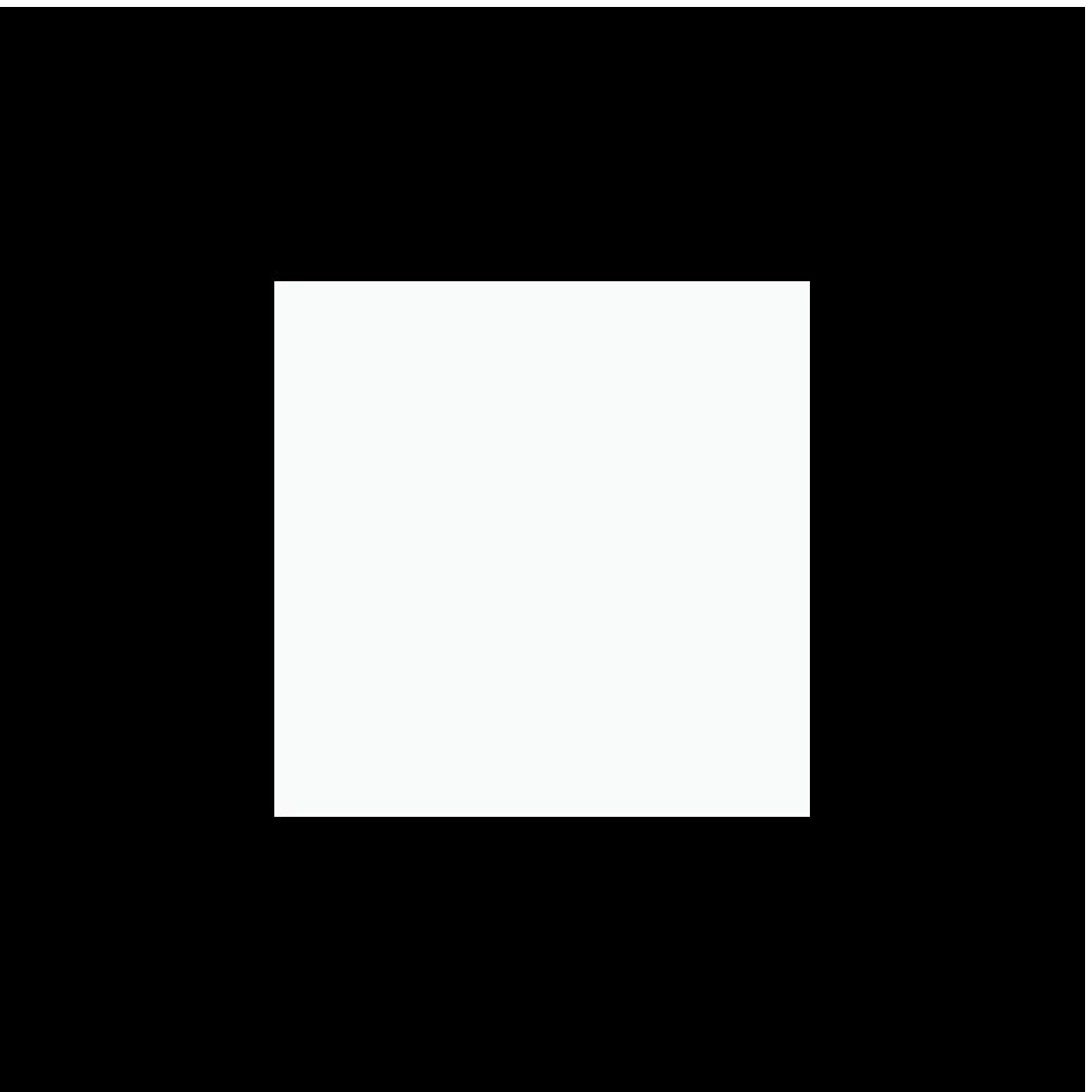} &
\includegraphics[width=1.5in,height=1.5in]{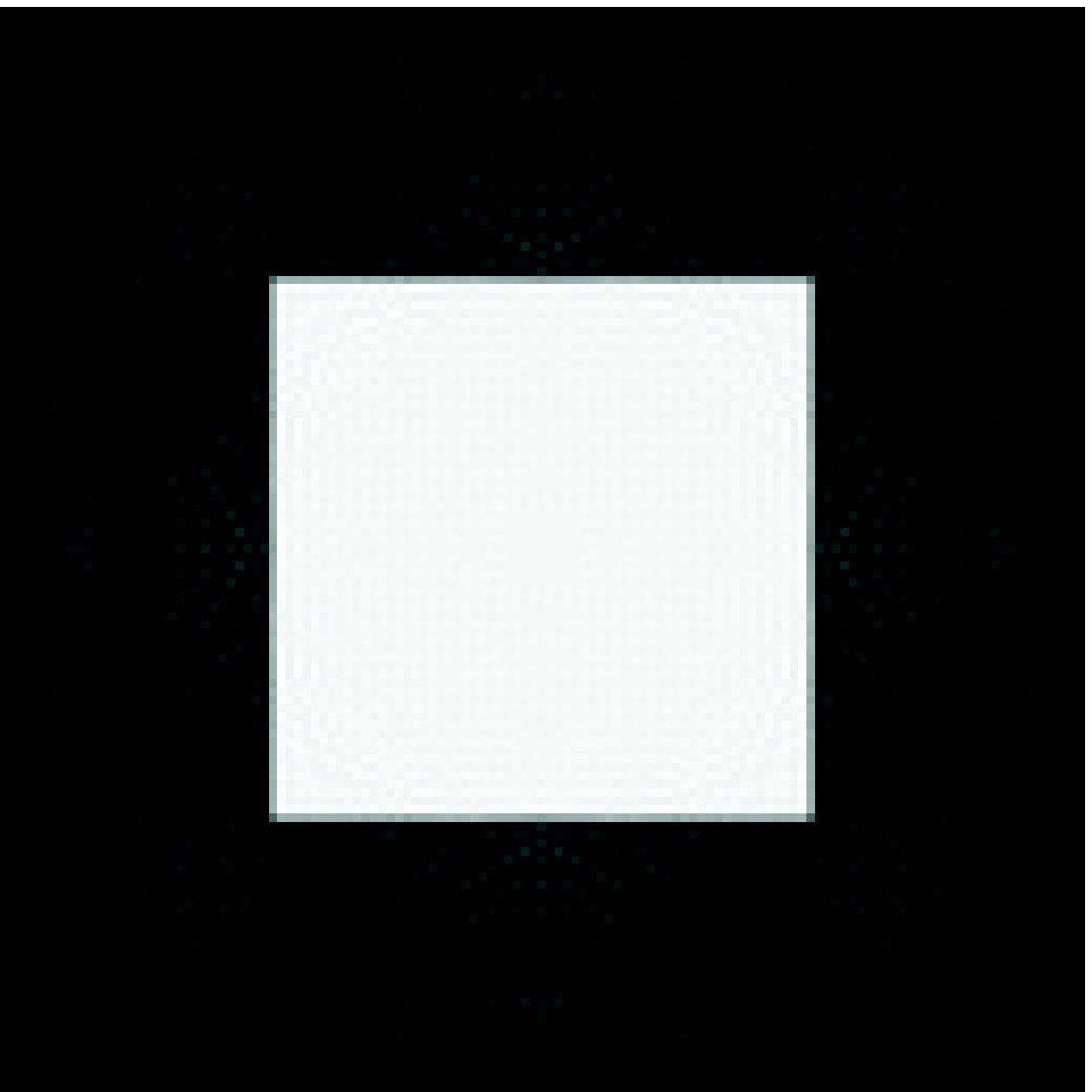} &
\includegraphics[width=1.5in,height=1.5in]{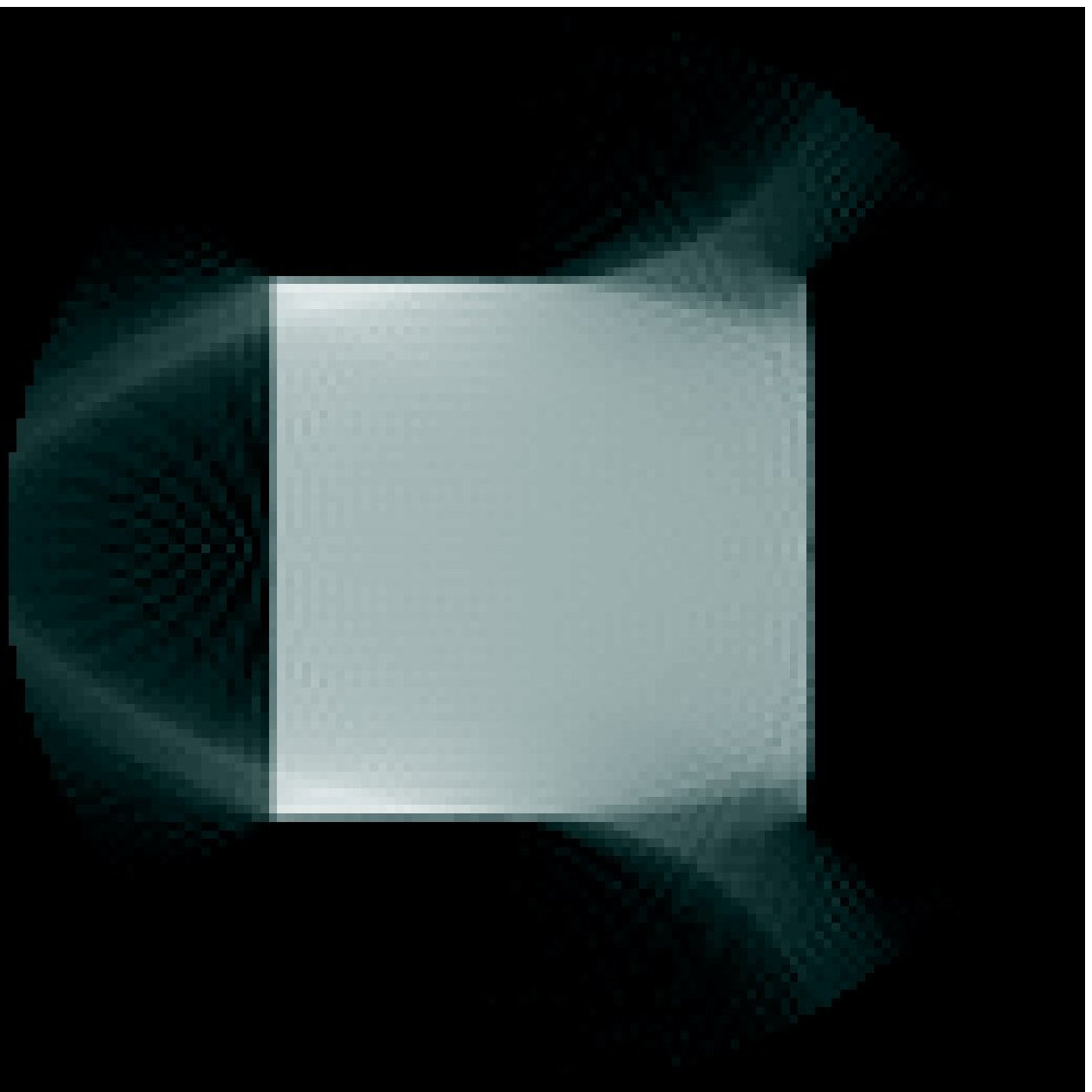}\\
(a) & (b) & (c)\\
\end{tabular}
\end{center}
\caption{Reconstruction from incomplete data using closed form inversion formula in $2D$; detectors are located
on the left half circle of radius 1.05 (a) phantom
(b) reconstruction from complete data
(c) reconstruction from the incomplete data}\label{F:partial}
\end{figure}

In the constant speed case, the rays are straight, and thus the visibility
condition has a simple test:
\begin{proposition}(e.g., \cite{XWAK,XWAK2,HKN})\label{P:vis_const}
If the speed is constant, a point $x_0$ is in the visible region, if and only
if any line passing through $x_0$ intersects at least once the acquisition
surface $S$ (and thus a detector location).
\end{proposition}
Figure \ref{F:partial} illustrates this statement. It shows a square phantom
and its reconstruction from complete data and from the data collected on the
half-circle $S$ surrounding the left half of object. The parts of the
interfaces where the normal to the interface does not cross $S$ are blurred.

\paragraph{Visibility for acoustically inhomogeneous media}\indent

When the speed of sound is variable, an analog of Proposition
\ref{P:vis_const} holds, with lines replaced by rays.

\begin{proposition}(e.g., \cite{HKN,StefUhlTAT,nguyen_stab})\label{P:vis_var}
A point $x_0$ is in the visible region, if and only if for any $\xi_0\neq 0$
at least one of the two geometric rays starting at $(x_0,\xi_0)$ and at
$(x_0,-\xi_0)$ (see (\ref{E:bichar})) intersects the acquisition surface $S$
(and thus a detector location).
\end{proposition}

The reader can now see an important difference between the acoustically
homogeneous and inhomogeneous media. Indeed, even if $S$ surrounds the support
of $f$ completely, trapped rays will never find their way to $S$, which will lead, as we know by now, to instabilities and blurring of some interfaces.

\begin{figure}[ht!]
\begin{center}
\begin{tabular}{cccc}
\includegraphics[width=1.5in,height=1.5in]{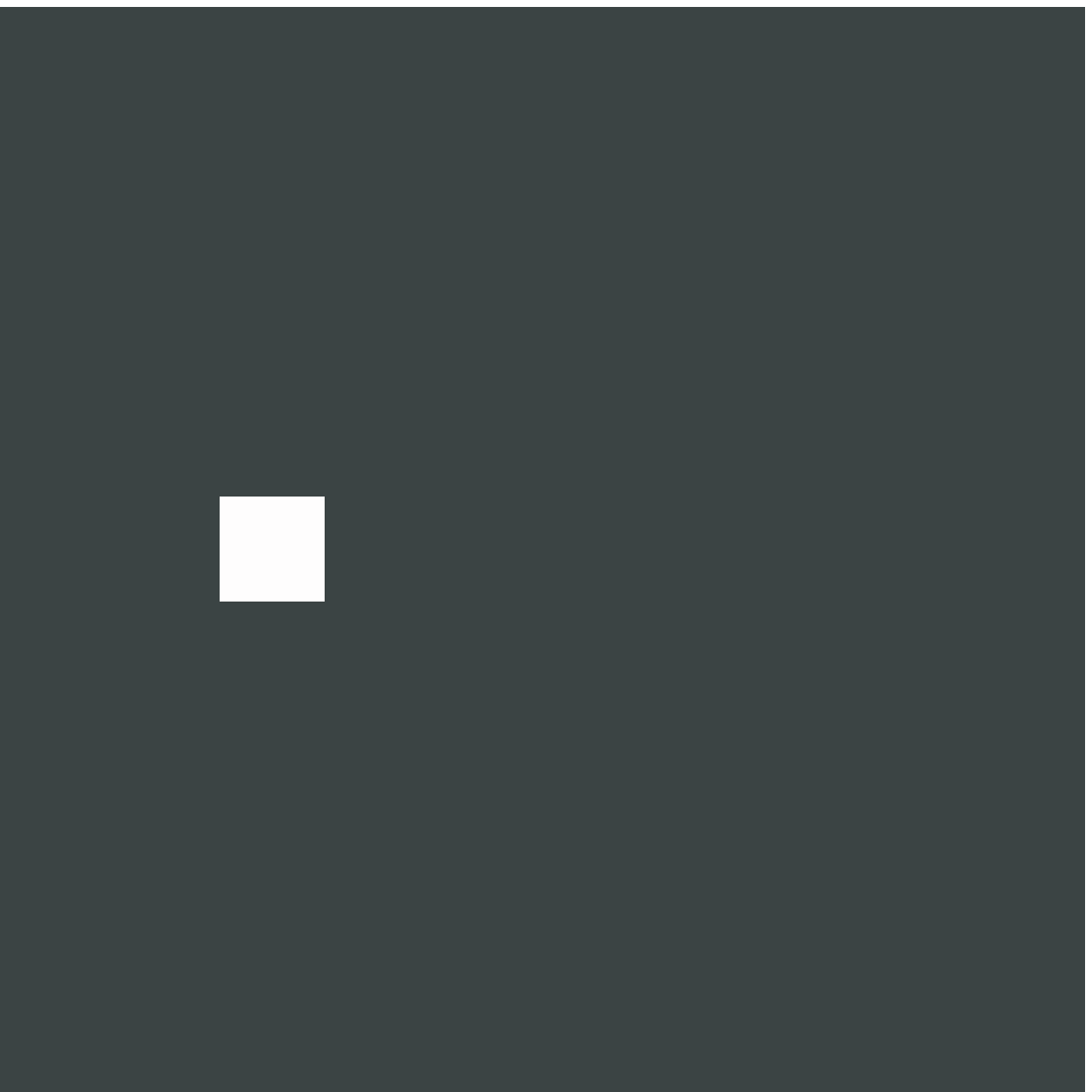} &
\includegraphics[width=1.5in,height=1.5in]{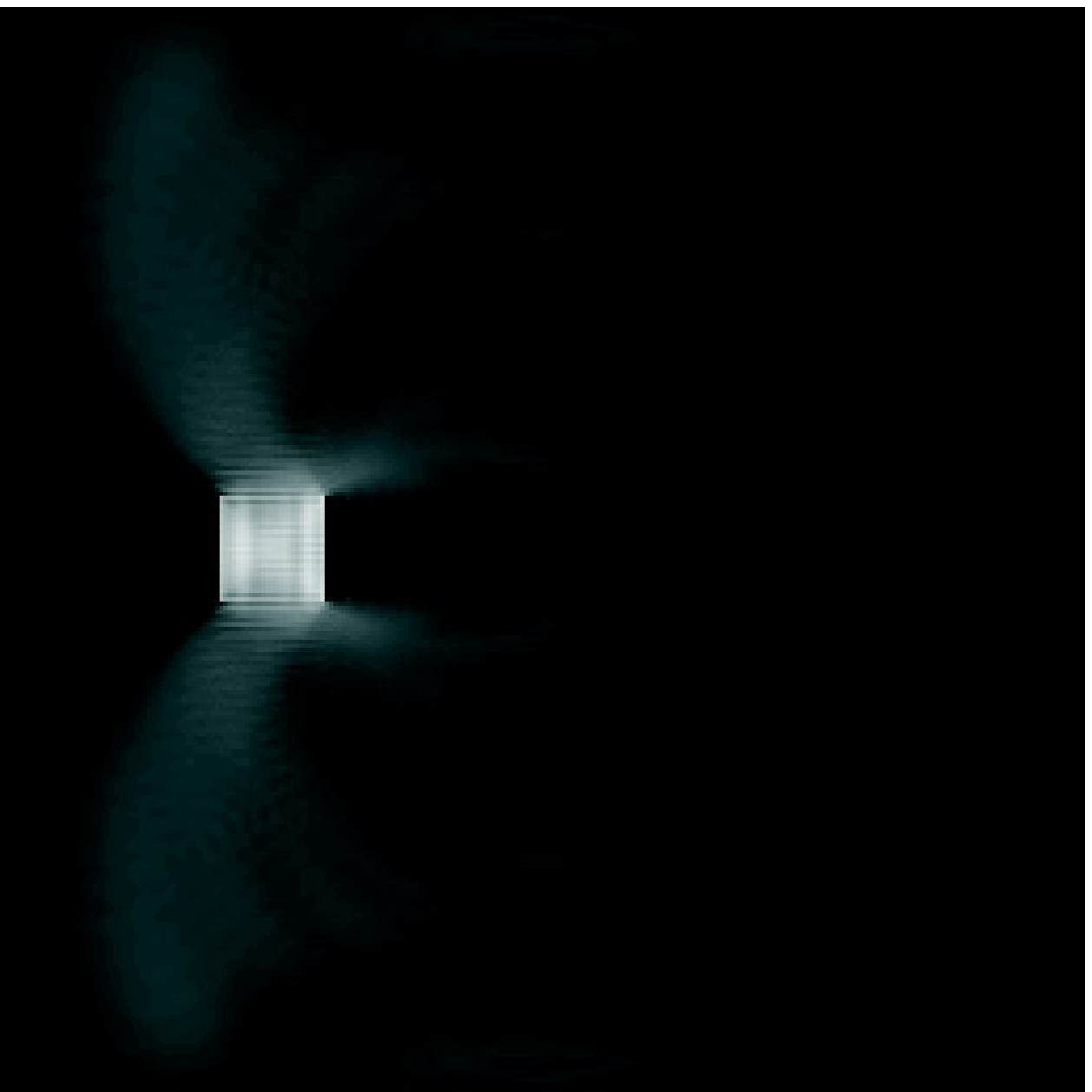} &
\includegraphics[width=1.5in,height=1.5in]{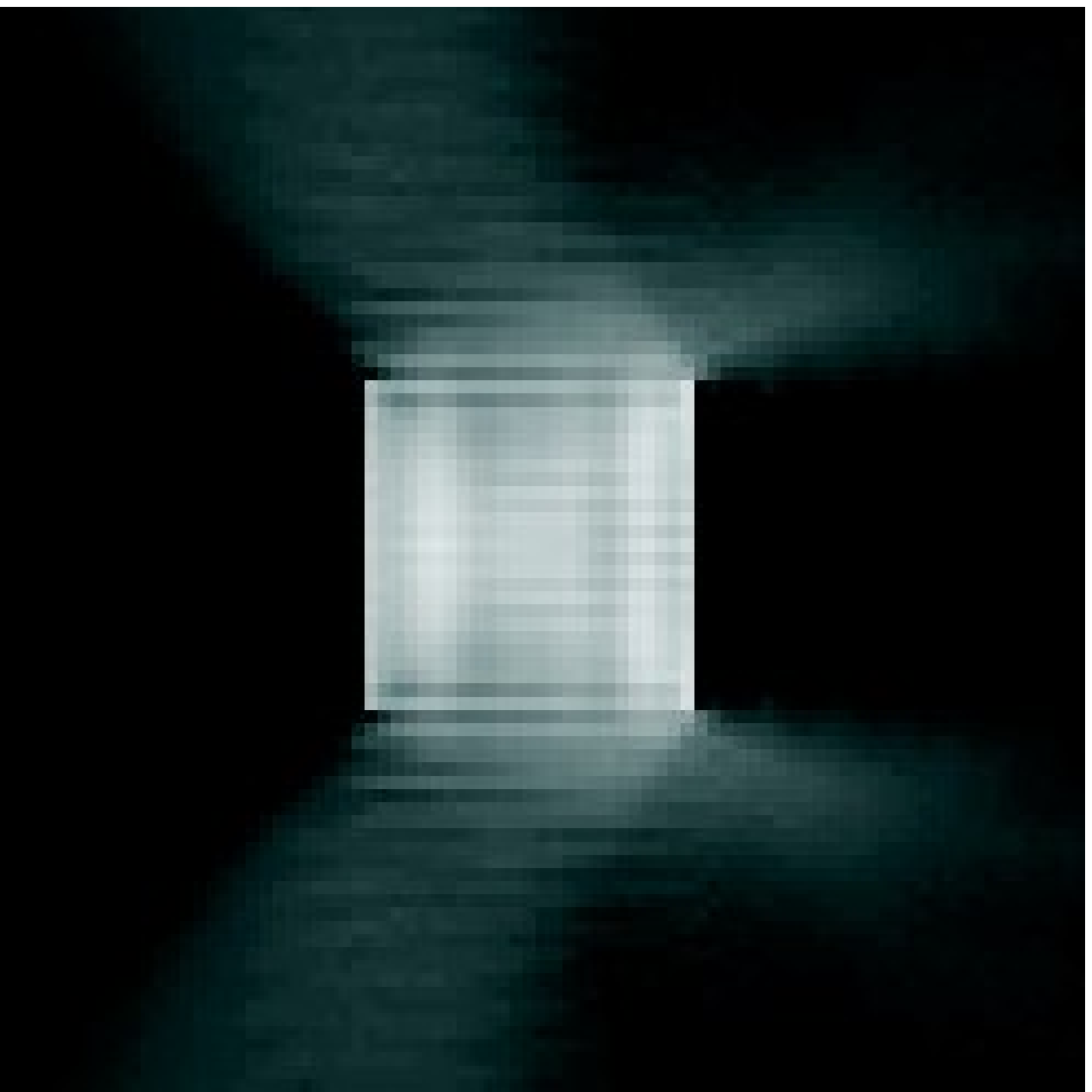} &
\includegraphics[width=1.5in,height=1.5in]{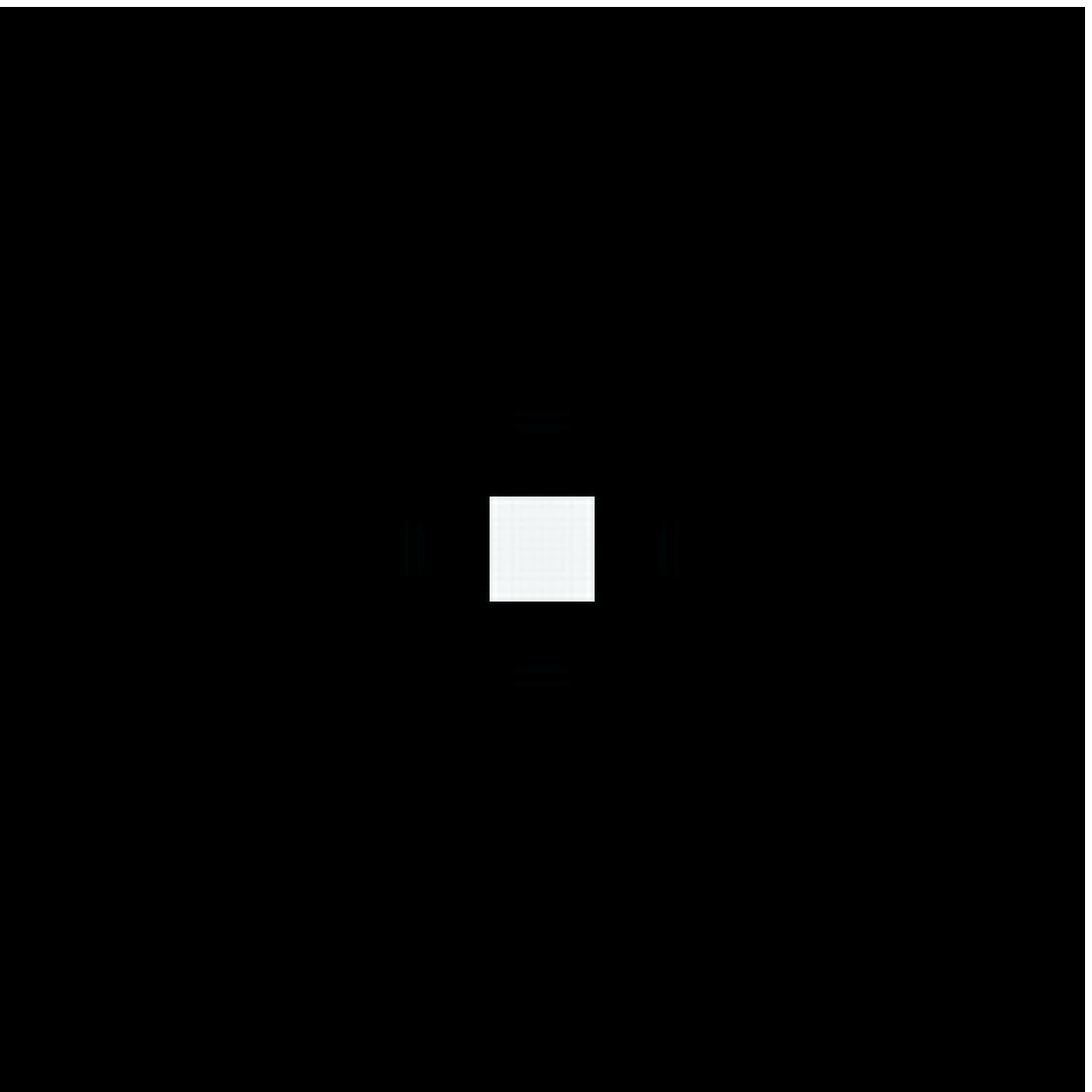}
\\
(a) & (b) & (c)&(d)
\end{tabular}
\end{center}
\caption{Reconstruction of a square phantom from full data in the presence of
a trapping parabolic speed of sound (the speed is radial with respect to the
center of the picture): (a) an off-center phantom; (b) its reconstruction; (c)
a magnified fragment of (b); (d) reconstruction of a centered square
phantom.}\label{F:trap_blur}
\end{figure}

Thus, {  presence of rays trapped
inside the acquisition surface creates effects of incomplete data type}. This is exemplified in Fig. \ref{F:trap_blur}
with a square phantom and its reconstruction shown in the presence of a
trapping (parabolic) speed. Notice that the square centered at the center of
symmetry of the speed is reconstructed very well (see (d)), since none of the
rays carrying its singularities is trapped.

\subsection{Range conditions}\label{S:range}

In this section we address the problem of describing the ranges of the forward
operators $\W$ (see (\ref{E:forward})) and $\M$ (see (\ref{E:Radon_S})), the
latter in the case of an acoustically homogeneous medium (i.e., for
$c=\mathrm{const}$). The ranges of these operators, similarly to the range of the Radon
and X-ray transforms (see \cite{Natt_old,Natt_new}), are of infinite
co-dimensions. This means that ideal data $g$ from a suitable function space
satisfy infinitely many mandatory identities. Knowing the range is useful for
many theoretical and practical purposes in various types of tomography
(reconstruction algorithms, error corrections, incomplete data completion,
etc.), and thus this topic has attracted a lot of attention (e.g.,
\cite{Leon_Radon,GGG1,GGG,GelfVil,Helg_Radon,Helg_groups,Kuch_AMS05,KuchQuinto,Natt_old,Natt_new,Pal_book,Q2006} and references therein).

As we will see in the next section, range descriptions in TAT are also intimately
related to recovery of the unknown speed of sound.

We recall \cite{Natt_old,Helg_Radon,GGG,GGG1} that for the standard Radon transform
$$
f(x)\to g(s,\omega)=\int\limits_{x\cdot\omega =s}f(x)dx, |\omega|=1,
$$
where $f$ is assumed to be smooth and supported in the unit ball $B=\{x\mid
|x|\leq 1\}$, the range conditions on $g(s,\omega)$ are:
\begin{enumerate}
\item {\em smoothness and support}: $g\in C^\infty_0\left([-1,1]\times
\mS\right)$, where $\mS$ is the unit sphere of vectors $\omega$,

\item {\em evenness}: $g(-s,-\omega)=g(s,\omega)$,

\item {\em moment conditions}: for any integer $k\geq 0$, the $k$th
moment
$$
G_k(\omega)=\int\limits_{-\infty}^{\infty}
s^k g(\omega,s)ds
$$
extends from the unit sphere $\mS$ to a homogeneous
polynomial of degree $k$ in $\omega$.
\end{enumerate}
The seemingly ``trivial'' evenness condition is sometimes the hardest to generalize to other transforms of Radon type, while it is often easier
to find analogs of the moment conditions. This is exactly what happens in TAT.

For the operators $\W,\M$ in TAT, some sets of range conditions of the moment
type had been discovered over the years \cite{LP1,LP2,AQ,Patch}, but complete
range descriptions started to emerge only since 2006
\cite{AKQ,Finch_range,AmbKuc_rang,AFK,AKK,KuKuTAT,AL}.

Range descriptions for the more general operator $\W$ are harder to obtain
than for $\M$, and complete range descriptions are not known for even dimensions
or for the case of the variable speed of sound.

Let us address the case of the spherical mean operator $\M$ first.

\subsubsection{The range of the spherical mean operator $\M$.}\indent

The support and smoothness conditions are not hard to come up with, at least when
$S$ is a sphere. By choosing appropriate length scale, we can assume
that the sphere is of radius $1$ and centered at the origin, and that the interior
domain $\Omega$ is the unit ball $B=\{x\mid |x|=1\}$. If $f$ is smooth and supported inside $B$ (i.e. $f\in C_0^\infty(B)$), then it is clear
that the measured data satisfies the following

{\em {  Smoothness and support conditions:} }
\begin{equation}\label{E:smoothsupport}
g\in C_0^\infty(S\times[0,2]).
\end{equation}

An analog of the moment conditions for $g(y,r):=\M f$ was implicitly present in \cite{LP1,LP2,AQ} and explicitly formulated as such in \cite{Patch}:

{\em {  Moment conditions:} for any integer $k\geq 0$, the
moment
\begin{equation}\label{E:moment}
M_k(y)=\int\limits_{0}^{\infty} r^{2k+d-1} g(y,r)dr
\end{equation}
extends from $S$ to an (in general, non-homogeneous) polynomial $Q_k(x)$ of
degree at most $2k$.}

These two types of conditions happen to be incomplete, i.e. infinitely many
others exist. The Radon transform experience suggests to look for an analog of
evenness conditions. And indeed, a set of conditions called orthogonality
conditions was found in \cite{Finch_range,AmbKuc_rang,AKQ}.

{\em {  Orthogonality conditions:} Let $-\lambda_k^2$ be the eigenvalue of
the Laplace operator $\Delta$ in $B$ with zero Dirichlet conditions and
$\psi_k$ be the corresponding eigenfunctions. Then the following orthogonality
condition is satisfied:

\begin{equation}\label{E:orthog}
\int\limits_{S\times [0,2]} g(x,t)\partial_\nu \psi_\lambda
(x)j_{n/2-1}(\lambda t)t^{n-1}dxdt=0.
    \end{equation}
    Here $j_p(z)=c_pz^{-p}J_p(z)$ is the so called spherical Bessel function.
}

The range descriptions obtained in $2D$ in \cite{AmbKuc_rang} and then in
general dimension in \cite{AKQ} showed that these three types of
conditions completely describe the range of the operator $\M$ on functions
$f\in C^\infty_0(B)$. At the same time, the results of
 \cite{Finch_range,AKQ} showed that the moment conditions can be dropped in
odd dimensions. It was then discovered in \cite{AFK} that the moment
conditions can be dropped altogether in any dimension, since they follow
from the other two types of conditions:

\begin{theorem}\label{T:rangeM}\cite{AFK}
Let $S$ be the unit sphere. A function $g(y,t)$ on the cylinder $S\times \R^+$
can be represented as $\M f$ for some $f\in C^\infty_0(B)$ if an only if
it satisfied the above smoothness and support and orthogonality conditions
(\ref{E:smoothsupport}),(\ref{E:orthog}).

The statement also holds in the finite smoothness case, if one replaces the
requirements by $f\in H^s_0(B)$ and $g\in H^{s+(n-1)/2}_0(S\times [0,2])$.
\end{theorem}

The range of the forward operator $\M$ has not been described when $S$ is not
a sphere, but, say, a convex smooth closed surface. The moment and
orthogonality conditions hold for any $S$, and appropriate smoothness and
support conditions can also been formulated, at least in the convex case.
However, it has not been proven that they provide the complete range
description.

It is quite possible that for non-spherical $S$ the moment conditions might have to be included into the range description.

A different range description of the Fredholm alternative type was
developed in \cite{Palam_funk} (see also \cite{FR3} for description of this
result).

\subsubsection{The range of the forward operator $\W$.}\indent

We recall that the operator $\W$ (see (\ref{E:forward})) transforms the
initial value $f$ in (\ref{E:wave_data}) into the observed on $S$ values $g$
of the solution. There exist Kirchhoff-Poisson formulas representing
the solution $p$, and thus $g=\W f$ in terms of the spherical means of $f$ (i.e., in terms of $\M f$). However, translating the result of Theorem \ref{T:rangeM} into the
language of $\W$ is not straightforward, since in even dimensions these
formulas are non-local \cite{CH,Evans} (pp. 682 and 801 correspondingly):

\begin{equation}\label{E:KirchPois_odd}
\W f(y,t)=\frac{\sqrt{\pi}}{2\Gamma(n/2)}\left(\frac{1}{t}\frac{\partial}{\partial t}\right)^{(n-3)/2} t^{n-2}\left(\M f\right)(y,t),\mbox{ for odd } n.
\end{equation}
and
\begin{equation}\label{E:KirchPois_even}
 \W f(y,t)=\frac{1}{\Gamma(n/2)}\left(\frac{1}{t}\frac{\partial}{\partial t}
\right)^{(n-2)/2}\int\limits^t_0 \frac{r^{n-1}\left(\M
f\right)(y,r)}{\sqrt{t^2-r^2} } dr,\mbox{ for even } n.
\end{equation}

The non-locality of the transformation for even dimensions reflects the
absence of Huygens' principle (i.e. absence of sharp rear fronts of waves) in
these dimensions; it also causes difficulties in establishing the complete
range descriptions. In particular, due to the integration in (\ref{E:KirchPois_even})
$\M f(y,t)$ does not vanish for large times $t$ anymore.
One
can try to use other known operators intertwining the two problems (see
\cite{AKQ} and references therein), some of which do preserve vanishing for
large values of $t$, but this so far has lead only to very clumsy range
descriptions.

However, for odd dimensions, the range description of $\W$ can be obtained. In
order to do so, given the TAT data $g(y,t)$, let us introduce an auxiliary
time-reversed problem in the cylinder $B\times[0,2]$:
\begin{equation}\label{E:wave_aux}
    \begin{cases}
    q_{tt}-\Delta q=0 \mbox{ for }(x,t)\in B\times[0,2]),\\
    q(x,2)=q_t(x,2)=0 \mbox{ for }x\in B,\\
    q(y,t)=g(y,t)\mbox{ for }(y,t)\in S\times[0,2]).
    \end{cases}
\end{equation}
We can now formulate the range description from \cite{FR3,Finch_range}:
\begin{theorem}\label{T:rangeW_odd}\cite{FR3,Finch_range}
For odd dimensions $n$ and $S$ being the unit sphere, a function $g\in
C^\infty_0(S\times[0,2])$ can be represented as $\W f$ for some $f\in
C^\infty_0(B)$ if and only if the following condition is satisfied:

 \begin{center}
 The solution $q$ of (\ref{E:wave_aux}) satisfies $q_t(x,0)=0$ for all $x\in B$.
 \end{center}

\end{theorem}
Orthogonality type and Fredholm alternative type range conditions, equivalent
to the one in the theorem above, are also provided in \cite{FR3,Finch_range}.

\subsection{Speed of sound reconstruction}\label{S:speed}

Unsurprisingly, all inversion procedures outlined in Section \ref{S:reconstr}
rely upon the knowledge of the speed of sound $c(x)$. Although often, e.g. in
breast imaging, the medium is assumed to be acoustically homogeneous, this is
not a good assumption in many other cases. It has been observed (e.g.,
\cite{JinWang,HKN}) that replacing even slightly varying speed of sound
with its average value might significantly distort the image; not only the
numerical values, but also the shapes of interfaces between the tissues
will be reconstructed incorrectly.
Thus, the question of estimating $c(x)$ correctly becomes important.
One possible approach \cite{JinWang} is to use an additional transmission ultrasound scan to
reconstruct the speed beforehand. The question arises of whether one could determine the speed
of sound $c(x)$ and the tomogram $f(x)$ (assuming that $f$ is not zero)
simultaneously from the TAT data. In fact, one needs only to determine $c(x)$
(without knowing $f$), since then inversion procedures of Section
\ref{S:reconstr} would apply to recover $f$.

At the first glance, this seems to be an overly ambitious project.
 Indeed, if we denote the forward operator $\W$ by $\W_c$, to
indicate its dependence on the speed of sound $c(x)$, then the problem becomes, given
the data $g$, to find both $c$ and $f$ from the equality
\begin{equation}\label{E:CandF}
    \W_c f=g.
\end{equation}
A similar situation arises in the SPECT emission tomography (see
\cite{Natt_old,Natt_new,Kuch_AMS05} and references therein), where the role of the speed
of sound is played by the unknown attenuation. It is known, however, that in
SPECT the attenuation can be recovered for a ``generic'' $f$.

What is the reason for such a strange situation? It looks like for
\underline{any} $c$ one could solve the equation (\ref{E:CandF}) for an $f$,
and thus no information about $c$ is contained in the data $g$. This argument
is incorrect for the following reason: the range of the forward operator, as
we know already from the previous section, has infinite co-dimension. Thus,
this range has a lot of space to ``rotate'' when $c$ changes. Imagine for an
instance that the rotation is so powerful that for different values of $c$ the
ranges have only zero (the origin) in common. Then, knowing $g$ in the range,
one would know which $c$ it came from. Thus, the problem of recovering the
speed of sound from the TAT data is closely related to the range descriptions.

Numerical inversions using algebraic iterative techniques (e.g., \cite{Anastasio_speed,Zhang}) show that
recovering both $c$ and $f$ might be indeed possible.

Unfortunately, very little is known at the moment concerning this problem.
Direct usage of range conditions attempted in \cite{HKN} has lead only to extremely
weak and not practically useful results so far. A revealing relation to the transmission eigenvalue problem
well known in inverse problems  (see \cite{Colton_trans} for the survey) was
recently discovered by D. Finch. Unfortunately, the
transmission eigenvalue problem remains still unresolved. However, one can
derive from this relation the following (still not too useful for TAT)
uniqueness of the speed of sound determination result, due to M. Agranovsky:
\begin{theorem}\label{T:speed_trans}
If two speeds satisfy the inequality $c_1(x)\geq c_2(x)$ for all $x\in \Omega$
and produce for some functions $f_1,f_2$ the same non-zero TAT data $g$ (i.e.,
$\W c_1f_1=g,\W c_2f_2=g$), then $c_1(x)=c_2(x)$.
\end{theorem}

It is known \cite[Corollary 8.2.3]{Isakov} that if a function $f(x)$ is
such that $\Delta f(x)\neq 0$ and for two acoustic speeds $c_1(x)$ and
$c_2(x)$ it produces the same TAT data $g$, then $c_1=c_2$.


It is clear that the problem of finding the speed of sound from the TAT data
still requires significant analysis.

\section{Reconstruction formulas and procedures}\label{S:reconstr}

Numerous formulas, algorithms and procedures for reconstruction of images
from TAT measurements have been developed by now. Most of these techniques
require the data being collected on a closed surface (closed curve in $2D$)
surrounding the object to be imaged. Such methods are discussed in Section
\ref{S:closed-surf}. We review methods that work under  the assumption of
constant speed of sound in Section \ref{SS:constantspeed}. The techniques applicable in the case
of the known variable speed of sound are considered in
Section \ref{S:varspeed}. Closed surface measurements cannot always be implemented, since in some practical situations the object cannot be
completely surrounded by the detectors. In this case, one has to resort to
various approximate reconstruction techniques as discussed in Section
\ref{S:openalgs}.

\subsection{Full data (closed acquisition surfaces)\label{S:closed-surf}}

\subsubsection{Constant speed of sound}\label{SS:constantspeed}

When the speed of sound within the tissues is a known constant, the TAT
problem can be reformulated (see Section \ref{S:models}) in terms of the values of the spherical  means
of the initial condition $f(x)$. These
means can be easily recovered from the measurements of the acoustic pressure using formulas (\ref{E:KirchPois_odd}) and (\ref{E:KirchPois_even}) (see the discussion in \cite{AQ}). In this case, image
reconstruction becomes equivalent to inverting the spherical mean transform $\M$. Thus, in what follows, we consider the problem of reconstructing a function $f(x)$
supported within the region bounded by a closed surface $S$ from known values of its spherical integrals $g(y,r)$ with centers on $S$:
\begin{equation}
g(y,r)=\int\limits_{\mathbb{S}^{n-1}}f(y+r\omega)r^{n-1}d\omega,
\qquad y\in S,
\label{E:spher-int}
\end{equation}
where $d\omega$ is the standard measure on the unit sphere.
\paragraph{Series solutions for spherical geometry}\indent

The first inversion procedures for the case of closed acquisition surfaces
were described in \cite{Norton1,Norton2}, where solutions were found for the
cases of circular (in $2D$) and spherical (in $3D$) surfaces, respectively.
These solutions were obtained by the harmonic decomposition of the measured
data and of the sought function $f(x)$, followed by equating coefficients of the
corresponding Fourier series. In particular, the $2D$ algorithm
of \cite{Norton1} pertains to the case when the detectors are located on a
circle of radius $R$. This method is based on the Fourier decomposition of $f$
and $g$ in angular variables:

\begin{equation}
f(x)=\sum_{-\infty}^{\infty}f_{k}(\rho)e^{ik\varphi},\quad x=(\rho\cos
(\varphi),\rho\sin(\varphi)) \label{nortonfourf}
\end{equation}
\[
g(y(\theta),r)=\sum_{-\infty}^{\infty}g_{k}(r)e^{ik\theta},\quad
y=(R\cos(\theta),R\sin(\theta)),
\]
where
\[
\left(  \mathcal{H}_{m}u\right)  (s)=2\pi\int_{0}^{\infty}u(t)J_{m}(st)tdt,
\]
is the Hankel transform and $J_{m}(t)$ is the Bessel function. As shown in
\cite{Norton1}, the Fourier coefficients $f_{k}(\rho)$ can be recovered from
the known coefficients $g_{k}(r)$ by the following formula:
\[
f_{k}(\rho)=\mathcal{H}_{m}\left(  \frac{1}{J_{k}(\lambda|R|)}\mathcal{H}
_{0}\left[  \frac{g_{k}(r)}{2\pi r}\right]  \right)  .
\]

This method requires division of the Hankel transform of the measured data by
the Bessel functions $J_{k}$, which have infinitely many zeros. Theoretically,
there is no problem: the range conditions (Section \ref{S:range}) on the exact data $g$ imply that the Hankel transform
$\mathcal{H}_{0}\left[  (2\pi r)^{-1}g_{k}(r)\right]$ has zeros that cancel
those in the denominator. However, since the measured data always contain
errors, the exact cancelation does not happen, and one needs a sophisticated
regularization scheme to guarantee that the error remains bounded.

This difficulty can be avoided (see, e.g. \cite{KuKuTAT}) by replacing the
Bessel function $J_{0}$ in the inner Hankel transform by the Hankel function
$H_{0}^{(1)}$. This yields the following formula for $f_{k}(\rho):$
\[
f_{k}(\rho)=\mathcal{H}_{k}\left(  \frac{1}{H_{k}^{(1)}(\lambda|R|)}\int
_{0}^{\infty}g_{k}(r)H_{0}^{(1)}(\lambda r)dr\right)  .
\]
Unlike $J_{m}$, Hankel functions $H_{m}^{(1)}(t)$ do not have zeros for any
real values of $t$, which removes the problems with division by
zeros \cite{Norton1}. (A different way of avoiding divisions by zero was found
in \cite{Haltm-circ})

This derivation can be repeated in $3D$, with the exponentials $e^{ik\theta}$
replaced by the spherical harmonics, and with cylindrical Bessel functions
replaced by their spherical counterparts. By doing this, one arrives at the
Fourier series method of \cite{Norton2} (see also \cite{MXW1}). The use of the
Hankel function $H_{0}^{(1)}$ above is similar to the way the spherical Hankel
function $h_{0}^{(1)}$ is utilized in \cite{Norton2} to avoid the divisions by
zero.

\paragraph{Eigenfunction expansions for a general geometry\label{SS:series}}\indent

The series methods described in the previous section rely on the separation of
variables that occurs only in spherical geometry. A different approach was
proposed in \cite{Kun_series}. It works for arbitrary closed surfaces, but is practical only for those with explicitly known eigenvalues and eigenfunctions of the Dirichlet Laplacian in the interior. Such surfaces include, in particular, spheres, half-sp
heres, cylinders, cubes and parallelepipeds, as well as the
surfaces of crystallographic domains.

Let $\lambda_{m}^{2}$ and $u_{m}(x)$ be the eigenvalues and an ortho-normal
basis of eigenfunctions of the Dirichlet Laplacian $-\Delta$ in the interior
$\Omega$ of a closed surface $S$:
\begin{align}
\Delta u_{m}(x)+\lambda_{m}^{2}u_{m}(x)  &  =0,\qquad x\in\Omega,\quad
\Omega\subseteq\mathbb{R}^{n},\label{Helmeq}\\
u_{m}(x)  &  =0,\qquad x\in S,\nonumber\\
||u_{m}||_{2}^{2}  &  \equiv\int\limits_{\Omega}|u_{m}(x)|^{2}dx=1.\nonumber
\end{align}
As before, one would like to reconstruct a compactly supported function $f(x)$
from the known values of its spherical integrals $g(y,r)$ (see
(\ref{E:spher-int})) with centers on $S$. Since $u_{m}(x)$ is the solution of
the Dirichlet problem for the Helmholtz equation with zero boundary conditions
and the wave number $\lambda_{m}$, this function admits the Helmholtz
representation
\begin{equation}
u_{m}(x)=\int_{S}\Phi_{\lambda_{m}}(|x-y|)\frac{\partial}{\partial n}
u_{m}(y)ds(y)\qquad x\in\Omega, \label{helmdiscr}
\end{equation}
where $\Phi_{\lambda_{m}}(|x-y|)$ is a free-space
Green's function of the Helmholtz equation (\ref{Helmeq}), and $n$ is the exterior normal to $S.$

The function $f(x)$ can be expanded into the series
\begin{align}
f(x) &  =\sum_{m=0}^{\infty}\alpha_{m}u_{m}(x), \mbox{ where}\label{fourierser}\\
\alpha_{m} &  =\int_{\Omega}u_{m}(x)f(x)dx.\nonumber
\end{align}
A reconstruction formula for $\alpha_{m}$ (and thus for $f(x))$ will result, if
one substitutes representation (\ref{helmdiscr}) into (\ref{fourierser}) and
interchanges the orders of integration:
\begin{equation}
\alpha_{m}=\int_{\Omega}u_{m}(x)f(x)dx=\int_{S}I(y,\lambda_{m})\frac{\partial
}{\partial n}u_{m}(y)dA(x),\label{serkoef1a}
\end{equation}
where
\begin{equation}
I(y,\lambda)=\int_{\Omega}\Phi_{\lambda}(|x-y|)f(x)dx=\int_{0}
^{\mathrm{diam\,\Omega}}g(y,r)\Phi_{\lambda}(r)dr.\label{Greenint}
\end{equation}
Now $f(x)$ can be obtained by summing the series (\ref{fourierser}). This method becomes
computationally efficient when the eigenvalues and eigenfunctions are known
explicitly, especially if a fast summation formula for the series
(\ref{fourierser}) is available. This is the case for a cubic acquisition
surface $S$, when the eigenfunctions are products of sine functions. The resulting $3D$ reconstruction algorithm is extremely fast and precise (see \cite{Kun_series}).

The above method has an interesting property. If the support of the source
$f(x)$ extends outside $\Omega,$ the algorithm still yields theoretically
exact reconstruction of $f(x)$ inside $\Omega$. Indeed, the value of the
expression (\ref{helmdiscr}) for all $x$ lying outside $\Omega$ is zero. Thus,
when one computes (\ref{serkoef1a}) for $x\in\mathbb{R}^{n}\setminus\Omega$, values of $f(x)$
are multiplied by zero and do not affect further computation in any way. This feature is shared by the time reversal method (see the corresponding paragraph in Section \ref{S:varspeed}). The closed form FBP type
reconstruction techniques considered in the next sub-section, do not have this
property. In other words, in presence of a source outside the measurement surface, reconstruction within $\Omega$ will be incorrect.

The reason for this difference is that all currently known closed form FBP-type formulas rely (implicitly or explicitly) on the assumption that the wave propagates outside $S$ in the whole free space and has no sources outside. On the other hand, the eig
enfunction expansion
method and the time reversal rely only upon the time decay of the wave inside $S$, which is not influenced by $f$
having a part outside $S$.

\paragraph{Closed form inversion formulas}\indent

Closed-form inversion formulas play a special role in tomography. They bring
about better theoretical understanding of the problem and frequently serve as
starting points for the development of efficient reconstruction algorithms. A
well known example of the use of explicit inversion formulas is the so-called
filtered backprojection (FBP) algorithm in X-ray tomography, which is derived
from one of the inversion formulas for the classical Radon transform (see, for
example \cite{Kak,Natt_old}).

The very existence of closed form inversion formulas for TAT had been in
doubt, till the first such formulas were obtained in odd dimensions by Finch
et al in \cite{FPR}, under the assumption that the acquisition surface $S$ is
a sphere. Suppose that the function $f(x)$ is supported within a ball of
radius $R$ and that the detectors are located on the surface $S=\partial B$ of
this ball. Then some of the formulas obtained in \cite{FPR} read as follows:
\begin{align}
f(x)  &  =-\frac{1}{8\pi^{2}R}\Delta_{x}\int\limits_{\partial B}
\frac{g(y,|y-x|)}{|y-x|}dA(y),\label{FPR3da}\\
f(x)  &  =-\frac{1}{8\pi^{2}R}\int\limits_{\partial B}\left(  \frac{1}{r}
\frac{\partial^{2}}{\partial r^{2}}g(y,r)\right)  \left.
{\phantom{\rule{1pt}{8mm}}}\right\vert _{r=|y-x|}dA(y),\label{FPR3d}\\
f(x)  &  =-\frac{1}{8\pi^{2}R}\int\limits_{\partial B}\left(  \frac{1}{r}
\frac{\partial}{\partial r}\left(  r\frac{\partial}{\partial r}\frac
{g(y,r)}{r}\right)  \right)  \left.  \phantom{\rule{1pt}{8mm} }\right\vert
_{r=|y-x|}dA(y), \label{FPR3db}
\end{align}
where $dA(y)$ is the surface measure on $\partial B$ and $g$ represents the
values of the spherical integrals (\ref{E:spher-int}).

These formulas have a FBP (filtered back-projection) nature. Indeed,
differentiation with respect to $r$ in (\ref{FPR3d}) and (\ref{FPR3db}) and
the Laplace operator in (\ref{FPR3da}) represent the filtration, while the
(weighted) integrals correspond to the backprojection, i.e. integration over
the set of spheres passing through the point of interest $x$ and centered on
$S$.

The so-called ``universal backprojection formula'' in $3D$ was found in
\cite{MXW2} (it is also valid for the cylindrical and plane acquisition
surfaces, see Section \ref{S:openalgs}). In our notation, this formula takes
the form
\begin{equation}
f(x)=\frac{1}{8\pi^{2}}\mathrm{div}\int\limits_{\partial B}n(y)\left(
\frac{1}{r}\frac{\partial}{\partial r}\frac{g(y,r)}{r}\right)  \left.
{\phantom {\rule {1pt}{8mm}}}\right\vert _{r=|y-x|}dA(y), \label{E:universal}
\end{equation}
or, equivalently,
\begin{equation}
f(x)=-\frac{1}{8\pi^{2}}\int\limits_{\partial B}\frac{\partial}{\partial
n}\left(  \frac{1}{r}\frac{\partial}{\partial r}\frac{g(y,r)}{r}\right)
\left.  {\phantom {\rule {1pt}{8mm}}}\right\vert _{r=|y-x|}dA(y),
\label{E:universal1}
\end{equation}
where $n(y)$ is the exterior normal vector to
$\partial B$. One can show \cite{MXW2,AKK,nguyen} that formulas (\ref{FPR3da})
through (\ref{E:universal}) are not equivalent on non-perfect data: the result
will differ if these formulas are applied to a function that does not belong
to the range of the spherical mean transform $\M$. A family of inversion
formulas valid in $\mathbb{R}^{n}$ for arbitrary $n\geq2$ was found in
\cite{Kunyansky}:
\begin{equation}
f(x)=\frac{1}{4(2\pi)^{n-1}}\mathrm{div}\int\limits_{\partial B}
n(y)h(y,|x-y|)dA(y), \label{E:kunyansky}
\end{equation}
where
\begin{align}
h(y,t)  &  =\int\limits_{\mathbb{R}^{+}}Y(\lambda t)\left[  \int
\limits_{0}^{2R}J(\lambda r)g(y,r)dr-J(\lambda t)\int\limits_{0}^{2R}Y(\lambda
r)g(y,r)dr\right]  \lambda^{2n-3}d\lambda,\label{genfilt}\\
J(t)  &  =\frac{J_{n/2-1}(t)}{t^{n/2-1}},\quad\quad Y(t)=\frac{Y_{n/2-1}
(t)}{t^{n/2-1}}, \label{E:specfun}
\end{align}
and $J_{n/2-1}(t)$ and $Y_{n/2-1}(t)$ are respectively the Bessel and
Neumann functions of order $n/2-1$. In $3D$, $J(t)$ and $Y(t)$ are simply
$t^{-1} \sin t$ and $t^{-1} \cos t $ and formulas (\ref{E:kunyansky}) and
(\ref{genfilt}) reduce to (\ref{E:universal1}).

In $2D$, equation (\ref{genfilt}) also can be simplified \cite{AKK}, which
results in the formula
\begin{equation}
f(x)=\frac{1}{2\pi^{2}}\mathrm{div}\int\limits_{\partial B}n(y)\left[
\int\limits_{0}^{2R}g(y,r)\frac{1}{r^{2}-|x-y|^{2}}dr\right]  dl(y),
\label{E:kun2d}
\end{equation}
where $\partial B$ now stands for the circle of radius $R$ and $dl(y)$ is the
standard arc length.

A different set of closed-form inversion formulas applicable in even dimensions
was found in \cite{Finch_even}. Formula (\ref{E:kun2d}) can be compared to the
following inversion formulas from  \cite{Finch_even}:
\begin{equation}
f(x)=\frac{1}{2\pi R}\Delta\int\limits_{\partial B}\int\limits_{0}
^{2R}g(y,r)\log(r^{2}-|x-y|^{2})\ dr\ dl(y), \label{E:Finch2d}
\end{equation}
or
\begin{equation}
f(x)=\frac{1}{2\pi R}\int\limits_{\partial B}\int\limits_{0}^{2R}
\frac{\partial}{\partial r}\left(  r\frac{\partial}{\partial r}\frac
{g(y,r)}{r}\right)  \log(r^{2}-|x-y|^{2})\ dr\ dl(y). \label{E:Finch2d1}
\end{equation}

Finally, a unified family of inversion formulas was derived in \cite{nguyen}.
In our notation, it has the following form:
\begin{align}
f(x)  &  =-\frac{4}{\pi R}\int\limits_{\partial B}\left(  \frac{\partial
}{\partial t}K_{n}(y,t)\right)  \left.  {\phantom {\rule {1pt}{8mm}}}
\right\vert _{t=|x-y|}\frac{<y-x,y-\xi>}{|x-y|}dA(y),\label{E:linhNd}\\
K_{n}(y,t)  &  =-\frac{1}{16(2\pi)^{n-2}}\int\limits_{\mathbb{R}^{+}}
\lambda^{2n-3}Y(\lambda t)\left(  \int\limits_{\mathbb{R}^{+}}J(\lambda
r)g(y,r)dr\right)  d\lambda\nonumber
\end{align}
where $\partial B$ is the surface of a ball in $\mathbb{R}^{n}$ of radius $R,$
functions $J$ and $Y$ are as in (\ref{E:specfun}),
and $\xi$ is an arbitrary fixed vector. In particular, in $3D$
\[
J(t)=\sqrt{\frac{2}{\pi}}\frac{\sin t}{t},J(t)=\sqrt{\frac{2}{\pi}}\frac{\cos
t}{t}
\]
and, after simple calculation, the above inversion formula reduces to
\begin{equation}
f(x)=-\frac{1}{8\pi^{2}R}\int\limits_{\partial B}\left(  \frac{\partial
}{\partial r}\frac{1}{r}\frac{\partial}{\partial r}\frac{g(y,r)}{r}\right)
\left.  {\phantom {\rule {1pt}{8mm}}}\right\vert _{r=|x-y|}\frac{<y-x,y-\xi
>}{|x-y|}dA(y). \label{E:linh3d}
\end{equation}
Different choices of vector $\xi$ in the above formula result in different
inversion formulas. For example, if $\xi$ is set to zero, the ratio
$\frac{<y-x,y-\xi>}{|x-y|}$ equals $R\cos\alpha,$ where $\alpha$ is the angle
between the exterior normal $n(y)$ and the vector $y-x;$ when combined with
the derivative in $t$ this factor produces the normal derivative, and the
inversion formula (\ref{E:linh3d}) reduces to (\ref{E:universal1}). On the
other hand, the choice of $\xi=x$ in (\ref{E:linh3d}) leads to a formula
\[
f(x)=-\frac{1}{8\pi^{2}R}\int\limits_{\partial B}\left(  r\frac{\partial
}{\partial r}\frac{1}{r}\frac{\partial}{\partial r}\frac{g(y,r)}{r}\right)
\left.  {\phantom {\rule {1pt}{8mm}}}\right\vert _{r=|x-y|}dA(y),
\]
which is reminiscent of formulas (\ref{FPR3da})-(\ref{FPR3db}).

\paragraph{Greens' formula approach and some symmetry considerations}\indent

Let us suppose for a moment that the acoustic detectors
could measure not only the pressure $p(y,t)$ at each point of the
acquisition surface $S$, but also the normal derivative $\partial p/\partial
n$ on $S$. Then the problem of reconstructing the initial pressure $f(x)$
becomes rather simple. Indeed, one can use the knowledge of the free-space
Green's function for the wave equation and  invoke the Green's theorem to
represent the solution $p(x,t)$ of (\ref{E:wave_data}) in the form of
integrals over $S$ involving $p(x,t)$ and its normal derivative and the
Green's function and its normal derivative. (This can be done in the Fourier
or time domains.) This would require infinite observation time, but in $3D$
the time $T(\Omega)$ will suffice, afte r which the wave escapes the region
of interest (a cut-off also would work approximately in $2D$. similarly to
the time-reversal method). This Green's function approach happens to be,
explicitly or implicitly, the starting point of all closed form inversions
described above. The trick is to rewrite the formula in such a way that the
unknown in reality normal derivative  $\partial p/\partial n$ disappears
from the formula.

This was achieved in \cite{Kunyansky} by reducing the question to some integrals involving special functions and making the key observation that the integral
\[
I_{\lambda}(x,y)=\int\limits_{\partial B}J(\lambda|x-z|)\frac{\partial}{\partial
n}Y(\lambda|y-z|)dA(z),\qquad x,y\in B\subset\mathbb{R}^{n}
\]
is a symmetric function of its arguments:
\begin{equation}
I_{\lambda}(x,y)=I_{\lambda}(y,x)\mbox{ for }x,y\in B\subset\mathbb{R}
.^{n}\label{E:my-symm}
\end{equation}
Similarly, the derivation of (\ref{E:linhNd}) in \cite{nguyen} employs the symmetry of the
integral
\[
K_{\lambda}(x,y)=\int\limits_{\partial B}J(\lambda|x-z|)Y(\lambda|y-z|)dA(z),\qquad
x,y\in B\subset\mathbb{R}^{n}.
\]

In fact, the symmetry holds for any integral
\[
W_{\lambda}(x,y)=\int\limits_{\partial B}U(\lambda|x-z|)V(\lambda|y-z|)dA(z),\qquad
x,y\in B\subset\mathbb{R}^{n},
\]
where $U(\lambda|x|)$ and $V(\lambda|x|)$ are any two radial solutions of Helmholtz equation

\begin{equation}
\Delta u(x)+\lambda^{2}u(x)=0.\label{E:helmh}
\end{equation}

It is straightforward to verify this symmetry when $S$ is a sphere and $B$ is the corresponding ball, and the points $x,y$ lie on the boundary $S$ only, rather than anywhere in $B$. This follows immediately from the rotational symmetry of $S$. The same i
s true for the normal derivatives on $S$ of $W_{\lambda}(x,y)$ in $x$ and $y$.

This boundary symmetry happens to imply the needed full symmetry (\ref{E:my-symm}) for $x,y\in B$.

Indeed, $W_{\lambda}(x,y)$
is a solution of the Helmholtz equation
separately as a function of $x$ and of $y.$ Let us introduce a family of
solutions $\{w_{n}(x)\}_{n=0}^{\infty}$ of (\ref{E:helmh})
in $B$, such that the members of this family form an orthonormal basis for all
solutions of the latter equation in $B$.
For example, the spherical waves, i.e. the
products of spherical harmonics and Bessel functions, can serve as such a basis.

Then $W_{\lambda}(x,y)$ can be
expanded n the following series:
\begin{equation}
W_{\lambda}(x,y)=\sum_{n=0}^{\infty}\sum_{m=0}^{\infty}b_{n,m}w_{m}
(y)w_{n}(x).\label{E:doubleser}
\end{equation}
Since $W_{\lambda}(x,y)$ is a solution to the Helmholtz
equation in $\partial B\times\partial B,$ coefficients $b_{n,m}$ are
completely determined by the boundary values of $W_{\lambda}$. Since the
boundary values are symmetric, the coefficients are symmetric, i.e.
$b_{n,m}=b_{m,n}$ which by (\ref{E:doubleser}) immediately implies
$W_{\lambda}(x,y)=W_{\lambda}(y,x)$ for all pairs $(x,y)\in B\times B$.

This consideration extends to infinite cylinders
and planes. This
explains why the ``universal backprojection formula''
(\ref{E:universal1}) is valid also for infinite cylinders and planes
\cite{MXW2}. Since the sort of symmetry used is shared only by these three
surfaces, we believe it is unlikely that a closed-form formula could exist
for any other acquisition surface.

\paragraph{Algebraic iterative algorithms}\indent

Iterative algebraic techniques are among the favorite tomographic methods of
reconstruction and have been used in CT for quite a while
\cite{Natt_new,Natt_old,Kak}. They amount to discretizing the equation
relating the measured data with the unknown source, followed by iterative
solution of the resulting linear system. Iterative algebraic reconstruction
algorithms frequently produce better images than those obtained by other
methods. However, they are notoriously slow. In TAT, they have been used
successfully for reconstructions with partial data
(\cite{Paltaufiter,Anastasio_halftime,Anastasio_half}), see Section
\ref{S:openalgs}.

\paragraph{Parametrix approaches}\indent

Some of the earlier non-iterative reconstruction techniques \cite{Kruger} were
of approximate nature. For example, by approximating the integration spheres
by their tangent planes at the point of reconstruction and by applying one of
the known inversion formulas for the classical Radon transform, one can
reconstruct an approximation to the image. Due to the evenness symmetry in the
classical Radon projections (see Section \ref{S:range}), the normals to the
integration planes need only fill a half of a unit sphere, in order to make
possible the reconstruction from an open measurement surface. A more
sophisticated approach is represented by the so-called ``straightening''
methods \cite{PopSush,PopSush2} based on the approximate reconstruction of the
classical Radon projections from the values of the spherical mean transform
$\M f$ of the function $f(x)$ in question. These methods yield not a true
inversion, but rather what is called in micro-local analysis a
\textbf{parametrix}. Application of a parametrix reproduces the function $f$
with an additional, smoother term.
In other words, the locations (and often the sizes)
of jumps across sharp material interfaces, as well as the whole wave front set
$WF(f)$, are reconstructed correctly, while the accuracy of the lower spatial
frequencies cannot be guaranteed.
(Sometimes, the reconstructed function has a more general form $Af$, where $A$ is an elliptic pseudo-differential
operator \cite{Shubin,Horm,Str} of order zero. In this case, the sizes of the jumps across the interfaces might be altered.) Unlike the approximations resulting from the
discretization of the exact inversion formulas (in the situations when such
formulas are known), the parametrix approximations do not converge, when the
discretization of the data is refined and the noise is eliminated. Parametrix
reconstructions can be either accepted as approximate images, or used as starting points for iterative algorithms. See \cite{StefUhlTAT} for a recent discussion of parametrices.

These methods are closely related to the general scheme proposed in
\cite{Be,Beylkin} for the inversion of the generalized Radon transform with
integration over curved manifolds. It reduces the problem to a Fredholm
integral equation of the second kind, which is well suited for numerical
solution. Such an approach amounts to using a parametrix method as an
efficient pre-conditioner for an iterative solver; the convergence of such
iterations is much faster than that of algebraic iterative methods.

\paragraph{Numerical implementation and computational examples.}\indent

By discretizing exact formulas presented above, one can easily develop accurate
and efficient reconstruction algorithms. The $3D$ case is especially simple:
computation of derivatives in the formulas (\ref{FPR3da})-(\ref{E:universal1})
and (\ref{E:linh3d}) can be easily done, for instance by using finite differences; it is followed
by the backprojection (described by the integral over $\partial B$),
which requires prescribing quadrature weights for quadrature nodes that
coincide with the positions of the detectors. The backprojection step is
stable; the differentiation is a mildly unstable operation. The sensitivity to
noise in measurements across the formulas presented above seems to be roughly
the same. It is very similar to that of the widely used FBP algorithm of
classical X-ray tomography \cite{Natt_old,Natt_new}. In $2D$, the
implementation is just a little bit harder: the filtration step in formulas
(\ref{E:kun2d})-(\ref{E:Finch2d1}) can be reduced to computing two Hilbert
transforms (see \cite{KuKuTAT}), which, in turn, can be easily done in the
frequency domain.

The number of floating point operations (flops) required by such algorithms is
determined by the slower backprojection step. In $3D$, if the number of
detectors is $m^{2}$ and the size of the reconstruction grid is $m\times
m\times m$, the backprojection step (and the whole algorithm) will require
$O(m^{5})$ flops. In practical terms this amounts to several hours of
computations on a single processor computer for a grid of size $129\times
129\times129.$

In $2D$, the operation count is just $O(m^{3})$. As it is discussed in Section \ref{S:integrating}, the $2D$ problem needs to be solved, when integrating line detectors are used. In this situation, the $2D$ problem needs to be solved $m$ times in order t
o reconstruct the image, which
raises the total operation count to $O(m^{4})$ flops.

\begin{figure}[ht!]
\par
\begin{center}
\begin{tabular}
[c]{cc}
\includegraphics[width=1.5in,height=1.5in]{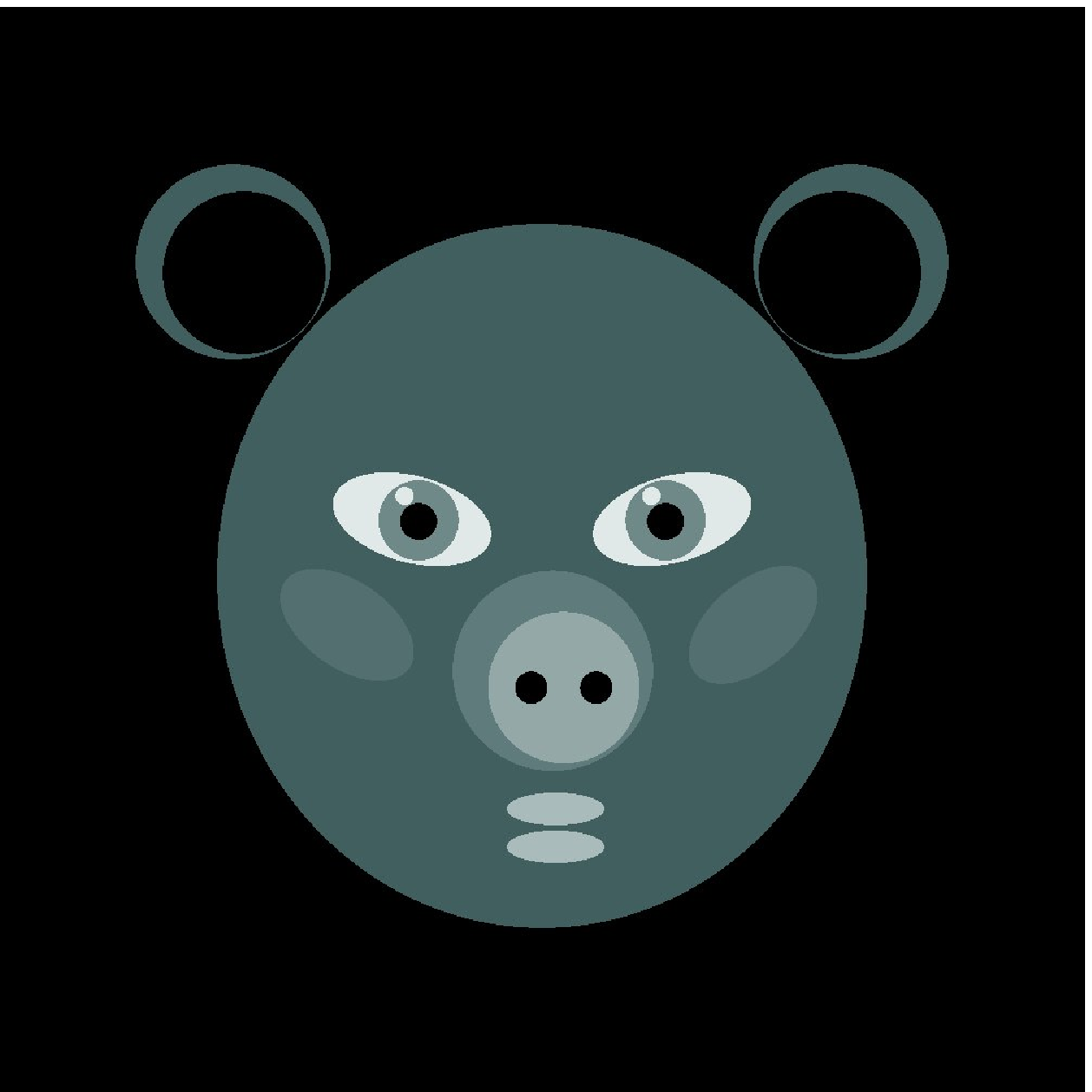} &
\includegraphics[width=1.5in,height=1.5in]{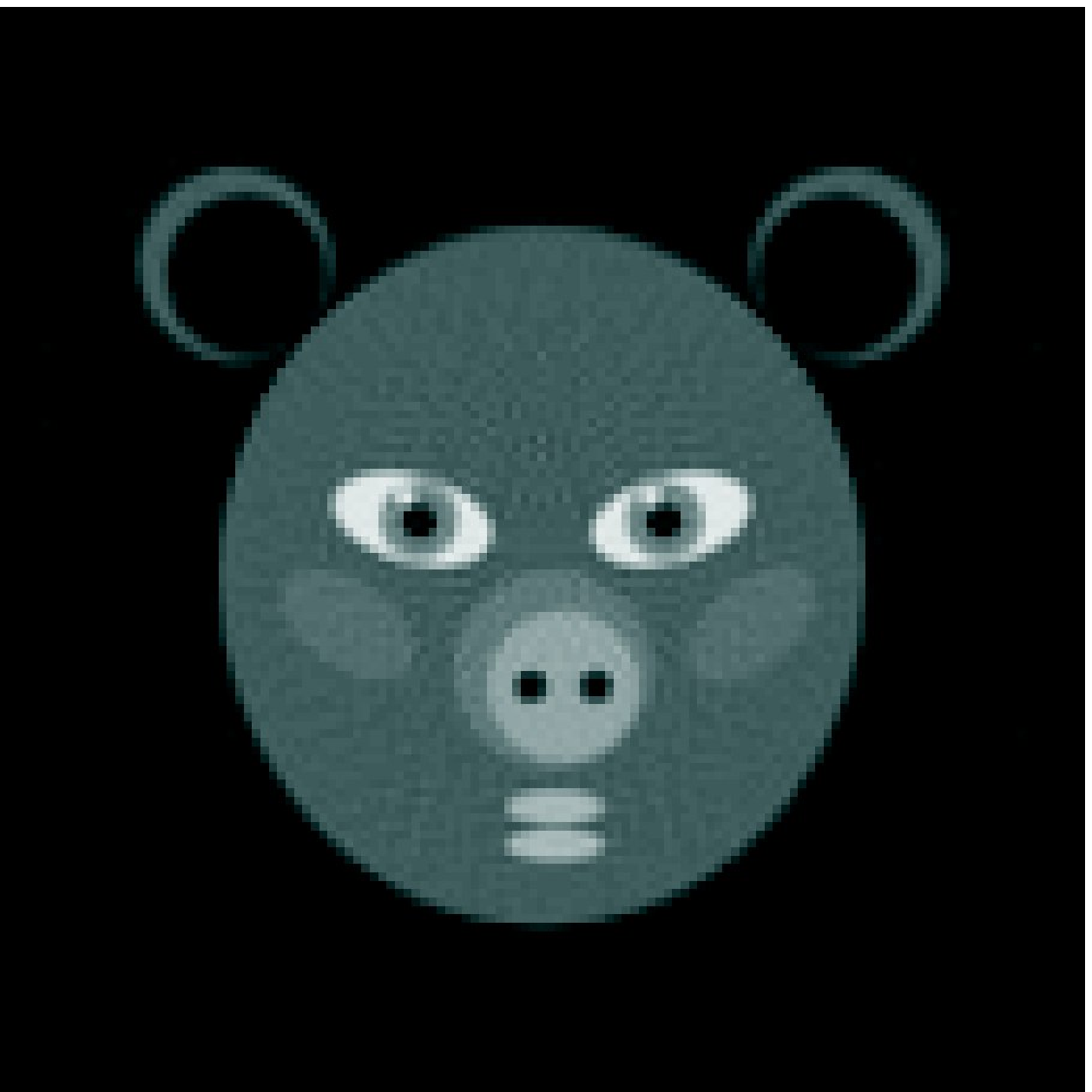}\\
(a) & (b)\\
& \\
\includegraphics[width=1.5in,height=1.5in]{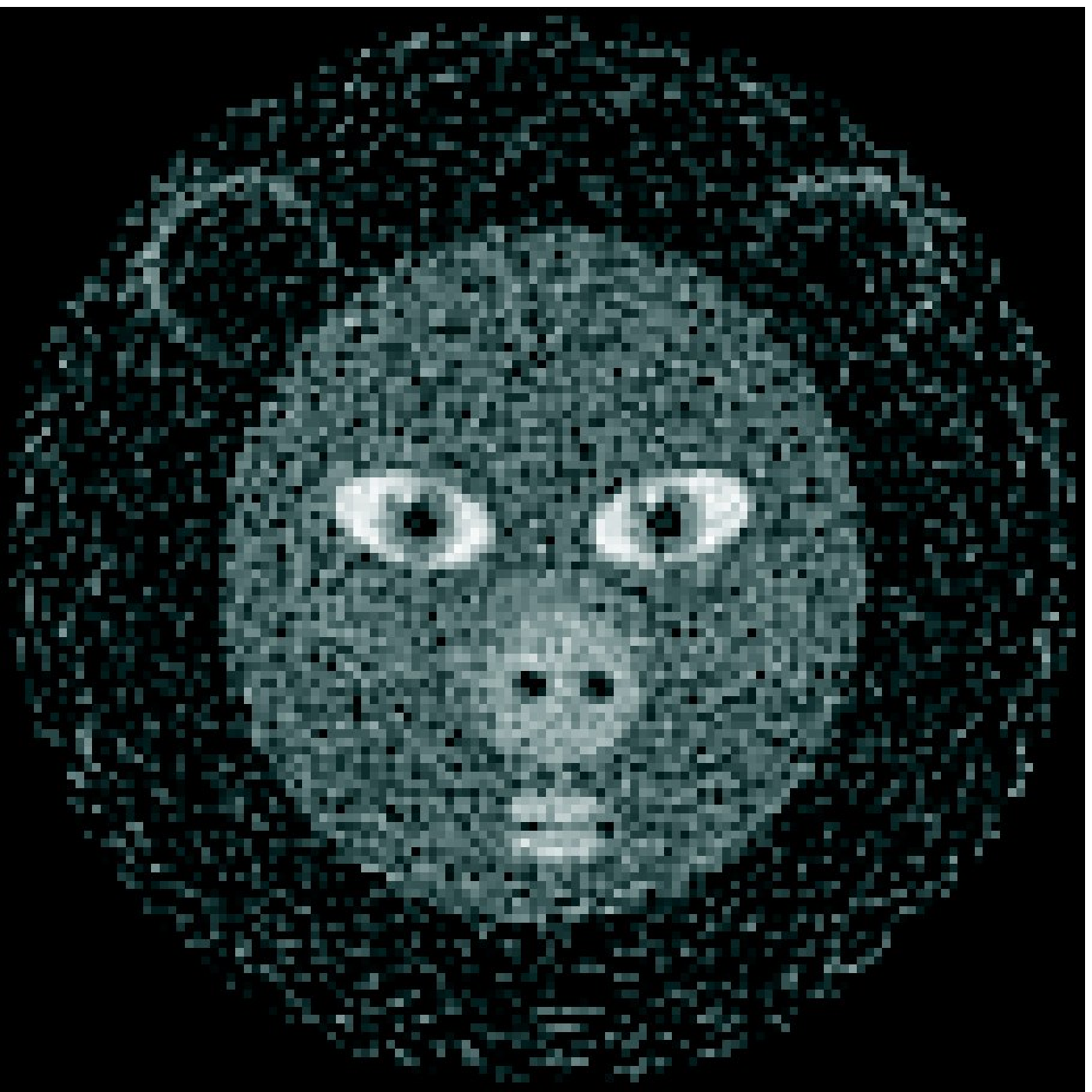} &
\includegraphics[width=1.5in,height=1.5in]{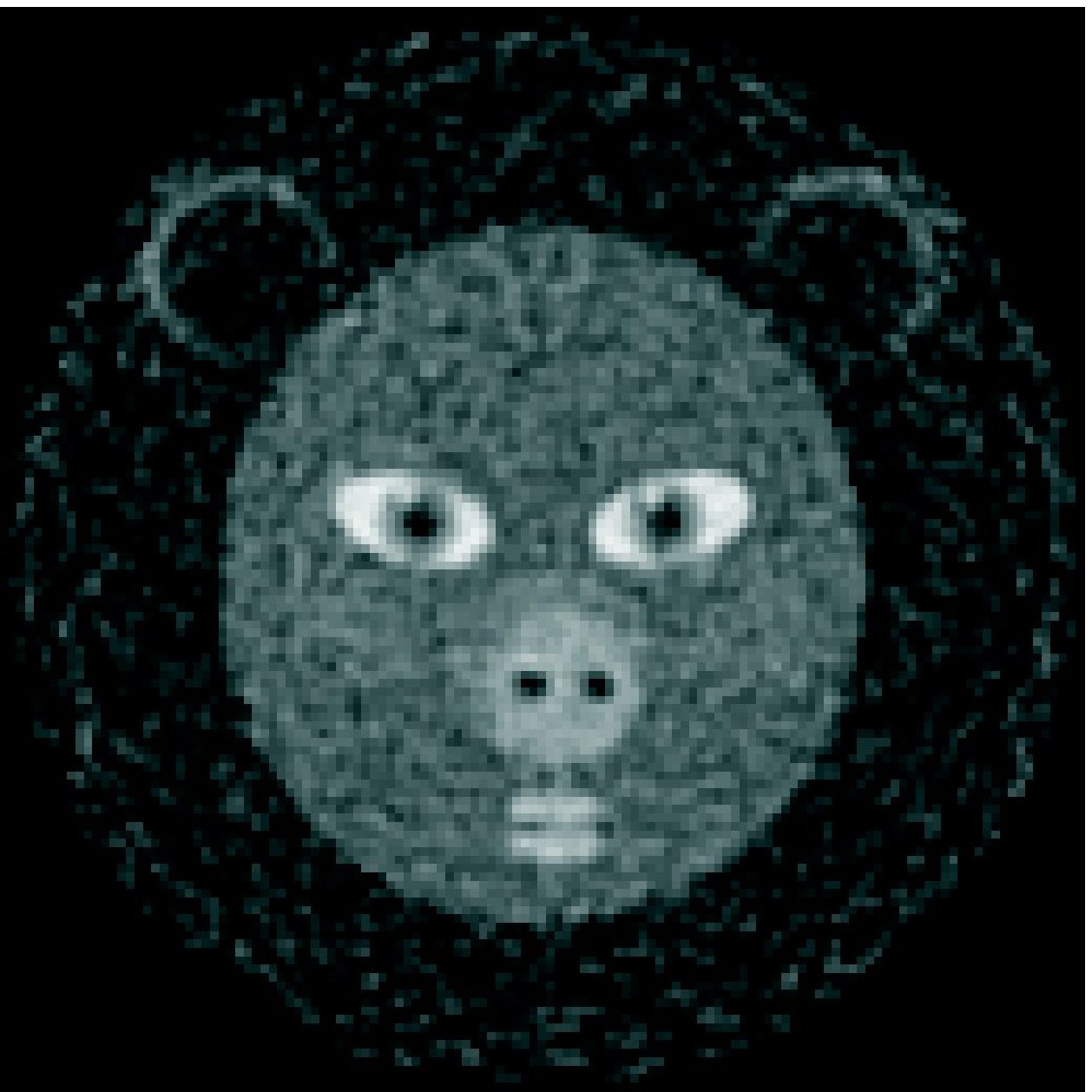}\\
(c) & (d)
\end{tabular}
\end{center}
\caption{Example of a reconstruction using formula (\ref{E:kun2d}): (a)
phantom; (b) reconstruction from accurate data; (c) reconstruction from the
data contaminated with 15\% noise; (d) reconstruction from the noisy data with
additional smoothing}
\label{F:2d}
\end{figure}

Figure \ref{F:2d} shows three examples of simulated reconstruction using
formula (\ref{E:kun2d}). The phantom we use (Figure \ref{F:2d}(a)) is a linear
combination of several characteristic functions of disks and ellipses. Part
(b) illustrates the image reconstruction within the unit circle
from 257 equi-spaced projections each containing 129 spherical integrals. The detectors were placed on the concentric circle of radius 1.05. The image
shown in Figure \ref{F:2d}(c) corresponds to the reconstruction from the
simulated noisy data that were obtained by adding to projections\ values of a
random variable scaled so that the $L^{2}$ intensity of the noise was 15\% of
the intensity of the signal. Finally, Figure \ref{F:2d}(d) shows how
application of a smoothing filter (in the frequency domain)\ suppresses the
noise; it also somewhat blurs the edges in the image.

\subsubsection{Variable speed of sound\label{S:varspeed}}

The reconstruction formulas and algorithms described in the previous
section work under the assumption that the speed of sound within the region of
interest is constant (or at least close to a constant). This assumption,
however, is not always realistic -- for example, if the region of interest
contains both soft tissues and bones, the speed of sound will vary
significantly. Experiments with numerical and physical phantoms
show \cite{HKN,JinWang} that if acoustic inhomogeneities are not taken into
account, the reconstructed image might be severely distorted. Not only the
numerical values could be reconstructed incorrectly, but so would the material
interface locations and discontinuity magnitudes.

Below we review some of the reconstruction methods that work in acoustically
inhomogeneous media. We will assume that the speed of sound $c(x)$ is known,
smooth, positive, constant for large $x$, and non-trapping. In practice, a
transmission ultrasound scan can be used to reconstruct $c(x)$ prior to
thermoacoustic reconstruction, as it is done in \cite{JinWang}.

\paragraph{Time reversal}\indent

Let us assume temporarily that the speed of sound $c$ is constant and the
spatial dimension is odd. Then Huygens' principle guarantees that the sound
wave will leave the region of interest $\Omega$ in time
$T=c/(\mathrm{diam\,}\Omega),$ so that $p(x,t)=0$ for all $x\in\Omega$ and
$t\geq T$. Now one can solve the wave equation back in time from $t=T$ to
$t=0$ in the domain $\Omega\times\lbrack T,0]$, with zero initial conditions
at $T$ and boundary conditions on $S$ provided by the data $g$ collected by
the detectors. Then the value of the solution at $t=0$ will coincide with the
initial condition $f(x)$ that one seeks to reconstruct. Such a solution of the
wave equation is easily obtained numerically by finite difference techniques
\cite{Grun,HKN}. The required number of floating point operations is actually
lower than that of methods based on discretized inversion formulas
($\mathcal{O}(m^{4})$ for time reversal on a grid $m\times m\times m$ in $3D$
versus $\mathcal{O}(m^{5})$ for inversion formulas), which makes this method
quite competitive even in the case of constant speed of sound.

Most importantly, however, the method is also applicable if the speed of sound
$c(x)$ is variable and/or the spatial dimension is even. In these cases, the
Huygens' principle does not hold, and thus the solution to the direct problem
will not vanish within $\partial \Omega$ in finite time. However, the solution
inside $\Omega$ will decay with time. Under the non-trapping condition, as it
is shown in (\ref{E:decay}) (see \cite{Egorov,Vainb,Vainb2}), the time decay
is exponential in odd dimensions, but only algebraic in even-dimensions.
Although, in order to obtain theoretically exact reconstruction, one would
have to start the time reversal at $T=\infty$, numerical experiments (e.g.,
\cite{HKN}) and theoretical estimates \cite{Hristova} show that in practice
it is sufficient to start at the values of $T$ when the signal becomes small
enough, and to approximate the unknown value of $p(x,T)$ by zero (a more
sophisticated cut-off is used in \cite{StefUhlTAT}, which leads to an equation with a contraction operator). This works
\cite{Grun,HKN} even in $2D$ (where decay is the slowest) and in
inhomogeneous media. However, when trapping occurs, the "invisible" parts blur
away (see Section \ref{S:incomplete} for the discussion).

\paragraph{Eigenfunction expansions.}\indent
\label{invers_variable}

An ``inversion formula'' that reconstructs the initial value $f(x)$ of
the solution of the wave equation from values on the measuring surface $S$
can be easily obtained using time reversal and Duhamel's
principle \cite{AK}. Consider in $\Omega$ the operator
$A=-c^{2}(x)\Delta$ with zero Dirichlet conditions on the boundary
$S=\partial\Omega$. This operator is self-adjoint, if considered in the
weighted space $L^{2}(\Omega;c^{-2}(x))$. Let us denote by $E$ the operator of
harmonic extension, which transforms a function $\phi$ on $S$ to a harmonic
function on $\Omega$ which coincides with $\phi$ on $S$. Then $f$ can be
reconstructed \cite{AK} from the data $g$ in (\ref{E:wave_data})
by the following formula:
\begin{equation}
f(x)=(Eg|_{t=0})-\int\limits_{0}^{\infty}A^{-\frac{1}{2}}\sin{(\tau
A^{\frac{1}{2}})}E(g_{tt})(x,\tau)d\tau, \label{E:reconstruction_variable}%
\end{equation}
which is valid under the non-trapping condition on $c(x)$.
However, due to the involvement of functions of the operator $A$, it is not clear how useful this formula can be.

One natural way to try to implement numerically the formula
(\ref{E:reconstruction_variable}) is to use the eigenfunction expansion of the
operator $A$ in $\Omega$ (assuming that such expansion is known). This quickly
leads to the following procedure \cite{AK}. The function $f(x)$ can be
reconstructed inside $\Omega$ from the data $g$ in (\ref{E:wave_data}), as the
following $L^{2}(B)$-convergent series:
\begin{equation}
f(x)=\sum\limits_{k}f_{k}\psi_{k}(x), \label{E:coef_variable}
\end{equation}
where the Fourier coefficients $f_{k}$ can be recovered from the data using one of the following formulas:
\begin{equation}
\begin{array}
[c]{c}
f_{k}=\lambda_{k}^{-2}g_{k}(0)-\lambda_{k}^{-3}\int\limits_{0}^{\infty}
\sin{(\lambda_{k}t)}g_{k}^{\prime\prime}(t)dt,\\
f_{k}=\lambda_{k}^{-2}g_{k}(0)+\lambda_{k}^{-2}\int\limits_{0}^{\infty}
\cos{(\lambda_{k}t)}g_{k}^{\prime}(t)dt,\mbox{ or }\\
f_{k}=-\lambda_{k}^{-1}\int\limits_{0}^{\infty}\sin{(\lambda_{k}t)}
g_{k}(t)dt=-\lambda_{k}^{-1}\int\limits_{0}^{\infty}\int\limits_{S}
\sin{(\lambda_{k}t)}g(x,t)\overline{\frac{\partial\psi_{k}}{\partial n}
(x)}dxdt,
\end{array}
\label{E:coef_variable2}
\end{equation}
where
\[
g_{k}(t)=\int\limits_{S}g(x,t)\overline{\frac{\partial\psi_{k}}{\partial
n}(x)}dx.
\]

One notices that this is a generalization of the expansion method of
\cite{Kun_series} discussed in Section \ref{SS:series} to the case of a
variable speed of sound. Unlike the algorithm of \cite{Kun_series}, the
present method does not require the knowledge of the whole space Green's
function for $A$ (which is in this case unknown). However, computation of a
large set of eigenfunctions and eigenvalues followed by the summation of the
series (\ref{E:coef_variable}) at the nodes of the computational grid may
prove to be too time consuming.

It is worthwhile to mention again that the non-trapping condition is crucial
for the stability of any TAT reconstruction method in acoustically inhomogeneous
media. As it was discussed in Section \ref{S:vis_discuss}, trapping can
significantly reduce the quality of reconstruction. It is, however, most
probable that trapping does not occur much in biological objects.

\subsection{Partial (incomplete) data\label{S:openalgs}}

Reconstruction formulas and algorithms of the previous sections work
under the assumption that the acoustic signal is measured by detectors
covering a closed surface $S$ that surrounds completely the object of interest.
However, in many practical applications of TAT, detectors can be placed
only on a certain part of the surrounding surface. Such is the case, for
example, when TAT is used for breast screening -
one of the most promising applications of this modality. Thus, one needs methods and algorithms capable of accurate reconstruction of images from partial (incomplete) data, i.e. from the measurements made on open surfaces (or open curves in $2D$).

Most exact inversion formulas and
methods discussed above are based (explicitly or implicitly) on some sort of
the Green's formula, Helmholtz representation, or eigenfunction decomposition
for closed surfaces, and thus they cannot be extended to the case of partial
data. The methods that do work in this situation rely on approximation
techniques, as discussed below.

\subsubsection{Constant speed of sound}\label{S:openalgs_const}

Even the case of an acoustically homogeneous medium is quite challenging when
reconstruction needs to be done from partial data (i.e., when the acquisition
surface $S$ is not closed). As it was discussed in Section \ref{S:incomplete},
if the detectors located around the object in such a way that the
``visibility'' condition is not satisfied, accurate reconstruction is
impossible: the ``invisible'' interfaces will be smoothed out in the
reconstructed image. On the other hand, if the visibility condition is
satisfied, the reconstruction is only mildly unstable (similarly to the
inversion of the classic Radon transform) \cite{Palam_funk,StefUhlTAT}. If, in
addition, the uniqueness of reconstruction from partial data is guaranteed (which
is usually the case, see Section \ref{S:limited_unique}), one can hope to be
able to develop an algorithm that would reconstruct quality images.

Special cases of open acquisition surfaces are a plane or an infinite
cylinder, for which exact inversion formulas are known (see, for example,
\cite{Faw,XFW,And,Nils,GGG} for the plane and \cite{XXW} or for a cylinder).
Of course, the plane or a cylinder would have to be truncated in any
practical measurements. The resulting acquisition geometry will not satisfy
the visibility condition, and material interfaces whose normals do not
intersect the acquisition surface will be blurred.

Iterative algebraic techniques (see the corresponding paragraph in Section
\ref{SS:constantspeed}) were among the first methods successfully used for
reconstruction from surfaces only partially surrounding the object (e.g.,
\cite{Paltaufiter,Anastasio_halftime,Anastasio_half}). As it is mentioned in
Section \ref{SS:constantspeed}, such methods are very slow. For example,
reconstructions in \cite{Anastasio_halftime} required the use of a cluster
of computers and took 100 iterations to converge.

Parametrix type reconstructions in the partial data case were proposed in
\cite{PopSuch_half}. A couple of different parametrix-type algorithms were
proposed in \cite{Paltaufnew1,Paltaufnew}. They are based on applying one of
the exact inversion formulas for full circular acquisition to the available
partial data, with zero-filled missing data and some correction factors.
Namely, since the missing data is replaced by zeros, each line passing
through a node of the reconstruction grid will be tangent either to one or
to two circles of integration. Therefore some directions during the
backprojection step will be represented twice, and some only once. This, in
turn, will cause some interfaces to appear twice stronger then they should
be. The use of weight factors was proposed in \cite{Paltaufnew1,Paltaufnew}
in order to partially compensate for this distortion. In particular, in
\cite{Paltaufnew} smooth weight factors (depending on a reconstruction
point) are assigned to each detector in such a way that the total weight for
each direction is exactly one. This method is not exact; the error is
described by a certain smoothing operator. However, the singularities (or
jumps) in the image will be reconstructed correctly. As shown by numerical
examples in \cite{Paltaufnew}, such a correction visually significantly
improves the reconstruction. Moreover, iterative refinement is proposed in
\cite{Paltaufnew1,Paltaufnew} to further improve the image, and it is shown
to work well in numerical experiments.

Returning to non-iterative techniques, one should mention an interesting
attempt made in \cite{Sarah,Patch}) to generate the missing data using the
moment range conditions for $\M$ (see Section \ref{S:range}). The resulting
algorithm, however, does not seem to recover the values well; although, as
expected, it reconstructs all visible singularities.

An accurate $2D$ non-iterative algorithm for reconstruction from data measured
on an open curve $S$ was proposed in \cite{Kun_open}. It is based on
pre-computing approximations of plane waves in the region of interest $\Omega$
by the single layer potentials of the form
$$
\int\limits_{S}Z(\lambda|y-x|)\rho(y)dl(y),
$$
where $\rho(y)$ is the density of the potential, which needs to be chosen
appropriately, $dl(y)$ is the standard arc length, and $Z(t)$ is either the
Bessel function $J_{0}(t)$, or the Neumann function\ $Y_{0}(t)$. Namely, for
a fixed $\xi$ one finds numerically the densities $\rho_{\xi,J}(y)$ and
$\rho_{\xi,Y}(y)$ of the potentials
\begin{align}
W_{J}(x,\rho_{\xi,J})  &  =\int_{S}J_{0}(\lambda|y-x|)\rho_{\xi,J}
(y)dl(y),\label{potj}\\
W_{Y}(x,\rho_{\xi,Y})  &  =\int_{S}Y_{0}(\lambda|y-x|)\rho_{\xi,Y}(y)dl(y),
\label{poty}
\end{align}
where $\lambda=|\xi|,$ such that
\begin{equation}
W_{J}(x,\rho_{\xi,J})+W_{Y}(x,\rho_{\xi,Y})\thickapprox\exp(-i\xi\cdot
x)\mbox{ for all } x\in\Omega. \label{approx}
\end{equation}
Obtaining such approximations is not trivial. One can show that exact equality
in (\ref{approx}) cannot be achieved, due to different behavior at infinity of
the plane wave and the approximating single-layer potentials. However, as
shown by numerical examples in \cite{Kun_open}, if each point in $\Omega$ is
``visible'' from $S$, very accurate \emph{approximations} can be obtained,
while keeping the densities $\rho_{\xi,J}$ and $\rho_{\xi,Y}$ under certain
control.

Once the densities $\rho_{\xi,J}$ and $\rho_{\xi,Y}$ have
been found for all $\xi$, function $f(x)$ can be easily reconstructed. Indeed,
for the Fourier transform $\hat{f}(\xi)$ of $f(x)$
\begin{equation}
\hat{f}(\xi)=\frac{1}{2\pi}\int_{\Omega}f(x)\exp(-i\xi\cdot x)dx,\nonumber
\end{equation}
 one obtains, using (\ref{approx})
\begin{align}
\hat{f}(\xi)  &  \thickapprox\frac{1}{2\pi}\int_{\Omega}f(x)\left[
W_{J}(x,\rho_{\xi,J})+W_{Y}(x,\rho_{\xi,Y})\right]  dx\nonumber\\
&  =\frac{1}{2\pi}\int_{S}\left[  \int_{\Omega}f(x)J_{0}(\lambda
|y-x|)dx\right]  \rho_{\xi,J}(y)dl(y)\nonumber\\
&  +\frac{1}{2\pi}\int_{S}\left[  \int_{\Omega}f(x)Y_{0}(\lambda
|y-x|)dx\right]  \rho_{\xi,Y}(y)dl(y), \label{planewavemethod}
\end{align}
where the inner integrals are computed from the data $g$:
\begin{align}
\int_{\Omega}f(x)J_{0}(\lambda|y-x|)dx  &  =\int_{R^{+}}g(y,r)J_{0}(\lambda
r)dr,\label{gj}\\
\int_{\Omega}f(x)Y_{0}(\lambda|y-x|)dx  &  =\int_{R^{+}}g(y,r)Y_{0}(\lambda
r)dr. \label{gy}
\end{align}

Formula (\ref{planewavemethod}), in combination with (\ref{gj}) and (\ref{gy}),
yields values of $\hat{f}(\xi)$ for arbitrary $\xi$. Now
$f(x)$ can be recovered by numerically inverting the Fourier transform, or by
a reduction to a FBP inversion \cite{Kak,Natt_old} of the regular Radon
transform.

The most computationally expensive part of the algorithm, which is computing
the densities $\rho_{\xi,J}$ and $\rho_{\xi,Y}$, needs to be done only once
for a given acquisition surface. Thus, for a scanner with a fixed $S$, the
resulting densities can be pre-computed once and for all. The actual
reconstruction part then becomes extremely fast.

\begin{figure}[th]
\begin{center}
\begin{tabular}
[c]{ccc}
\includegraphics[width=1.5in,height=1.5in]{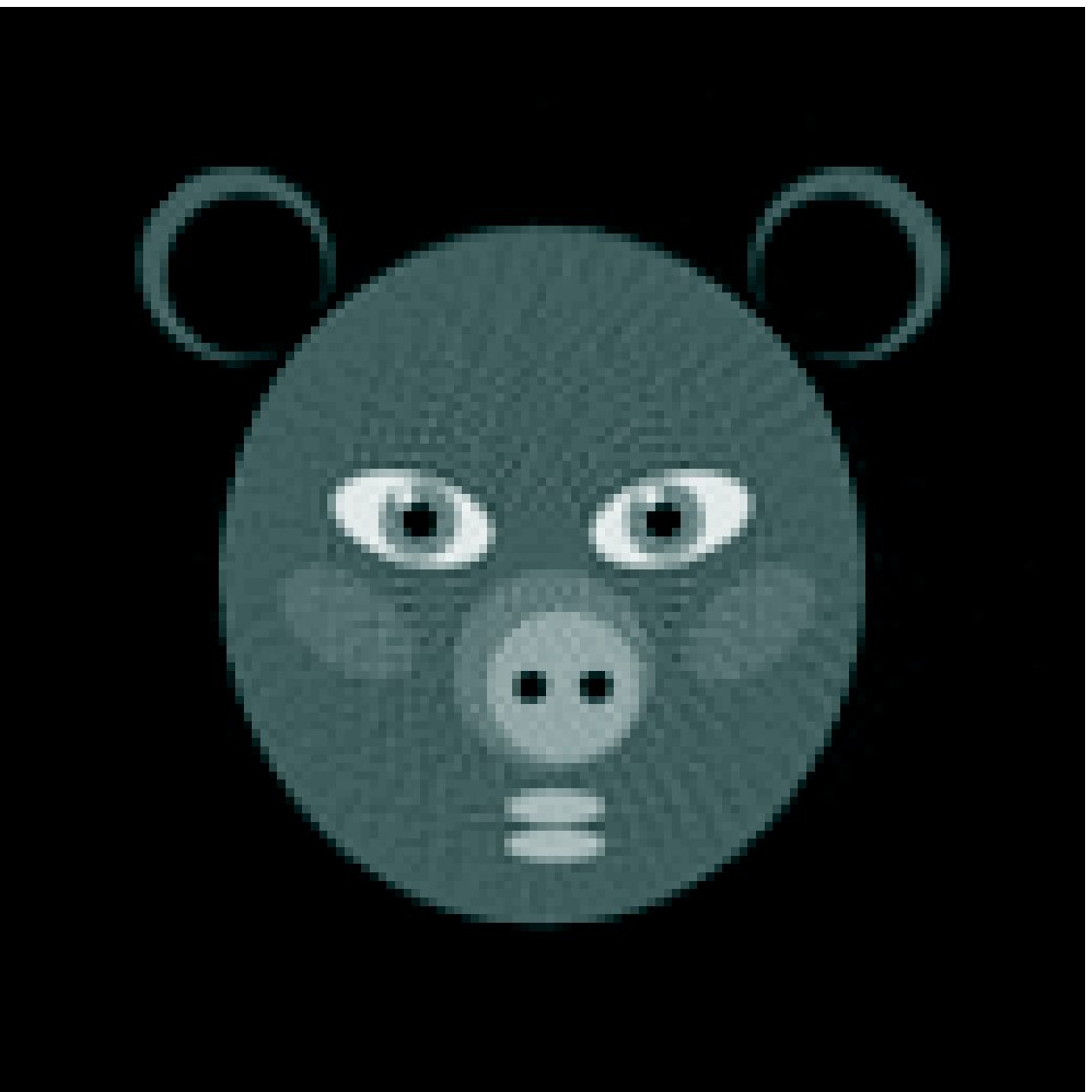} &
\includegraphics[width=1.5in,height=1.5in]{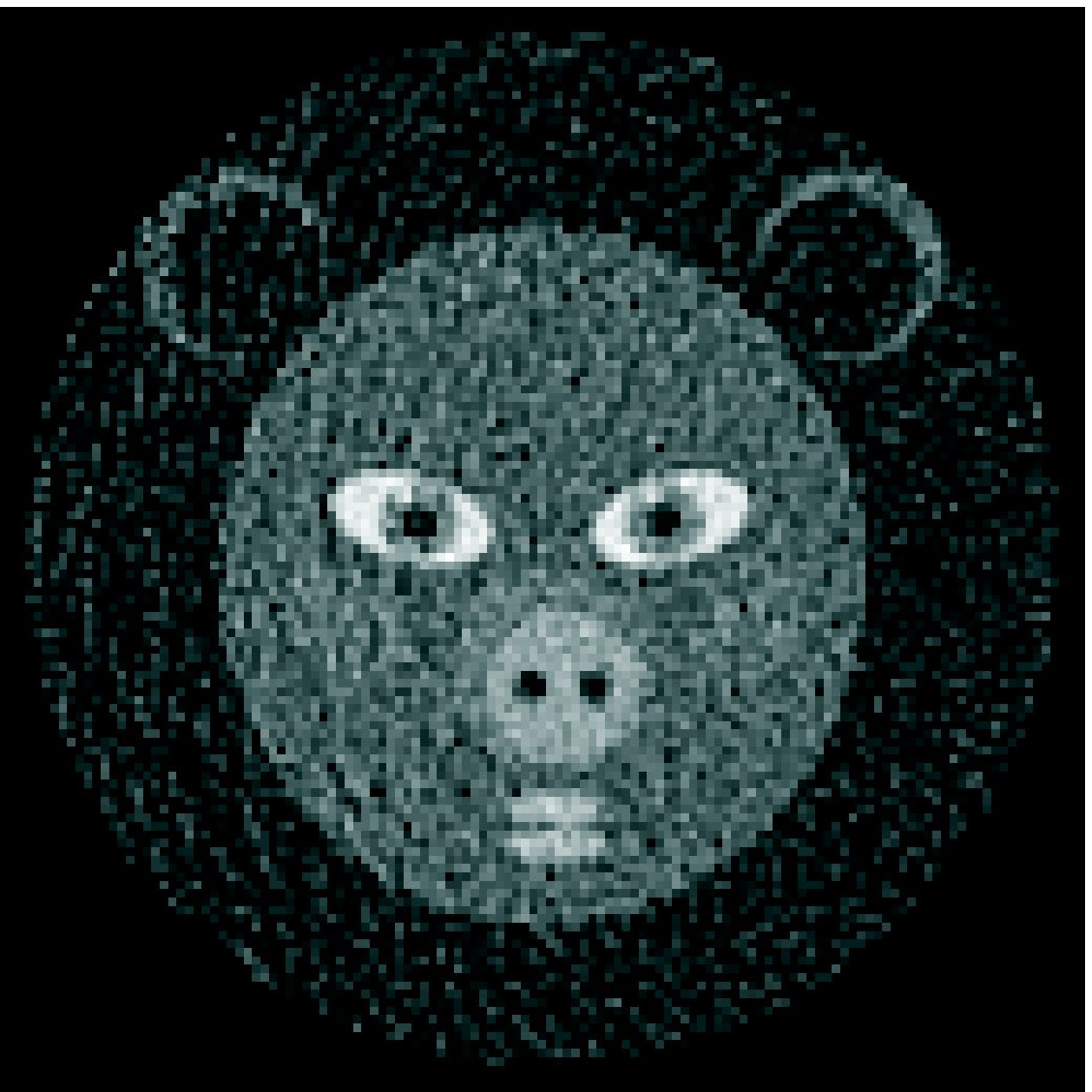} &
\includegraphics[width=1.5in,height=1.5in]{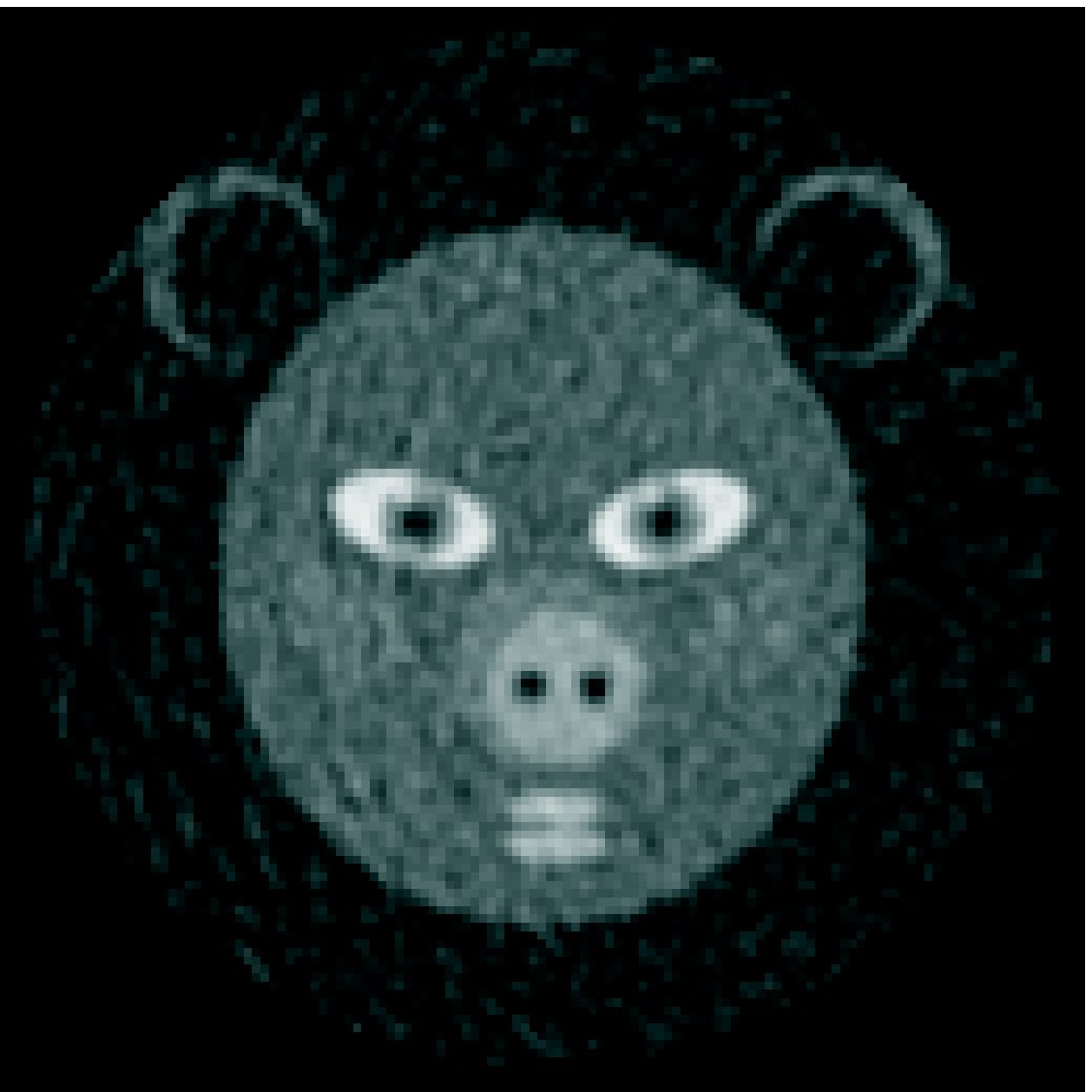}\\
(a) & (b) & (c)\\
&  &
\end{tabular}
\end{center}
\caption{Examples of reconstruction from incomplete data using
the technique of \cite{Kun_open}. Detectors are
located on the part of circular arc of radius 1.3 lying left of the line
$x_1=1$.  (a) reconstruction from accurate data
(b) reconstruction from the data with added 15\% noise
(c) reconstruction from noisy data with additional smoothing filter}
\label{F:incomplete}
\end{figure}

Examples of reconstructions from incomplete data using
this technique of \cite{Kun_open}) are shown in
Figure \ref{F:incomplete}. The images were reconstructed
within the unit square $[-1,1]\times[-1,1]$, while the detectors were placed
on the part of the concentric circle of radius 1.3 lying to the left of line
$x_1=1$. We used the same phantom as in Figure \ref{F:2d}(a));
the reconstruction from the data with added 15\% noise is shown in
Figure \ref{F:incomplete}(b); part (c) demonstrates the results of
applying additional smoothing filter to reduce the effects of noise
in the data.

\subsubsection{Variable speed of sound}\label{SS:num_var}

The problem of numerical reconstruction in TAT\ from the data measured on open
surfaces in the presence of a known variable speed of sound currently remains largely
open. One of the difficulties was discussed in Section \ref{S:incomplete}:
even if the speed of sound $c(x)$ is non-trapping, it can  happen that some of
the characteristics escape from the region of interest to infinity without
intersecting the open measuring surface. Then stable reconstruction of the
corresponding interfaces will become impossible. It should be possible,
however, to develop stable reconstruction algorithms in the case when the
whole object of interest is located in the visible zone.

The generalization of the method of \cite{Kun_open} to the case of variable
speed of sound is so far problematic, since this algorithm is based on the
knowledge of the open space Green's function for the Helmholtz equation. In
the case of a non-constant $c(x)$, this Green's function is position-depended,
and its numerical computation is likely to be prohibitively time-consuming.

A promising approach to this problem, currently under development, is to use time
reversal with the missing data replaced by zeros, or maybe by a more clever
extension (e.g., using the range conditions, as in \cite{Patch,Sarah}). This
would produce an initial approximation to $f(x)$, which one can try to refine
by fixed-point iterations; however, the pertinent questions concerning such an
algorithm remain open.

An interesting technique of using a reverberant cavity enclosing the target
to compensate for the missing data is described in \cite{Cox2:2007}.

\section{Final remarks and open problems}\label{S:remarks}
We list here some unresolved issues of mathematics of TAT/PAT, as well as
some developments that were not addressed in the main text.

\begin{enumerate}
\item The issue of uniqueness acquisition sets $S$ (i.e., such that
transducers distributed along $S$ provide sufficient information for TAT
reconstruction) can be considered to be resolved, for most practical
purposes. However, there remain significant unresolved theoretical
questions.  One of them consists of
proving an analog of Theorem \ref{T:AQ} for non-compactly
supported functions with a sufficiently fast (e.g., super-exponential) decay
at infinity. The original (and only known) proof of this
theorem uses microlocal techniques \cite{StefUhl,AQ} that significantly rely
upon the compactness of support. However, one hopes that
the condition of a fast decay should suffice for this result.
In particular, there is no proven analog of
Theorem \ref{T:ABK} for non-closed sets $S$ (unless $S$ is an open part of a
closed analytic surface).

Techniques developed in \cite{FPR} (see also \cite{AmbKuc_inj} for their
further use in TAT) might provide the right approach.

This also relates to still unresolved situation in dimensions $3$ and
higher. Namely, one would like to prove Conjecture \ref{C:n-dim}.

\item Concerning the inversion methods, one notices that closed form
formulas are known only for spherical, cylindrical, and planar acquisition
surfaces. The question arises whether closed form inversion formulas could be
found for any other closed surface? It is the belief of the authors that the
answer to this question is negative.

Another feature of the known closed form formulas that was mentioned before
is that they do not work correctly if the support of the sought function
$f(x)$ lies partially outside the acquisition surface.
Time reversal and eigenfunction
expansion methods do not suffer from this deficiency. The question arises
whether one could find closed form formulas that reconstruct the function
inside $S$ correctly, in spite of it having part of
its support outside. Again, the authors believe that the answer is negative.

\item Besides algebraic iterative approaches, there are no reliable reconstruction methods in the case of the detectors partially surrounding the target, if the medium is acoustically inhomogeneous (see Section \ref{SS:num_var}). This contrasts with the acoustically homogeneous situation (Section \ref{S:openalgs_const}).

\item The complete range description of the forward operator $\W$ in even
dimensions is still not known. It is also not clear whether one can obtain
complete range descriptions for non-spherical observation sets $S$ or for a
variable sound speed. The moment and orthogonality conditions do hold in the
case of a constant speed and arbitrary closed surface, but they do not
provide a complete description of the range. For acoustically inhomogeneous
media, an analog of orthogonality conditions exists, but it also does not
describe the range completely.

\item The problem of unique determination of the speed of sound from TAT
data is largely open.

\item As it was explained in the text, knowing full Cauchy data of the
pressure $p$ (i.e., its value and the value of its the normal derivative) on
the observation surface $S$ leads to unique determination and simple
reconstruction of $f$. However, the normal derivative is not measured by
transducers and thus needs to be either found mathematically or measured
in a different experiment. Thus, feasibility of techniques \cite{Klibanov,AmmariPAT}
relying on full Cauchy data  requires further mathematical and
experimental study.

\item In the standard X-ray CT, as well as in SPECT, the so called
\textbf{local tomography} technique \cite{FRS,FRS1,KLM,FFRS} is often very
useful. It allows one to emphasize in a stable way singularities (e.g.,
tissue interfaces) of the reconstruction, even in the case of incomplete
data (in the latter case, the invisible parts will be lost). An analog of
local tomography can be easily implemented in TAT, for instance, by
introducing an additional high-frequency filter in the FBP type formulas.

\item The mathematical analysis of TAT presented in the text did not take
into account the issue of modeling and compensating for the acoustic
attenuation. This subject is addressed in
\cite{Maslov,Anastasio_atten,Burgh_atten,Patch_atten,Kowar}, but probably
cannot be considered completely resolved.

\item The initial pressure $f(x)$ that was the center of all discussions in
the chapter (as well as in most papers devoted to TAT/PAT), is related, but
is not exactly identical to the optical features of interest of the tissue.
The issue of recovering the actual optical parameters of the tissue after
the initial pressure $f(x)$ is found is non-trivial and is addressed,
probably for the first time, in \cite{Bal}.

\item This chapter as well as most other papers devoted to TAT/PAT is
centered on the initial pressure $f(x)$. This quantity is related, but
is not exactly identical to the relevant optical features of the tissue.
The problem of recovering the actual optical parameters of tissue
(after $f(x)$ is found) is non-trivial and is addressed,
probably for the first time, in \cite{Bal}.

\item The TAT technique discussed in the chapter uses active interrogation
of the medium. There is a discussion in the literature of a passive version
of TAT, where no irradiation of the target is involved \cite{passive}.

\end{enumerate}

\section*{Acknowledgments}
The work of both authors was partially supported by the NSF DMS grant
0908208. The first author was also supported by the NSF DMS grant 0604778
and by the KAUST grant KUS-CI-016-04 through the IAMCS. The work of the
second author was partially supported by the DOE grant DE-FG02-03ER25577.
The authors express their gratitude to NSF, DOE, KAUST, and IAMCS for the
support.



\begin{thebibliography}{999}
\bibitem{FR2} D. Finch and  Rakesh. Recovering a function from its
spherical mean values in two and three dimensions. In \cite{CRC}, 77-88.
\end{thebibliography}
\end{document}